\documentclass[12pt]{article}
\usepackage{graphicx} % Required for inserting images

\usepackage{amsthm,amsfonts,amssymb,amsmath,exscale,relsize,derivative,amsthm}
\usepackage{natbib}
\usepackage[colorlinks,citecolor=blue,urlcolor=blue,filecolor=blue,backref=page]{hyperref}
\usepackage{graphicx}

\usepackage{tocbibind}
\usepackage[toc,page]{appendix}
\providecommand{\U}[1]{\protect\rule{.1in}{.1in}}
\newtheorem{theorem}{Theorem}

\newtheorem{lemma}{Lemma}
\newtheorem{corollary}{Corollary}

\newtheorem{remark}{Remark}
\newtheorem{solution*}{Solution}

\DeclareMathOperator*{\argmax}{argmax}

\def\U{\mathcal{U}}

\def\1{\mathbb{1}}

\newcommand\simiid{\stackrel{iid}{\sim}}
\newcommand\simind{\stackrel{ind}{\sim}}

\title{Posterior Contraction Rate and Asymptotic Bayes Optimality for One Group  Global-Local Shrinkage Priors in Sparse Normal Means Problem}

\author{Sayantan Paul and Arijit Chakrabarti}
\date{}

\begin{document}

\maketitle

%% *** Frontmatter *** 

%%%%%%%%%%%%%%%%%%%%%%% file template.tex %%%%%%%%%%%%%%%%%%%%%%%%%
%
% This is a template file for the LaTeX package SVJour2 for the
% Springer journal "Annals of the Institute of Statistical Mathematics"
%
%                                    Springer Heidelberg 2006/07/26
%
% Copy it to a new file with a new name and use it as the basis
% for your article. Delete % as needed.
%
%%%%%%%%%%%%%%%%%%%%%%%%%%%%%%%%%%%%%%%%%%%%%%%%%%%%%%%%%%%%%%%%%%%
%
% First comes an example EPS file -- just ignore it and
% proceed on the \documentclass line
% your LaTeX will extract the file if required

%
% \usepackage{times}
% \usepackage[LY1]{fontenc}
% \usepackage[LY1,mtbold]{mathtime}
% %
% \usepackage{marvosym}
% \newcommand{\envelope}{(\raisebox{-.5pt}{\scalebox{1.45}{\Letter}}\kern-1.7pt)}
% %
% %\usepackage{mathptmx}
% \usepackage{times}
% \usepackage[LY1]{fontenc}
% \usepackage[LY1,mtbold]{mathtime}
% \usepackage{biograph}
% %
% \input{glov3.doi}
% \let\SprJLogo=\relax
% \DOIyear{2002}
% \OFyear{2003}
% \DOImsnr{6789}
% \journalnumber{10463}
% \idline{12: 345--678}{345}
% \papertype{generic article}
%
% \usepackage{mathptmx}      % use Times fonts if available on your TeX system
%
% insert here the call for the packages your document requires
%\usepackage{latexsym}
% etc.
%
% please place your own definitions here and don't use \def but
% \newcommand{}{}
%

% \documentclass[ba,preprint]{imsart}% use this for supplement article

\begin{abstract}

We investigate some asymptotic optimality properties of inference on the mean vector in the well-known normal means model under sparsity, using broad classes of one-group global-local shrinkage priors which include the horseshoe. We assume that the level of sparsity is unknown. We consider inference when the global shrinkage parameter is estimated in empirical Bayesian way and also for a full Bayes approach. \\
We show that the empirical Bayes posterior distributions (using the simple data-based estimator (\cite{van2014horseshoe}) and MMLE  (\cite{van2017adaptive}) of the global shrinkage parameter) of the mean vector contract around the true mean vector at a near minimax rate with respect to the squared $L_2$ loss. We prove similar results in the full Bayes approach. The results for the first empirical Bayes estimate hold for the very general class of priors considered earlier in \cite{ghosh2017asymptotic} and those corresponding to the MMLE and full Bayes approaches hold for the large sub-class of three parameter Beta normal priors. These results significantly generalize those of \cite{van2014horseshoe} and \cite{van2017adaptive} proved for the horseshoe prior. \\
We next consider simultaneous testing for the mean vector to infer which coordinates are actually non-zero. We study asymptotic optimality of a testing rule based on a full Bayes approach using the class of priors considered in \cite{ghosh2017asymptotic} when the prior on the global shrinkage parameter is allowed to be chosen arbitrarily on a suitable support. We work under the asymptotic framework of \cite{bogdan2011asymptotic} where the loss in using a multiple testing rule is defined as the total number of misclassifications, and it is assumed that the true mean is being generated from a two-groups spike and slab prior. We prove, in a first result of its kind in this problem, that the Bayes risk of our full-Bayes decision rule asymptotically attains the risk of the Bayes Oracle defined in \cite{bogdan2011asymptotic}, up to some constant. This reinforces the argument that one-group shrinkage priors can be very reasonable alternatives to two-groups priors for inference in sparse problems. \\
	
\end{abstract}

%% ** Keywords **

%% ** Mainmatter **

%\section{}\label{}

% \begin{figure} 
% \includegraphics{<eps-file>}% place <eps-file> in ./img  subfolder
% \caption{}
% \label{}
% \end{figure}

% \begin{table} 
% *****************
% \begin{tabular}{lll}
% \end{tabular}
% *****************
% \caption{}
% \label{}
% \end{figure}

%%%%%%%%%%%%%%%%%%%%%%%%%%%%%%%%%%%%%%%%%%%%%%
%% Supplementary Material, if any, should   %%
%% be provided in {supplement} environment  %%
%% with title and short description.        %%
%%%%%%%%%%%%%%%%%%%%%%%%%%%%%%%%%%%%%%%%%%%%%%
%\begin{supplement}
%\stitle{???}
%\sdescription{???.}
%\end{supplement}

%% ** The bibliograhy **
%\bibliographystyle{ba}
%\bibliography{<bib-data-file>}% place <bib-data-file> in ./bib folder 

% ** Acknowledgements **
%\begin{acks}[Acknowledgments]
%\end{acks}

\section{Introduction} 
Analysis and inference on data with a large number of variables have become problems of great importance for researchers in the last few decades. Such data regularly come from various fields like genomics, finance, astronomy and economics, just to name a few. We are interested in the high-dimensional asymptotic setting when the total number of parameters in the problem grows at least at the same rate as the number of observations. We study in this paper the well-known normal means model described next. We have an observation vector $\mathbf{X} \in \mathbb{R}^{n}$, $\mathbf{X}=(X_1,X_2,\cdots, X_n)$, where each observation $X_i$ is expressed as,
\begin{equation*} \label{eq:1.1}
    X_i=\theta_i+\epsilon_i,i=1,2,\cdots,n, \tag{1}
\end{equation*}
where the parameters $\theta_1,\theta_2,\cdots,\theta_n$ are unknown fixed quantities
and $\epsilon_i$'s are independent $\mathcal{N}(0,1)$ random variables. Our goal is to infer about the mean vector $\boldsymbol {\theta} = (\theta_1,\ldots,\theta_n)$ as $n \rightarrow \infty$. This model has been widely studied in the literature, see, e.g.,  \cite{johnstone2004needles}, \cite{efron2004large} and \cite{jiang2009general}, and arises naturally in many potential applications including astronomical and other signal and image processing, microarray data analysis and nonparametric function estimation. Our focus in this paper is inference on $\boldsymbol {\theta}$ when it is sparse in the sense that the fraction of its non-zero coordinates is negligible compared to $n$ when $n$ is large.\\
\hspace*{0.5cm} In the Bayesian setting, a natural way to model a sparse $\boldsymbol{\theta}$ is by a two-groups prior like a spike and slab prior distribution (see, e.g., \cite{mitchell1988bayesian}), which is a mixture of a distribution degenerate at $0$ and a non-degenerate absolutely continuous  distribution $F$ over $\mathbb{R}$, given as
\begin{equation*} \label{eq:1.2}
    \theta_i \simiid{} (1-p) \delta_{0}+p F, i=1,2,\cdots,n, \tag{2}
\end{equation*}
where $p$ is small and denotes the probability that the true value of $\theta_i$ is non-zero. Analyzing this model
can be computationally prohibitive, especially in
high-dimensional problems and complex parametric frameworks. An alternative Bayesian approach where one uses a continuous one-group prior for modeling the parameters has become very popular in sparse settings. Such priors require substantially reduced computational effort for posterior analysis than their two-groups counterparts. Some examples of one-group priors are the t-prior due to \cite{tipping2001sparse}, the Laplace prior of \cite{park2008bayesian}, the normal-exponential-gamma prior due to \cite{brown2010inference}, the horseshoe prior of \cite{carvalho2009handling}, \cite{carvalho2010horseshoe}, the three-parameter beta normal prior of \cite{armagan2011generalized}, the generalized double Pareto prior introduced by \cite{armagan2013generalized} and the Dirichlet-Laplace prior of \cite{bhattacharya2015dirichlet}. All these priors can be expressed as “global-local" scale mixtures of normals as
\begin{equation*} \label{eq:1.3}
    \theta_i|\lambda_i,\tau \simind{} \mathcal{N}(0,\lambda^2_i\tau^2), \lambda^2_i \simind{} \pi_1(\lambda^2_i), \tau \sim \pi_2(\tau). \tag{3}
\end{equation*}
where $\lambda_i$, $i=1,\ldots,n$ are the \textit{local shrinkage} parameters used for detecting the true signals and $\tau$ is the \textit{global shrinkage} parameter used to control the overall sparsity in the model. \cite{polson2010shrink} recommended that $\pi_1$ should be heavy-tailed in order to accommodate large signals and $\pi_2$ should give large probability around zero to induce overall sparsity.\\
%As recommended by \cite{polson2010shrink}, in sparse settings, the prior $\pi_1(.)$ should have thick tails while $\pi_2(.)$ should have substantial mass near zero. In this way, one achieves the dual objectives of inducing an overall sparsity while maintaining the ability to accommodate large signals in the model. The horseshoe prior is one example of such a prior. 
\hspace*{0.3cm} We are interested in answering certain questions about asymptotic optimality of inference on $\boldsymbol{\theta}$ using some general classes of one-group priors. Such questions have been investigated in some depth in the literature when the global shrinkage parameter is used as a tuning parameter, i.e., in \eqref{eq:1.3}, $\pi_2(\cdot)$ is taken to be a degenerate probability. However, some interesting questions are still unanswered when one does an empirical Bayesian treatment thereof or a full Bayesian analysis using a non-degenerate prior. Our work is a modest attempt to answer some of these questions. In the next few paragraphs, we motivate and describe our work while touching upon the relevant literature. We have come up with totally new results for answering some of our questions, while the rest are generalizations/extensions of some existing works requiring novel technical arguments. \\
\hspace*{0.5cm} The systematic theoretical study of asymptotic optimality of inference using one-group priors started in \cite{datta2013asymptotic}. They considered the problem of simultaneous hypothesis testing to decide which coordinates of $\boldsymbol {\theta}$ are non-zero, under the asymptotic framework of \cite{bogdan2011asymptotic} (described in detail in Section \ref{sec-4}). As in \cite{bogdan2011asymptotic}, \cite{datta2013asymptotic} worked in the setting where the true means are generated from a two-groups prior of the form \eqref{eq:1.2} where the slab part is taken to be normally distributed with zero mean and a large variance. \cite{datta2013asymptotic} used the horseshoe prior to model the mean parameters in this context and theoretically studied a simultaneous testing rule based on the resulting posterior. The rule is formally described in Section \ref{sec-4} and was inspired by a rule originally suggested in \cite{carvalho2010horseshoe}. \cite{datta2013asymptotic} proved that the Bayes risk of the rule based on the horseshoe asymptotically attains, up to a constant, the risk of the Bayes Oracle, namely the Bayes rule in the two-groups setting (with respect to an additive misclassification loss). This was proved by choosing the global shrinkage parameter of the horseshoe of the same order as $p$, assuming it to be known. This suggests that carefully chosen one-group priors can be an alternative to the more complicated two-groups priors for inference in sparse problems, as was already observed in simulations in \cite{carvalho2010horseshoe}. In a different direction, \cite{van2014horseshoe} used the horseshoe prior in the problem of estimating $\boldsymbol{\theta}$ by choosing the global shrinkage parameter $\tau$ in terms of the proportion of non-zero means assuming it to be known. Their study was under the framework where the true mean vector $\boldsymbol{\theta}_0$ is assumed to be sparse in the ``nearly black" sense of \cite{donoho1992maximum}, that is, $\boldsymbol{\theta}_0 \in l_0[q_n]$ where $l_0[q_n]= \{\boldsymbol{\theta} \in \mathbb{R}^n: \# (1 \leq i \leq n: \theta_i \neq 0) \leq q_n \}$ with $q_n=o(n)$. They proved that by choosing $\tau=(\frac{q_n}{n})^{\alpha}$, for $\alpha \geq 1$, the risk of the  Bayes estimate ( with respect to the squared $L_2$ loss) corresponding to the horseshoe prior asymptotically attains the minimax risk in this problem \textit{up to a multiplicative constant}. It was proved further that the corresponding posterior distribution contracts around the true mean at the minimax rate. \cite{ghosh2016asymptotic} and \cite{ghosh2017asymptotic} substantially generalized the results of \cite{datta2013asymptotic} and \cite{van2014horseshoe} for broad classes of one-group priors including the horseshoe, but still choosing the parameter $\tau$ using knowledge of the level of sparsity.  It is evident from such results
that the choice of the \textit{global shrinkage} parameter $\tau$ based on the known level of sparsity plays a key role in ensuring optimality for both types of inferences.\\
\hspace*{0.5cm} The proportion of non-zero means may however be unknown in practice. To address this issue, in a first work of its kind, \cite{van2014horseshoe} proposed a simple estimator $\widehat{\tau}$ of $\tau$ (defined in equation \eqref{eq:2.4}) obtained using the observed data and plugging it in the formula for the Bayes estimate (given $\tau$) of the mean vector using horseshoe prior. A main motivation behind this choice of $\widehat{\tau}$ is one interpretation of $\tau$ as the fraction (upto a logarithmic factor) of non-zero co-ordinates of the true mean, see, e.g., \cite{van2017adaptive} in this context. Although the more natural and commonly used empirical Bayes estimator of $\tau$, namely the MMLE of \cite{van2017adaptive} (discussed in the next pargraph) is also available now, the estimator of \cite{van2014horseshoe} has its appeal because of its simplicity and interpretability and is worth studying. \cite{van2014horseshoe} showed that the horseshoe estimator used with their proposed estimator of $\tau$ attains a near-minimax rate for squared $L_2$ loss. The term ``near-minimax" is explained in detail in Section \ref{sec-2}. \cite{ghosh2017asymptotic} extended this result for the class of priors satisfying \eqref{eq:1.3} where the local shrinkage parameter is modeled as
\begin{equation*} \label{eq:1.4}
    \pi_1(\lambda^2_i)=K (\lambda^2_i)^{-a-1} L(\lambda^2_i), \tag{4}
\end{equation*}
where $K \in (0,\infty)$ is the constant of proportionality, $a$ is a positive real number and $L :(0,\infty) \to (0,\infty)$ is measurable non-constant slowly
varying function (defined in section \ref{sec-2}) satisfying \hyperlink{assumption1}{Assumption 1} of section \ref{sec-2}. This class of priors contains three parameter beta normal (which contains, e.g., the horseshoe, Strawderman-Berger prior and nonrmal-exponential-gamma priors), generalized double Pareto, and half t priors, just to name a few. In a follow-up paper, \cite{van2017adaptive} proved that the empirical Bayes posterior distribution for the horseshoe using the simple estimator of $\tau$ under discussion contracts around the truth at a near minimax rate. The natural question on whether this optimality holds for the broader class of priors satisfying \eqref{eq:1.3} and \eqref{eq:1.4} is still not addressed in the literature. We want to investigate this in this paper. Towards this, a key observation is Theorem 2 of \cite{ghosh2017asymptotic} which shows that
 the upper bound on the expected squared $L_2$ loss of the empirical Bayes estimate of the global-local priors is asymptotically of a near-minimax order. Given this, all that is needed to get the desired optimality result (for the general class of priors) is an upper bound of the correct order on the expected posterior variance corresponding to such priors. We are able to establish such an upper bound whereby we have our Theorem \ref{Thm-1} proving a near-minimax contraction rate for the empirical Bayes posterior corresponding to this general class of priors. Proof of this upper bound requires invoking several new arguments involving an important division of the support of $\widehat{\tau}$ and in particular, deriving a sharper upper bound (compared to that obtained in \cite{ghosh2017asymptotic}) on the posterior mean of the shrinkage coefficient (defined in Section \ref{sec-2}) as in Lemma \ref{lem3}. In addition, adopting a technique of \cite{van2014horseshoe} for proving an important upper tail probability bound for the proposed estimator also proved useful.\\ %It is explained in detail in Section \ref{sec-3} how the techniques of this paper for proving this result differ from those of \cite{ghosh2017asymptotic}.  
\hspace*{0.5cm} As mentioned earlier, \cite{van2017adaptive} introduced and studied a different estimator of $\tau$, namely the maximizer (with respect to $\tau \in [\frac{1}{n},1]$) of the marginal likelihood function of $\mathbf{X}$ given $\tau$ (defined in section \ref{sec-3.2}) in case of horseshoe prior. This estimator is referred to as the Maximum Marginal Likelihood Estimator(MMLE) of $\tau$, where the maximization over $\tau \in [\frac{1}{n},1]$ ensures that the estimator does not collapse to zero. MMLE seems to be a more natural estimator in the empirical Bayesian sense and may turn out to be a better estimator in terms of controlling the mean squared error (see the simulations of \cite{van2017adaptive}) if the true signals are not large enough. \cite{van2017adaptive} proved near-minimax posterior contraction rate of the empirical Bayes posterior around the truth by using the MMLE. Crucial to their proof was showing that the MMLE corresponding to the horseshoe prior lies in a certain small interval (depending on the level of sparsity) with probability tending to 1 (see condition (\hyperlink{C1}{C1}) in Section \ref{sec-3.2}). A natural question is whether a similar phenomenon holds for the MMLE when the mean parameter is modelled by a more general class of one group priors. Interestingly, we are able to show that this condition is indeed true when the mean vector is modelled by the three parameter beta normal priors. This can be proved by using interesting technical observations and arguments to adapt the basic line of the proof of \cite{van2017adaptive} to this more general context. As a result, it then follows easily (Theorem \ref{Thm1}) that a near-minimax contraction rate result also holds for the empirical Bayes posterior using the MMLE of $\tau$ for the above mentioned class of priors.
\\ 
\hspace*{0.5cm}
Another natural approach to model $\tau$ when the level of sparsity is unknown is the full Bayes one where $\tau$ is assumed to have a non-degenerate prior distribution. \cite{van2017adaptive} first considered such a full Bayes approach for inference using the horseshoe prior. They came up with two conditions on the prior density of $\tau$ (discussed in section \ref{sec-3.3}), under which the posterior distribution of $\boldsymbol{\theta} $ contracts around the true value at a near minimax rate. A key step in their proof (see their Lemma 3.6) is to show that the posterior distribution of $\tau$ asymptotically puts most of its mass in a region where the empirical Bayes estimates discussed before reside with high probability asymptotically. One of our goals is to investigate whether such results hold under the general class of priors mentioned above. We have established by proving Theorem \ref{Thm-3} that such optimality indeed holds for the three parameter beta normal class of priors, under conditions on the prior density of $\tau$ similar to those of \cite{van2017adaptive}. Towards that, we first write the full Bayes posterior mean as a weighted average of the posterior mean of $\boldsymbol{\theta}$ given $\tau$, the weight being the posterior distribution of $\tau$ given data. We then quantify the tail behaviour of the posterior distribution of $\tau$, which is a key ingredient of our proof as in the corresponding proof for \cite{van2017adaptive}. Towards this, the properties of the marginal log-likelihood given $\tau$ and its derivative (as in lemmas 1-6 in Appendix derived in the context of proving optimal contraction rate of empirical Bayes posterior distribution using MMLE) are very useful. The proof finally obtains by studying the tail behaviour of the posterior distribution of the mean vector given $\tau$, using techniques of \cite{ghosh2017asymptotic} and Theorem \ref{Thm-1}. \\
%More details are given in Section \ref{sec-3}. \\
\hspace*{0.5cm}  \cite{van2014horseshoe} commented that in order to attain optimal posterior contraction, a sharp peak near zero like the horseshoe prior might not be essential. \cite{ghosh2017asymptotic} provided support to this statement by proving such optimality results for a broad class of priors which do not necessarily have sharp peaks near zero, but are appropriately heavy-tailed and have sufficient mass around zero by appropriately choosing $\tau$. Theorem \ref{Thm-3} of this paper reinforces that intuition by establishing near-optimal posterior contraction rate for a very broad class of priors when one assumes any arbitrary prior on $\tau$ that is supported on an appropriately chosen interval. \\
%proves this to be true by proving optimal contraction results for a general class of priors of the form \eqref{eq:1.3} satisfying \eqref{eq:1.4} in the full Bayes setting which have a similar tail property and ensure sufficient mass around zero by appropriately choosing $\tau$.\\
 \hspace*{0.5cm} We then consider the problem of testing simultaneously whether each coordinate of $\boldsymbol{\theta}$ is zero or not based on the observations coming from \eqref{eq:1.1} in the sparse regime. Assuming that the true means are generated from a spike and slab prior described in Section \ref{sec-4}, \cite{bogdan2011asymptotic} proved that the Bayes risk of the optimal rule (with respect to an additive misclassification loss), namely the Bayes Oracle, for this problem is asymptotically of the form ${np\bigg[2 \Phi(\sqrt{{C}})-1 \bigg]} (1+o(1))$ where $n$ denotes the number of hypotheses, $p$ is the proportion of non-zero means and $C$ is a constant defined in \hyperlink{assumption2}{Assumption 2} in Section \ref{sec-4}. There have been several attempts in the literature to investigate asymptotic optimality of multiple testing rules based on one group global-local priors in this setting. The results in \cite{datta2013asymptotic}, \cite{ghosh2016asymptotic} and \cite{ghosh2017asymptotic}, used the global shrinkage parameter $\tau$ either as a tuning parameter or it was estimated in an empirical Bayesian way. \cite{ghosh2016asymptotic} reported a simulation study using one-group priors when a non-degenerate prior was placed on $\tau$. The results were promising, in that, for small values of $p$, the full Bayes rule practically matches the Bayes Oracle in terms of risk. Motivated by this, we consider the class of one-group priors satisfying \eqref{eq:1.3} and \eqref{eq:1.4} where the prior on $\tau$ is supported on a carefully chosen small interval. The details are given in ({\hyperlink{C4}{C4}}) of Section \ref{sec-4}. Our theoretical calculations reveal that condition  ({\hyperlink{C4}{C4}}) is sufficient for the Bayes risk of the induced decision rule (defined in \eqref{eq:4.9} in Section \ref{sec-4}) to attain the Optimal Bayes risk up to a multiplicative constant.  Moreover, the multiplicative constant obtained by us is smaller than those obtained by \cite{ghosh2016asymptotic} and \cite{datta2013asymptotic} for wide ranges of sparsity levels. We further show that for problems where the sparsity is pronounced, for ``horseshoe-type" priors, the decision rule attains the Bayes Oracle exactly asymptotically. Our results are adaptive and hold true for a wide range of sparsity levels. They also highlight the fact that one has the liberty to use any non-degenerate prior on $\tau$ supported on our proposed range for ensuring near-optimal results using one-group priors.
 To the best of our knowledge, these are the first optimality results of their kind in the literature in the full Bayes approach in this hypothesis testing problem. The proofs of our optimality results hinge upon careful choice on the range of the prior density of $\tau$ along with the use of monotonicity of the posterior expected value of the shrinkage coefficient with respect to $\tau$.\\ %They substantially reinforce the intuition that suitably chosen one-group priors can be meaningful alternatives to two-groups priors for sparse Bayesian modelling. 
\hspace*{0.5cm} The rest of the paper is organized as follows. Section \ref{sec-2} discusses the the main results on (near) optimal posterior concentration rates. In subsection \ref{sec-3.1}, we discuss the results
using the empirical Bayes estimate of $\tau$ as proposed by \cite{van2014horseshoe}. Subsection \ref{sec-3.2} also contains the results
using the empirical Bayes estimate of $\tau$, but due to \cite{van2017adaptive}. Results of the 
full Bayes approach are presented in subsection \ref{sec-3.3}. In Section \ref{sec-4}, optimality results in the problem of simultaneous hypothesis testing are discussed. Section \ref{sec-5} contains a discussion and proposes possible extensions of our work. %Proof of Theorem \ref{Thm-1} is in section \ref{sec-2}. Proofs of Theorems \ref{Thm-6}-\ref{Thm-8} are provided in section \ref{sec-4}. 
Proofs of most of the results can be found in the Appendix.

\subsection{Notations}
\label{sec-1.1}
For any two sequences $\{a_n\}$ and $\{b_n\}$ with with $\{b_n\} \neq 0$ for all large $n$, we write $a_n \asymp b_n$ to denote $0< \liminf_{n \to \infty} \frac{a_n}{b_n} \leq \limsup_{n \to \infty } \frac{a_n}{b_n} < \infty$ and $a_n \lesssim b_n$ to denote for sufficiently large $n$, there exists a constant $c>0$ such that $a_n \leq cb_n$. Finally, $a_n=o(b_n)$ denotes $\lim_{n \to \infty} \frac{a_n}{b_n}=0$. 
$\phi(\cdot)$ denotes the density of a standard normal distribution. In the full Bayes procedure, we use the notation $\tau_n(q_n)$ which denotes $(q_n/n)\sqrt{\log (n/(q_n))}$.

\section{Contraction rate of posterior distribution using global-local shrinkage priors}
\label{sec-2}
In this section, we discuss our results on the posterior contraction rate for certain classes of global-local shrinkage priors we consider, when the level of sparsity is assumed unknown. The results of the empirical Bayes and the full Bayes approaches are presented in three separate subsections.\\
\hspace*{0.5cm} We start with some preliminary definitions and  assumptions will be used throughout. Recall that a function $L(\cdot): (0,\infty) \to (0,\infty)$ is said to be slowly varying, if for any $\alpha>0,\frac{L(\alpha x)}{L(x)} \to 1$  as $x \to \infty$. For the theoretical development of the paper, we consider slowly varying functions that satisfy \hyperlink{assumption1}{Assumption 1}
below.
\\
\textbf{\hypertarget{assumption1}{Assumption 1.}} \\
\textbf{(\hypertarget{A1}{A1})} There  exists  some $c_0(>0)$ such that $L(t) \geq c_0 \; \forall t \geq t_0$, for some $t_0>0$, which depends on both $L$ and $c_0$. \\
\textbf{(\hypertarget{A2}{A2})} There exists some $M \in (0,\infty)$ such that $\sup_{t \in (0,\infty)} L(t) \leq M$. \\
These assumptions are similar to those of \cite{ghosh2017asymptotic}. \\
\hspace*{0.5cm}
The minimax rate is considered to be a useful benchmark for the rate of contraction of posterior distributions. According to Theorem
2.5 of \cite{ghosal2000convergence}, posterior distributions cannot contract around the
truth faster than the minimax risk. Thus it is of considerable interest to investigate if the posterior distribution generated by a Bayesian procedure contracts around the truth at the minimax rate or not. When $\boldsymbol{\theta}_0 \in l_0[q_n]$,  \cite{donoho1992maximum} proved that the minimax risk under the usual squared $L_2$ loss $||  \boldsymbol{\tilde{\theta}} -\boldsymbol{\theta} ||^2=\sum_{i=1}^{n} (\tilde{\theta}_i-\theta_i)^2$ is asymptotically of the form
\begin{equation*}
    \inf_{\tilde{\theta}} \sup_{\boldsymbol{\theta}_0 \in l_0[q_n]} \mathbb{E}_{\boldsymbol{\theta}_0}||  \boldsymbol{\tilde{\theta}} -\boldsymbol{\theta} ||^2= 2q_n \log (\frac{n}{q_n}) (1+o(1)), 
\end{equation*}
when $q_n \to \infty $ as $ n \to \infty$ and $q_n=o(n)$. Many results on the posterior concentration rate available in the literature are in terms of $q_n \log n$. 
See,e.g., \cite{van2017adaptive} and some references therein. This is referred to as the ``near-minimax" rate in this literature. It may be however noted that the near-minimax rate $q_n \log n$ is equivalent to the minimax rate if $q_n \sim n^{\beta}$ for some $0<\beta<1$.
In both the empirical Bayes and full Bayes approaches, our results in the next subsections show that under the broad classes of priors, the posterior distributions of $\boldsymbol{\theta}$ contract around the truth at a near minimax rate, which substantially generalize the results of  \cite{van2014horseshoe} and \cite{van2017adaptive} proved for the horseshoe.
\subsection{Empirical Bayes Procedure-plug-in Estimator}
\label{sec-3.1}

When the proportion of non-zero means is unknown, \cite{van2014horseshoe} proposed using
an empirical Bayes estimate of $\tau$. The estimator is defined as
\begin{equation*} \label{eq:2.4}
	\widehat{\tau}=\text{max} \biggl\{\frac{1}{n}, \frac{1}{c_2n}\sum_{i=1}^{n}1 \{|X_i|>\sqrt{c_1\log{n}}\} \biggr\}, \tag{5}
\end{equation*}
where  $c_1 \geq 2$ and $c_2 \geq 1$  are constants. Observe that $ \widehat{\tau}$ always lies between $\frac{1}{n}$ and $1$. Let $T_{\widehat{\tau}}(\mathbf{X})$ be the Bayes estimate of $\boldsymbol{\theta}$ given $\tau$, namely $T_{\tau}(\mathbf{X})= E(\boldsymbol{\theta}|\mathbf{X},\tau)$ evaluated at $\tau=\widehat{\tau}$. \cite{van2014horseshoe} showed that the empirical Bayes estimate $T_{\widehat{\tau}}(\mathbf{X})$ corresponding to the horseshoe prior (by taking $a=0.5$ and $L(t)=t/(t+1)$ in \eqref{eq:1.4}) asymptotically attains the near minimax $L_2$ risk. This fact was generalized by \cite{ghosh2017asymptotic} who established the same asymptotic rate for the class of priors \eqref{eq:1.3} satisfying \eqref{eq:1.4} when $a \geq 0.5$  provided $q_n \propto n^{\beta},0<\beta<1$. \cite{van2017adaptive} also proved that the empirical Bayes posterior (i.e., the posterior distribution of the mean vector given $\tau$ evaluated at $\tau =\widehat{\tau}$) for the horseshoe contracts around the truth at a near minimax rate. Our first result, namely Theorem \ref{Thm-1} presented below, is a generalization of this result. The proof of this result is given in the Appendix.

\begin{theorem}
	\label{Thm-1}
	Suppose $\mathbf{X} \sim \mathcal{N}_n(\boldsymbol{\theta}_0,\mathbf{I}_n)$ where $\boldsymbol{\theta}_0 \in l_0[q_n]$ with $q_n \propto$  $n^{\beta}$ where $0<\beta<1$. Consider the class of priors
	\eqref{eq:1.3} satisfying \eqref{eq:1.4} with $a \geq \frac{1}{2}$ where $L(\cdot)$ satisfies \hyperlink{assumption1}{Assumption 1}. Then the empirical Bayes
	posterior contracts around the true parameter $\boldsymbol{\theta}_0$ at the near minimax rate, namely,  
	\begin{equation*}
	\sup_{\boldsymbol{\theta}_0 \in l_0[q_n]} \mathbb{E}_{\boldsymbol{\theta}_0} \Pi_{\widehat{\tau}} \bigg(  \boldsymbol{\theta}:||\boldsymbol{\theta}-\boldsymbol{\theta}_0 ||^2>    M_n q_n \log n |\mathbf{X} \bigg) \to 0 \hspace{0.1cm},
	\end{equation*}
	for any $M_n \to \infty$.
\end{theorem}

%\begin{proof}
%For obtaining the first part, we make use of Markov's inequality along with Theorem \ref{Thm-1} and Theorem 2 of \cite{ghosh2017asymptotic}. On the other hand,  the use of Markov's inequality together with Theorem \ref{Thm-1} is sufficient for proving the second part.
%\end{proof}
\begin{remark}
\label{rem-1}

   For proving Theorem \ref{Thm-1}, first we need to show that, the total posterior variance using the empirical Bayes estimator of $\tau$	corresponding to this class of priors satisfies, as $n \to \infty$,
	\begin{equation*} \label{eq:2.2.1}
	\sup_{\boldsymbol{\theta}_0 \in l_0[q_n]} \mathbb{E}_{\boldsymbol{\theta}_0} \sum_{i=1}^{n} Var(\theta_i|X_i,\widehat{\tau}) \lesssim q_n \log n, \tag{6}
	\end{equation*}
where $Var(\theta_i|X_i,\widehat{\tau})$ denotes $Var(\theta_i|X_i,{\tau})$ evaluated at $\tau=\widehat{\tau}$. In order to obtain \eqref{eq:2.2.1}, 
   we divide our calculations for the non-zero and zero means, namely in Step 1 and Step 2, respectively. For both cases, we need to invoke substantially new technical arguments to achieve our results.  For the non-zero means, using an interesting division of the range of $|X_i|$ (keeping in mind the minimax rate) and the monotonicity of the posterior mean of the shrinkage coefficients, we are able to exploit some known results of \cite{ghosh2017asymptotic} for the case when $\tau$ is used as a tuning parameter. For Step 2, we divide our calculations based on $a \in [\frac{1}{2},1)$ and $a \geq 1$. When $a \in [\frac{1}{2},1)$, using some novel handling of cases for different ranges of $|X_i|$ and $\widehat{\tau}$, we are again able to exploit the monotonicity of the expected posterior shrinkage coefficients and known results for the fixed $\tau$ cases. However, for $a \geq 1$, we need to substantially modify our arguments, since using similar arguments as before produces an upper bound on the posterior variance which exceeds the near minimax rate. To take care of this,
   We need to come up with a sharper upper bound on the posterior shrinkage coefficient compared to that in \cite{ghosh2017asymptotic} in Lemma \ref{lem3} presented below and employ a very handy upper bound on the posterior variance (given in equation (18) of the Appendix). 
%The novelty in proving Lemma \ref{lem3} lies in dividing the range of integration into two regions and treating the two regions separately for obtaining the upper bound. Further details are available in the remark following the proof of Lemma \ref{lem3}. 
\end{remark}
We end this subsection by stating and proving an important lemma related to the upper bound on the posterior mean of the shrinkage coefficient which is crucial for the proof of the theorem \ref{Thm-1}.
\begin{lemma}
	\label{lem3}
	Suppose $\mathbf{X} \sim \mathcal{N}_n(\boldsymbol{\theta}_0,\mathbf{I}_n)$. Consider the class of priors
	\eqref{eq:1.3} satisfying \eqref{eq:1.4} with $a \geq 1$. Then for any $\tau \in (0,1)$ and any $x_i \in \mathbb{R}$,
	\begin{align*} \label{L-1.1}
		E(1-\kappa_i|x_i,\tau)  & \leq \bigg(  \tau^2 e^{\frac{x^2_i}{4}}+ K \tau^2 \int_{1}^{\infty}  \frac{1}{({1+t \tau^2})^{\frac{3}{2}}} t^{-a} L(t) e^{\frac{x^2_i}{2}\cdot \frac{t \tau^2}{1+t \tau^2}} dt \bigg) (1+o(1)) , \tag{7}
	\end{align*}
	where the term $o(1)$ depends only on $\tau$ and $\lim_{\tau \to 0}o(1)=0$ and $\int_{0}^{\infty} t^{-a-1} L(t) dt=K^{-1}$.
	
\end{lemma}

\begin{proof}[Proof of Lemma \ref{lem3}]
	Posterior distribution of $\kappa_i$ given $x_i$ and $\tau$ is,
	\begin{equation*}
		\pi(\kappa_i|x_i,\tau) \propto \kappa_i^{a-\frac{1}{2}} (1-\kappa_i)^{-a-1}L(\frac{1}{\tau^2}(\frac{1}{\kappa_i}-1)) \exp{\frac{(1-\kappa)x^2_i}{2}} \hspace*{0.1cm}, 0<\kappa_i<1 .
	\end{equation*}
	Using the transformation $t=\frac{1}{\tau^2}(\frac{1}{\kappa}-1)$, $E(1-\kappa_i|x_i,\tau)$ becomes
	\begin{equation*}
		E(1-\kappa_i|x_i,\tau)= \frac{\tau^2 \int_{0}^{\infty} (1+t \tau^2)^{-\frac{3}{2}} t^{-a} L(t) e^{\frac{x^2_i}{2}\cdot \frac{t \tau^2}{1+t \tau^2}} dt}{\int_{0}^{\infty} (1+t \tau^2)^{-\frac{1}{2}} t^{-a-1} L(t) e^{\frac{x^2_i}{2}\cdot \frac{t \tau^2}{1+t \tau^2}} dt}\hspace*{0.1cm}.
	\end{equation*}
	Note that, 
	\begin{align*}
		\int_{0}^{\infty} (1+t \tau^2)^{-\frac{1}{2}} t^{-a-1} L(t) e^{\frac{x^2_i}{2}\cdot \frac{t \tau^2}{1+t \tau^2}} dt & \geq \int_{0}^{\infty} (1+t \tau^2)^{-\frac{1}{2}} t^{-a-1} L(t) dt= K^{-1}(1+o(1)),
	\end{align*}
 where $o(1)$ depends only on $\tau$ and tends to zero as $\tau \to 0$.
	The equality in the last line follows due to the Dominated Convergence Theorem. Hence
	\begin{equation} \label{L-1.2}
		E(1-\kappa_i|x_i,\tau) \leq K(A_1+A_2)(1+o(1)) \hspace*{0.1cm}, \tag{8}
	\end{equation}
	
	where
	\begin{align*}
		A_1 &= \int_{0}^{1} \frac{t \tau^2}{1+t \tau^2} \cdot \frac{1}{\sqrt{1+t \tau^2}} t^{-a-1} L(t) e^{\frac{x^2_i}{2}\cdot \frac{t \tau^2}{1+t \tau^2}} dt \hspace*{0.1cm}.
	\end{align*}
	and $A_2=  \int_{1}^{\infty} \frac{t \tau^2}{1+t \tau^2} \cdot \frac{1}{\sqrt{1+t \tau^2}} t^{-a-1} L(t) e^{\frac{x^2_i}{2}\cdot \frac{t \tau^2}{1+t \tau^2}} dt$.
	Since, for any $t \leq 1$ and $\tau \leq 1$, $\frac{t \tau^2}{1+ t \tau^2} \leq \frac{1}{2}$, using the fact $\int_{0}^{\infty} t^{-a-1} L(t) dt=K^{-1}$ (obtained from \eqref{eq:1.4}),
	\begin{equation*} \label{L-1.3}
		A_1 \leq  K^{-1} \tau^2 e^{\frac{x^2_i}{4}}\hspace*{0.1cm}. \tag{9}
	\end{equation*}
	
	Using the \eqref{L-1.2} and \eqref{L-1.3} along with the form of $A_2$, for any $\tau \in (0,1)$ and $x_i \in \mathbb{R}$, we get the result of the form \eqref{L-1.1}.
\end{proof}
\begin{remark}
	 Lemma \ref{lem3} is a refinement of Lemma 2 of \cite{ghosh2017asymptotic}. The proof is obtained by a careful inspection and division of the range of the integral in $(0,1)$ and $[1,\infty)$. Note that, the upper bound of $E(1-\kappa_i|x_i,\tau)$ within the interval $(0,1)$ is sharper than that of \cite{ghosh2017asymptotic}. This order is important for obtaining the optimal posterior contraction rate for the empirical Bayes treatment when $a \geq 1$ and hence is a key ingredient for proving the Theorem \ref{Thm-1}. In this proof,
the second term in the r.h.s. of \eqref{L-1.1} is either used unchanged or a further upper bound on it needs to be employed as may be seen in the proof of \textbf{Case-2} for Theorem \ref{Thm-1} (given in Appendix).
  \end{remark}
\subsection{Empirical Bayes Procedure-Maximum Marginal Likelihood Estimator(MMLE)}
\label{sec-3.2}
As an alternative to the simple estimator in \cite{van2014horseshoe}, \cite{van2017adaptive} proposed the Maximum Marginal Likelihood Estimator (MMLE) of $\tau$, defined as
\begin{align*} \label{eq:2.3}
    \widehat{\tau}_n &= \argmax_{\tau \in [\frac{1}{n},1]} \prod_{i=1}^{n} \phi(x_i-\theta_i) g_{\tau}(\theta_i) d \theta_i.
\end{align*}
Here $\phi(\cdot)$ denotes the density of a standard normal random variable and $g_{\tau}(\theta_i)$ is the marginal pdf of $\theta_i$ given $\tau$ defined as
\begin{equation*}
    g_{\tau}(\theta_i)= \int_{0}^{\infty} \frac{1}{\lambda_i\tau}\phi(\frac{\theta_i}{\lambda_i \tau}) \pi_1(\lambda^2_i) d \lambda^2_i.
\end{equation*}
\cite{van2017adaptive} studied MMLE when the unknown mean vector is modeled by the horseshoe prior. They proved that the empirical Bayes posterior distribution of the mean vector (with $\tau$ estimated by the MMLE) contracts around the truth at the near-minimax rate. We want to investigate if this result can be generalized when the mean vector is modeled by a broader class of one-group priors of which the horseshoe is a member. 
We recall at this point that in the literature of global-local priors used for sparse modelling, the parameter $\tau$ has to play in a sense the role of the fraction of non-zero parameters and is usually of the order of the level of sparsity or thereabouts. See, for e.g., \cite{datta2013asymptotic}, \cite{van2014horseshoe}, \cite{ghosh2017asymptotic} etc. This suggests that any empirical Bayes estimate that turns out to be close to such an optimal value with high probability might turns out to be also optimal. We also note that, in sparse situation the fraction tends to zero. So, it means the maximization can be restricted in an interval which contains the optimal value. { 
This motivates us to redefine the MMLE by the maximization w.r.t. $\tau$ in an interval $[\frac{1}{n}, \frac{1} {\log n}]$, i.e. 
\begin{align*} \label{eq:2.2}
    \widehat{\tau}_n &= \argmax_{\tau \in [\frac{1}{n},\frac{1} {\log n}]} \prod_{i=1}^{n} \phi(x_i-\theta_i) g_{\tau}(\theta_i) d \theta_i \tag{10},
\end{align*}
This also ensures that the optimal value of $\tau$ is of a smaller order than $\frac{1}{\log n}$.}
We are able to establish in our Theorem \ref{Thm1} (stated at the end of this section) that the empirical Bayes posterior distribution of the mean vector corresponding to the MMLE as given in \eqref{eq:2.2} contracts around the truth at the near-minimax rate for the class of three parameter beta normal priors (with $a =\frac{1}{2}$ and $\alpha \geq \frac{1}{2}$). This class of priors was introduced by \cite{armagan2011generalized} and is defined as
\begin{align*} \label{eq:2.7}
    \theta_i|\lambda^2_i,\tau^2 &\simind \mathcal{N}(0,\lambda^2_i\tau^2), \\
    \lambda^2_i &\simind \pi_1(\lambda^2_i), \tag{11} \\ 
    \tau & \sim \pi_2(\tau),
\end{align*}
where 
\begin{equation*}
    \pi_1(\lambda^2_i)=\frac{\Gamma(a+\alpha)}{\Gamma(a)\Gamma(\alpha)}  (\lambda^2_i)^{\alpha-1}  (1+\lambda^2_i)^{-(\alpha+a)} 
\end{equation*}
%{\tt put definition of three parameter beta normal here}
for $\alpha>0, a>0$.
\cite{ghosh2016asymptotic} has shown that the above class of priors can be represented as \eqref{eq:1.3} where $\pi_1(.)$ satisfies \eqref{eq:1.4} with $L(\lambda^2_i)=(1+\frac{1}{\lambda^2_i})^{-(a+\alpha)}$ with $a > 0$ and $\alpha > 0$ and $K=\frac{\Gamma(a+\alpha)}{\Gamma(a)\Gamma(\alpha)} \cdot$  Three parameter beta normal priors is a rich class of priors containing horseshoe
 for $a=\frac{1}{2}=\alpha$, the Strawderman–Berger prior with $\alpha=1, a=\frac{1}{2}$ and the
normal–exponential–gamma priors with $\alpha=1, a>0$ as some special cases.
 
 \hspace*{0.5cm} The proof of Theorem \ref{Thm1} crucially depends on proving that the MMLE based on the three parameter beta normal family satisfies condition (\hyperlink{C1}{C1}) below. This condition is exactly the same as that of \cite{van2017adaptive}, which was proved in that paper to be satisfied by the MMLE based on the horseshoe. Interestingly, we have observed while proving theorem \ref{Thm-1} that the simple estimator of \cite{van2014horseshoe} also satisfies a similar condition. We now state condition (\hyperlink{C1}{C1}).\\ \\
 \textbf{(\hypertarget{C1}{C1})} There exists a constant $C>0$ such that the empirical Bayes estimator of $\tau$, say $\widehat{\tau}_n \in [\frac{1}{n}, C\tau_n(q_n)]$, with $P_{\theta_0}$-probability tending to one, uniformly in $\theta_0 \in l_0[q_n]$. \\ \\
 Lemma \ref{lem1}, presented below, shows that the MMLE satisfies condition (\hyperlink{C1}{C1}). Proof of this lemma is given in the Appendix.
\begin{lemma}
    \label{lem1}
    The MMLE defined in \eqref{eq:2.2} satisfies (\hyperlink{C1}{C1}).
\end{lemma}
\begin{remark}
\label{rem-3}
    In order to prove Lemma \ref{lem1}, we need to obtain an upper bound of the derivative (with respect to $\tau$) of the logarithm of the marginal likelihood function of $\mathbf{X}$ given $\tau$ and figure out the region of $\tau$ in which the upper bound is non-negative. Roughly speaking, the upper bound is negative when $\frac{\tau}{\zeta_{\tau}} \gtrsim \frac{q_n}{(n-q_n)}$. As a result, for maximizing the likelihood (or the log-likelihood) $\widehat{\tau}_n$ has to satisfy $\frac{\widehat{\tau}_n}{\zeta_{\widehat{\tau}_n}} \lesssim \frac{q_n}{(n-q_n)}$, i.e., $\widehat{\tau}_n \lesssim \tau_n(q_n) $, under the assumption $q_n=o(n)$. We have proved in a series of lemmas(lemmas 1-6 of the Appendix,  similar to those of \cite{van2017adaptive}) results related to the upper bound of the derivative (with respect to $\tau$) of the logarithm of the marginal likelihood function of $\mathbf{X}$ given $\tau$. Proofs of all of these lemmas are given in the Appendix. The proof starts by showing that
     the derivative can be expressed as $\frac{1}{\tau}\sum_{i=1}^{n}m_{\tau}(x_i)$ 
     of suitably defined function $m_{\tau}(x_i)$. The next step is to decompose the sum based on $x_i$ coming from zero and non-zero means and to study their behaviours separately. We have proved that the average of $m_{\tau}(x_i)$ terms
     for zero means behaves like its expectation, uniformly in $\tau$. On the other hand, the function $m_{\tau}(x_i)$ is uniformly bounded by a constant, which is used in the case of non-zero means. All of these arguments are essential for proving Lemma \ref{lem1}. \cite{van2017adaptive} had shown that all of these arguments hold for horseshoe prior and stated a condition exactly as same as (\hyperlink{C1}{C1}). However, our calculations show that the same holds even for the three parameter beta normal class of priors containing horseshoe as a special case. In spite of the lemmas of \cite{van2017adaptive} being available for the horseshoe prior, it is a non-trivial task to extend their results to the three parameter beta normal class of priors due to the more general form of the marginal density of $\mathbf{X}$ given $\tau$. In our case the $m_{\tau}(x_i)$ function is significantly different from that of \cite{van2017adaptive} and it requires new technical observations and manipulations to obtain results similar to those of \cite{van2017adaptive} in this context. In this way, Lemma \ref{lem1} generalizes the Theorem 3.1 of \cite{van2017adaptive}.
\end{remark}
%After the lemma is established, now we are in a position to state the result related to the optimal posterior concentration rate of the empirical Bayes posterior around the truth when the mean parameter is modeled by the three parameter beta normal priors. The proof of the theorem is given in the Appendix (\cite{paul2022posterior}).
\begin{theorem}
    \label{Thm1}
    Suppose $\mathbf{X} \sim \mathcal{N}_n(\boldsymbol{\theta}_0,\mathbf{I}_n)$ where $\boldsymbol{\theta}_0 \in l_0[q_n]$ with $q_n \propto$  $n^{\beta}$ where $0<\beta<1$. Consider the class of priors
	\eqref{eq:1.3} satisfying \eqref{eq:1.4} where $L(t)=(1+1/t)^{-(\frac{1}{2}+\alpha)}$ with $\alpha \geq \frac{1}{2}$.
    Then, for any  empirical Bayes estimator of $\tau$, say $\widehat{\tau}_n$, defined in \eqref{eq:2.2} that satisfies (\hyperlink{C1}{C1}), the empirical Bayes
	posterior contracts around the true parameter $\boldsymbol{\theta}_0$ at the near minimax rate, namely,  
\begin{equation*}
		\sup_{\boldsymbol{\theta}_0 \in l_0[q_n]} \mathbb{E}_{\boldsymbol{\theta}_0} \Pi_{\widehat{\tau}_n} \bigg(  \boldsymbol{\theta}:||\boldsymbol{\theta}-\boldsymbol{\theta}_0 ||^2>    M_n q_n \log n |\mathbf{X} \bigg) \to 0 \hspace{0.1cm},
\end{equation*}
	for any $M_n \to \infty$.
\end{theorem}
The proof of the theorem is given in the Appendix.
\subsection{A Full Bayes Approach}
\label{sec-3.3}
In this section, we consider a full Bayes approach where we assume a non-degenerate prior $\pi_{2n}(\tau)$ on the global shrinkage parameter $\tau$. We are making explicit the dependence on $n$ of the prior on $\tau$.
%We want to investigate the posterior contraction rate of the resulting posterior in this approach. 

Motivated by \cite{van2017adaptive}, we consider priors $\pi_{2n}(\tau)$ which satisfy the conditions (\hyperlink{C2}{C2})
and (\hyperlink{C3}{C3}) below.\\
\textbf{(\hypertarget{C2}{C2})} $\pi_{2n}(.)$ is supported on $[\frac{1}{n},\frac{1}{\log n}]$, where $\frac{1}{\log n}$ is defined in \eqref{eq:2.2}.     \\
%\textbf{(\hypertarget{C2}{C2})} $\pi_{2n}(.)$ is supported on $[\frac{1}{n},1]$.     \\
%\textbf{(\hypertarget{C3}{C3})} Let $t_n= \frac{2a\sqrt{\pi}}{K} \tau_n(q_n)$ where $K=\frac{\Gamma(a+\alpha)}{\Gamma(a)\Gamma(\alpha)}$ as given in the previous subsection. Then $\pi_{2n}(.)$ satisfies
%	\begin{equation*}
%		\bigg(\frac{q_n}{n}\bigg)^{M_1}  \int_{\frac{t_n}{2}}^{t_n} \pi_{2n}(\tau) d \tau \gtrsim e^{-c q_n}, \hspace{0.1cm} \text{for some } c\leq a, M_1 \geq 1 \hspace{0.1cm}.
%	\end{equation*}
\textbf{(\hypertarget{C3}{C3})} Let $t_n= \frac{\sqrt{\pi}}{K} \tau_n(q_n)$ where $K=\frac{\Gamma(\frac{1}{2}+\alpha)}{\sqrt{\pi}\Gamma(\alpha)}$ as given in the previous subsection. Then $\pi_{2n}(.)$ satisfies
	\begin{equation*}
		\bigg(\frac{q_n}{n}\bigg)^{M_1}  \int_{\frac{t_n}{2}}^{t_n} \pi_{2n}(\tau) d \tau \gtrsim e^{-c q_n}, \hspace{0.1cm} \text{for some } c\leq \frac{1}{2}, M_1 \geq 1 \hspace{0.1cm}.
	\end{equation*}

%Condition (\hyperlink{C2}{C2}) is exactly the same as that of \cite{van2017adaptive}, whereas (\hyperlink{C3}{C3}) is a slightly modified version of the condition therein but is asymptotically equivalent to it. \\
{ 
\hspace*{0.5cm} It may be noted that the empirical Bayes estimate of $\tau$ considered in \eqref{eq:2.2} always lie within $[\frac{1}{n},\frac{1}{\log n}]$ and we have proved in Theorem \ref{Thm1} of subsection \ref{sec-3.2} that the empirical Bayes posterior contracts around the truth at a near minimax rate. This is a motivation behind (\hyperlink{C2}{C2}) i.e. the restriction of the support of $\pi_{2n}(\cdot)$ to the interval $[\frac{1}{n},\frac{1}{\log n}]$. (\hyperlink{C3}{C3}) is a slightly modified version of the condition therein but is asymptotically equivalent to it.\\
%See also section 3.2 of \cite{van2017adaptive} in this context. \\ %Additionally, this truncation also solves the computational issues as mentioned by \cite{carvalho2009handling},  \cite{scott2010bayes},  \cite{bogdan2008comparison}, \cite{datta2013asymptotic} and \cite{van2017adaptive} when the estimate of $\tau$ is very small using an empirical Bayes approach. Since, the support of $\pi_{2n}(\cdot)$ is exactly the same as that of the empirical Bayes estimate defined in \eqref{eq:2.4}, one can also expect the same problem to be solved for the full Bayes approach, too.} \\
\hspace*{0.5cm} When $q_n$ is known, \cite{van2014horseshoe} proved that the horseshoe posterior attains the minimax rate when the global shrinkage parameter $\tau$ is chosen of the order of $\tau_n(q_n) = (q_n/n)\sqrt{\log (n/(q_n))}$. {\cite{van2014horseshoe} also stated that a choice of $\tau$ that goes to zero slower than $\tau_n(q_n)$, say the $(\frac{q_n}{n})^{\alpha}, 0<\alpha<1$, will lead to a sub-optimal rate of posterior contraction in the squared $L_2$ sense}. \cite{ghosh2017asymptotic} had the same observations for a broad class of priors.
Condition (\hyperlink{C3}{C3}) guarantees that there will be sufficient prior mass around the \textit{optimal} value of $\tau$  so that the full Bayes posterior distribution may be expected to have a high probability around the optimal values and may lead to a near-optimal rate of contraction. This condition is satisfied by many prior densities, except in the very sparse cases when $q_n \lesssim \log n$. Using the same arguments as in \cite{van2017adaptive}, it can be shown that the half Cauchy distribution supported on $[\frac{1}{n},\frac{1}{\log n}]$ having density $\pi_{2n}(\tau)=(\arctan(\frac{1}{\log n})-\arctan(\frac{1}{n}))^{-1}(1+\tau^2)^{-1}\mathbf{1}_{\tau \in [\frac{1}{n},\frac{1}{\log n}]}$ and uniform prior on  $[\frac{1}{n},\frac{1}{\log n}]$ with density $\pi_{2n}(\tau)=1/(\frac{1}{\log n}-\frac{1}{n})\mathbf{1}_{\tau \in [\frac{1}{n},\frac{1}{\log n}]}$ satisfy condition (\hyperlink{C3}{C3}) provided $q_n \gtrsim \log n$. \\ %Similarly, satisfies condition (\hyperlink{C2}{C2}) under the same assumption.\\
\hspace*{0.5cm} Note that, using (\hyperlink{C2}{C2}), the posterior mean of $\theta_i$ in the full Bayes approach can be written as
\begin{align*} \label{eq:2.5}
	\widehat{\theta}_i = E(\theta_i|\mathbf{X}) &= \int_{\frac{1}{n}}^{\frac{1}{\log n}} E(\theta_i|\mathbf{X},\tau) \pi_{2n}(\tau|\mathbf{X}) d\tau \\
	&= \int_{\frac{1}{n}}^{\frac{1}{\log n}} E(1-\kappa_i|X_i,\tau) X_i \cdot \pi_{2n}(\tau|\mathbf{X}) d\tau \\
	&= \int_{\frac{1}{n}}^{\frac{1}{\log n}}   T_{\tau}(X_i) \pi_{2n}(\tau|\mathbf{X}) d\tau. \tag{12}
\end{align*}
 Thus the full Bayes posterior mean is expressed as a weighted average of the posterior mean given $\tau$, the weight being the posterior density of $\tau$ given data.\\ 
 \hspace*{0.5cm} The next lemma is crucial for deriving the optimal posterior contraction rate for the three parameter beta normal priors. This states that the full Bayes posterior of $\tau$ asymptotically concentrates its mass in a neighbourhood of zero of size at most a constant multiple of $t_n$ away from zero. 
 %The proof of the lemma is provided in the Appendix (\cite{paul2022posterior}).
  \begin{lemma}
     \label{lem2}
     If (\hyperlink{C2}{C2})
and (\hyperlink{C3}{C3}) are satisfied, then
\begin{align*}
    \sup_{\theta_0 \in l_0[q_n]} E_{\theta_0} \Pi(\tau \geq 5t_n |\mathbf{X}) \to 0.
\end{align*}
 \end{lemma}
 \begin{proof}[Proof of Lemma \ref{lem2}]
     Using a similar set of arguments used in Lemma 3.6 of \cite{van2017adaptive}, with $P_{\theta_0}$-probability tending to one, 
     \begin{align*} \label{eq:2.9}
       \Pi(\tau \geq 5t_n |\mathbf{X}) & \lesssim \frac{e^{-0.5q_n}}{\int_{\frac{t_n}{2}}^{t_n} \pi_{2n}(\tau) d \tau}, \tag{13}
     \end{align*}
     which tends to zero using (\hyperlink{C3}{C3}).
     As a result, for $b_n=(\frac{q_n}{n})^{M_1}$ and $h(\mathbf{X})= \Pi(\tau \geq 5t_n |\mathbf{X})$,
     \begin{align*}
         E_{\theta_0} h(\mathbf{X}) &= E_{\theta_0}[h(\mathbf{X}) 1_{\{h(\mathbf{X}) \lesssim b_n\} }] +E_{\theta_0}[h(\mathbf{X}) 1_{\{h(\mathbf{X}) \gtrsim b_n\} }] \\
         & \lesssim b_n +P_{\theta_0}[h(\mathbf{X}) \gtrsim b_n]
     \end{align*}
     The proof follows from this since by definition the first term goes to zero and on the other hand, the second term goes to zero due to \eqref{eq:2.9}.
 \end{proof}
 %It is clear from the proofs of Theorems \ref{Thm-3} and \ref{Thm-4} that the above representation indeed greatly simplifies the arguments for establishing the near minimax concentration of the mean parameter vector of our interest.\\
 \begin{remark}
     For obtaining arguments similar to  Lemma 3.6 of \cite{van2017adaptive}, we need to come up with an upper bound on the derivative of the log-likelihood of $\mathbf{X}$ given $\tau$. This is obtained with the help of Lemmas 1-6 given in the Appendix. Based on the range of $\tau$, one needs to find another upper bound on this upper bound and use this to obtain \eqref{eq:2.9} with the help of Bayes' theorem.
 \end{remark}
\hspace*{0.5cm} The next theorem proves that the hierarchical Bayes posterior estimate $\widehat{\boldsymbol{\theta}}$ of $\boldsymbol{\theta}$ also attains near minimax $L_2$ risk. The proof of the theorem is given in the Appendix.
\begin{theorem}
    \label{Thm-3}
	Suppose $\mathbf{X} \sim \mathcal{N}_n(\boldsymbol{\theta}_0,\mathbf{I}_n)$ where $\boldsymbol{\theta}_0 \in l_0[q_n]$ with $q_n \propto$  $n^{\beta}$ where $0<\beta<1$. Consider the class of priors
	\eqref{eq:1.3} satisfying \eqref{eq:1.4} where $L(t)=(1+1/t)^{-(\frac{1}{2}+\alpha)}$ with $\alpha \geq \frac{1}{2}$. If the prior on $\tau$ satisfies conditions (\hyperlink{C2}{C2}) and (\hyperlink{C3}{C3}), then the full Bayes
	posterior contracts around the true parameter $\boldsymbol{\theta}_0$ at the near minimax rate, namely,
	\begin{equation*}
		\sup_{\boldsymbol{\theta}_0 \in l_0[q_n]} \mathbb{E}_{\boldsymbol{\theta}_0} \Pi \bigg(  \boldsymbol{\theta}:||\boldsymbol{\theta}-\boldsymbol{\theta}_0 ||^2>    M_n q_n \log n |\mathbf{X} \bigg) \to 0 ,
	\end{equation*}
	for any $M_n \to \infty$.
\end{theorem}
%\cite{van2017adaptive} established the same rate of convergence for horseshoe prior. Theorem \ref{Thm-3} generalizes their result for the three parameter beta normal priors.
\section{Asymptotic optimality of a full Bayes approach in a problem of hypothesis testing}
\label{sec-4}
Suppose we have observations $X_i,i=1,2,\cdots,n$, modeled through \eqref{eq:1.1} and our problem is to test simultaneously the hypotheses $H_{0i}: \theta_i =0$ against $H_{1i}: \theta_i \neq 0$ for $i=1,2,\cdots, n$. Our interest is in the sparse asymptotic regime when $n$ is large and the proportion of non-zero means is expected to be small. Formally, we will be using the asymptotic framework of \cite{bogdan2011asymptotic} in this section. As in \cite{bogdan2011asymptotic}, we introduce a set of iid binary latent variables $\nu_1,\nu_2, \cdots, \nu_n$ where $\nu_i=0$ refers to the event that $H_{0i}$ is true and $\nu_i=1$ denotes that $H_{1i}$ is true with $P(\nu_i =1)=p$ for $i=1,2,\cdots,n$. Given $\nu_i=0$, $\theta_i$ has a degenerate distribution at $0$ and when $\nu_i=1$, $\theta_i$ is assumed to have a $\mathcal{N}(0,\psi^2)$ distribution with $\psi^2>0$. In this Bayesian setting, the priors on $\theta_i$'s are given by 
\begin{align*} \label{eq:4.1}
	\theta_i & \simiid{} (1-p) \delta_{0}+p \mathcal{N}(0,\psi^2), i=1,2,\cdots,n. \tag{14}
\end{align*}
The marginal distribution of $X_i$ for two-groups model is thus of the form
\begin{align*} \label{eq:4.2}
	X_i & \simiid{} (1-p) \mathcal{N}(0,1)+ p \mathcal{N}(0,1+\psi^2), i=1,2,\cdots,n. \tag{15}
\end{align*}
Our testing problem is equivalent to test simultaneously
\begin{equation*} \label{eq:4.3}
	H_{0i}: \nu_i =0 \hspace*{0.2cm}  \text{versus} \hspace*{0.2cm} H_{1i}: \nu_i =1, \hspace*{0.2cm} \text{for}  \hspace*{0.2cm} i=1,2,\cdots,n. \tag{16}
\end{equation*}
\hspace*{0.5cm}   The overall loss in a multiple testing rule is defined as the total number of misclassifications made by the test which is nothing but a sum $0-1$ losses corresponding to each individual test. If $t_{1i}$ and $t_{2i}$ denote the probabilities of type I and type II errors respectively for the $i^{\text{th}}$ test, the Bayes risk $R$ of a multiple testing procedure is given by
\begin{equation*} \label{eq:4.4}
	R= \sum_{i=1}^{n} [(1-p) t_{1i}+ p t_{2i}]  . \tag{17}
\end{equation*}
 \cite{bogdan2011asymptotic} showed that the multiple testing rule  minimizing the Bayes risk \eqref{eq:4.4} rejects $H_{0i}$ if 
\begin{equation*} \label{eq:4.5}
	\pi(\nu_i=1|X_i) >\frac{1}{2} ,\hspace*{0.2cm} \text{i.e.}  \hspace*{0.2cm} X^2_i >c^2 \tag{18}
\end{equation*}
and accepts $H_{0i}$ otherwise, for each $i=1,2,\cdots n$.
Here $\pi(\nu_i=1|X_i)$ denotes the posterior probability of $H_{1i}$ to be true and $c^2 \equiv c^2_{p,\psi}= \frac{1+\psi^2}{\psi^2} (\log (1+\psi^2)+ 2 \log f)$ with $f=\frac{1-p}{p}$. Due to the presence of unknown model parameters $p$ and $\psi^2$ in \eqref{eq:4.5}, this rule is termed the Bayes Oracle. For defining their asymptotic framework, \cite{bogdan2011asymptotic} introduced two parameters $u \equiv u_n =\psi^2_n$ and $v \equiv v_n=u_n f^2_n$ and also assumed the following.  \\
%\vskip 6pt

\textbf{\hypertarget{assumption2}{Assumption 2.}} $p_n \to 0$, $u_n =\psi^2_n \to \infty$ and $\frac{\log v_n}{u_n} \to C \in (0,\infty)$ as $n \to \infty$. 
\vskip 3pt
For simplicity of notation, we will henceforth remove the dependence of the model parameters on $n$. Under \hyperlink{assumption2}{Assumption 2}, the asymptotic expression for the Bayes risk of the Bayes Oracle \eqref{eq:4.5}, denoted by $R^{\text{BO}}_{\text{Opt}}$ is of the form
\begin{equation*}  \label{eq:4.6}
	R^{\text{BO}}_{\text{Opt}}= np{\bigg[2 \Phi(\sqrt{{C}})-1 \bigg]} (1+o(1)).  \tag{19}
\end{equation*}
\hspace*{0.5cm} A multiple testing rule is defined to Asymptotically Bayes Optimal under Sparsity (ABOS) if the ratio of its Bayes risk and $R^{\text{BO}}_{\text{Opt}}$ converges to $1$ under the asymptotic framework of 
\hyperlink{assumption2}{Assumption 2}.\\
%\hspace*{0.05cm} As mentioned earlier, \cite{carvalho2010horseshoe} proposed using a testing rule based on the horseshoe prior in the sparse normal means problem. They observed that $E(1-\kappa_i|\text{data})$ has roughly the same role as the posterior probability $\pi(\nu_i=1|\text{data})$ and based on that introduced an intuitive testing rule based on the former which performed very well in their simulations (See Section 3.4 of \cite{carvalho2010horseshoe}). Their proposed decision rule is as follows:  
%\begin{equation*}
	%\text{reject} \hspace*{0.2cm} H_{0i} \hspace*{0.2cm} \text{if} \hspace*{0.2cm} E(1-\kappa_i|\text{data}) > \frac{1}{2}, i=1,2,\cdots,n.
%\end{equation*}
\hspace*{0.5cm}As mentioned in the introduction, \cite{datta2013asymptotic} first studied from a decision theoretic viewpoint, the asymptotic optimality of a testing rule, based on the horseshoe prior in this problem. Treating the global parameter $\tau$ as a tuning parameter, they studied the following decision rule for horseshoe prior,
\begin{equation*} \label{eq:4.7}
	\text{reject} \hspace*{0.2cm} H_{0i} \hspace*{0.2cm} \text{if} \hspace*{0.2cm} E(1-{\kappa}_i|X_i, \tau) > \frac{1}{2}, i=1,2,\cdots,n. \tag{20}
\end{equation*}
They showed that under \hyperlink{assumption2}{Assumption 2}, their proposed decision rule \eqref{eq:4.7} attains the optimal Bayes risk \eqref{eq:4.6} up to some multiplicative constant if $\tau=O(p)$. Later \cite{ghosh2016asymptotic} obtained
similar results (under the condition that $\lim_{n \to \infty} \frac{\tau}{p} \in (0,\infty)$) using the same testing rule as above but modeling $\theta_i$ by a general class of one-group global-local shrinkage priors where $\tau$ is treated as a tuning parameter. By choosing $\tau \sim p$ and considering a prior of the form \eqref{eq:1.3} satisfying \eqref{eq:1.4} with $a=\frac{1}{2}$ with $L(\cdot)$ satisfying  \hyperlink{assumption1}{Assumption 1}, \cite{ghosh2017asymptotic} proved that the induced decision rule is ABOS for the subclass of priors, named as
horseshoe-type priors by \cite{ghosh2017asymptotic}.  This result demonstrates that 
appropriately chosen one-group prior can be a reasonable alternative to the two-groups prior in this problem.\\
\hspace*{0.5cm} When $p$ is unknown, \cite{ghosh2017asymptotic} considered an empirical Bayes version of \eqref{eq:4.7} for the same class of priors as above and used the following decision rule, 
\begin{equation*} \label{eq:4.8}
	\text{reject} \hspace*{0.2cm} H_{0i} \hspace*{0.2cm} \text{if} \hspace*{0.2cm} E(1-{\kappa}_i|X_i, \widehat{\tau} ) > \frac{1}{2}, i=1,2,\cdots,n, \tag{21}
\end{equation*}
where $\widehat{\tau}$ is defined in \eqref{eq:2.4}. Under \hyperlink{assumption2}{Assumption 2}, they proved that for $p_n \propto n^{-\epsilon} ,\; 0< \epsilon<1$ and for $a=\frac{1}{2}$ with $L(\cdot)$ defined in \eqref{eq:1.4} satisfying  \hyperlink{assumption1}{Assumption 1}, the decision rule \eqref{eq:4.8} is also ABOS.\\
\hspace*{0.5cm}We are interested in  a full Bayes approach to the same problem where one puts a non-degenerate prior on $\tau$. In that case the decision rule modifies to 
\begin{equation*} \label{eq:4.9}
	\text{reject} \hspace*{0.2cm} H_{0i} \hspace*{0.2cm} \text{if} \hspace*{0.2cm} E(1-{\kappa}_i|\mathbf{X} ) > \frac{1}{2}, i=1,2,\cdots,n. \tag{22}
\end{equation*} 
Note that, the rule \eqref{eq:4.9} depends on the entire data set and is clearly different from \eqref{eq:4.7} and \eqref{eq:4.8}. The study of asymptotic optimality of a rule like \eqref{eq:4.9} is still missing in the literature, even for the horseshoe prior. 
 \cite{ghosh2016asymptotic} had done some simulation studies using the decision rule \eqref{eq:4.9} when a standard half Cauchy prior was used for $\tau$. Their results obtained through simulations motivate us to study the decision rule \eqref{eq:4.9} in a formal decision-theoretic manner. Our goal is to find a general class of priors on $\tau$ under which the decision rule \eqref{eq:4.9} behaves in an optimal or near-optimal manner in terms of Bayes risk.\\
\hspace*{0.5cm} We start with a general class of priors $\pi_{2n}(\cdot)$  satisfying the following condition: \\
\vskip 0.1pt
\textbf{({\hypertarget{C4}{C4}})}
$\int_{\frac{1}{n}}^{\alpha_n}\pi_{2n}(\tau) d \tau=1$ for some $\alpha_n$ such that $\alpha_n \to 0$ and $n \alpha_n \to \infty$ as $n \to  \infty$.\\
\vskip 0.1pt
%The motivation behind this initial choice comes from the  results of \cite{datta2013asymptotic}, \cite{ghosh2016asymptotic}, and \cite{ghosh2017asymptotic} where $\tau$ is used as a tuning parameter to be fixed using the knowledge of $p$. 
It may be recalled that for proving the optimality of the decision rule based on \eqref{eq:4.7}, a choice of $\tau=O(p)$ is sufficient. In the present situation, $p$ is unknown, so intuition suggests that it might be a good idea to start with a prior on $\tau$ which gives sufficient prior probability around such values of $\tau$. We are basically interested in situations where the unknown $p \propto n^{-\epsilon}, 0<\epsilon \leq 1$. So the choice of the lower boundary of the support as $\frac{1}{n}$ at least ensures positive prior mass around $\frac{1}{n}$ and $\frac{1}{n}=O(p)$ as $n \to \infty$ when $p \propto n^{-\epsilon}, 0<\epsilon \leq 1$. Given this, the specific choices of $\alpha_n$ that ensure asymptotic optimality become clear from Theorems \ref{Thm-6} and \ref{Thm-7} and the discussions thereafter. In order to keep the support of $\tau$ as broad as possible, we assume $n \alpha_n \to \infty$ as $n \to  \infty$.\\
%\hspace*{0.05cm} First, we state the reason behind a specific choice of the lower bound on the range of the prior density of $\tau$ in ({\hyperlink{C3}{C3}}).
%Since $p$ denotes the theoretical proportion of non-nulls in the population, its least possible value is $\frac{1}{n}$, hence it acts as the lower bound of $\pi_n(\tau)$ in ({\hyperlink{C3}{C3}}).
%A specific choice of $\alpha_n$ will be provided later. Since $p$ is unknown, our target is to choose $\alpha_n$ in such a way that depends on the number of hypotheses $n$ only.
%The next theorem provides the upper bound on the probability of type I error using the decision rule \eqref{eq:4.9}. This also shows the impact on the choice of $\alpha_n$ in this testing problem. Proof of the theorem can be found in section \ref{sec-6}. This condition ({\hyperlink{C3}{C3}}) indicates that instead of giving importance to the choice of prior, one can use any non-degenerate prior in a small interval of $\tau$ to derive some meaningful results related to the optimal Bayes risk in one-group setup.

\begin{theorem}
	\label{Thm-6}
	Let $X_1,X_2,\cdots,X_n$ be iid with distribution  \eqref{eq:4.2}. Suppose we test \eqref{eq:4.3} using decision rule \eqref{eq:4.9} induced by the class of priors \eqref{eq:1.3} with $\pi_{2n}(\cdot)$ satisfying ({\hyperlink{C4}{C4}}) and $\pi_1(\cdot)$ 
 satisfying \eqref{eq:1.4} with $a \in (0,1)$, where $L(\cdot)$ satisfies  \hyperlink{assumption1}{Assumption 1}. Then the probability of type-I error of the decision rule \eqref{eq:4.9}, denoted  $t^{\text{FB}}_{1i}$, satisfies
	\begin{equation*}
		t^{\text{FB}}_{1i}  \lesssim \frac{{\alpha}^{2a}_n}{\sqrt{\log (\frac{1}{\alpha^2_n}) }} (1+o(1)),
	\end{equation*}
	where $o(1)$ term is independent of $i$ and depends only on $n$ such that $\lim_{n \to \infty}o(1)=0$.
 \end{theorem}
\hspace*{0.5cm}Let us now state a result on the upper bound on the probability of type II error induced by the decision rule \eqref{eq:4.9}. %Proof of this theorem can be found at the end of this section.
\begin{theorem}
	\label{Thm-7}
	Consider the set-up of Theorem \ref{Thm-6}. Under \hyperlink{assumption2}{Assumption 2}, for any fixed $\eta \in (0,1)$, $\delta \in (0,1)$ and any fixed $\rho > \frac{2}{\eta (1-\delta)}$, the probability of type-II error induced by the decision rule \eqref{eq:4.9}, denoted  $t^{\text{FB}}_{2i}$ satisfies
	\begin{equation*}
	 		t^{\text{FB}}_{2i} \leq
		\bigg[2 \Phi \bigg(\sqrt{ \frac{a \rho C}{\epsilon}}\bigg)-1 \bigg](1+o(1)).
	\end{equation*}
	where $o(1)$ term is independent of $i$ and depends on $n, \eta, \delta$ and $\rho$ such that $\lim_{n \to \infty}o(1)=0$.
\end{theorem}

\begin{remark}
    Combining Theorems \ref{Thm-6} and \ref{Thm-7}, we immediately see that the Bayes risk of the decision rule \eqref{eq:4.9}, denoted $R^{\text{FB}}_{\text{OG}}$ satisfies
    \begin{equation*}
        R^{\text{FB}}_{\text{OG}} \leq n \bigg[(1-p)\frac{{\alpha}^{2a}_n}{\sqrt{\log (\frac{1}{\alpha^2_n}) }} (1+o(1)) +p \bigg[2 \Phi \bigg(\sqrt{ \frac{a \rho C}{\epsilon}}\bigg)-1 \bigg](1+o(1))\bigg ],
    \end{equation*}
   where $o(1)$ term is independent of $i$ and depends on $n, \eta, \delta$ and $\rho$ such that $\lim_{n \to \infty}o(1)=0$. We are now in a position to choose $\alpha_n, p$ and $a$ appropriately so that this upper bound can be made of the same order as $R^{\text{BO}}_{\text{Opt}}$ in \eqref{eq:4.6}. The result is stated in Theorem \ref{Thm-8} below. 
\end{remark}

\begin{theorem}
	\label{Thm-8}
	Let $X_1,X_2,\cdots,X_n$ be iid with distribution  \eqref{eq:4.2}. Suppose we test \eqref{eq:4.3} using decision rule \eqref{eq:4.9} induced by the class of priors \eqref{eq:1.3} with $\pi_{2n}(\cdot)$ satisfying ({\hyperlink{C4}{C4}}) and $\pi_1(\cdot)$ 
 satisfying \eqref{eq:1.4} with $a \in [0.5,1)$, where $L(\cdot)$ satisfies  \hyperlink{assumption1}{Assumption 1}. 
Also assume $\log (\frac{1}{\alpha_n})=\log n-\frac{1}{2} \log \log n+c_n$, where $c_n \to \infty$ such that $c_n=o(\log \log n)$ as $n \to \infty$. Then under \hyperlink{assumption2}{Assumption 2}, for any fixed $\eta \in (0,1)$, $\delta \in (0,1)$ and any fixed $\rho > \frac{2}{\eta (1-\delta)}$, the Bayes risk of the decision rule \eqref{eq:4.9}, denoted $R^{\text{FB}}_{\text{OG}}$ satisfies
	\begin{equation*}
			 R^{\text{FB}}_{\text{OG}} \leq
		np \bigg[2 \Phi \bigg(\sqrt{ \frac{a \rho C}{\epsilon}} \bigg)-1 \bigg](1+o(1)),
	\end{equation*}
	where $o(1)$ term is independent of $i$ and depends on $n, \eta, \delta, \rho$ such that $\lim_{n \to \infty}o(1)=0$. \\
	For the extremely sparse case i.e. for $\epsilon=1$ with $a=0.5$,
	\begin{equation*}
		\lim_{n \to \infty}	\frac{ R^{\text{FB}}_{\text{OG}}}{ R^{\text{BO}}_{\text{Opt}}} = 1.
	\end{equation*}
	
\end{theorem}

\section{Concluding remarks}
\label{sec-5}
In this paper, we study inference on a high-dimensional sparse normal means model when the proportion of non-zero means is unknown. 
Assuming the proportion of non-zero means
to be known, several results have been established both in terms of optimal contraction rate (\cite{van2014horseshoe}, \cite{van2016conditions}, \cite{ghosh2017asymptotic}) and optimal Bayes risk (\cite{datta2013asymptotic}, \cite{ghosh2016asymptotic}, \cite{ghosh2017asymptotic}) in the literature on one-group global-local priors.
In this work, we have established that for both the empirical Bayes and full Bayes approaches, the posterior distributions of the mean vector contract around the truth at a near minimax rate for broad classes of priors (using the simple estimator of \cite{van2014horseshoe} and the MMLE of \cite{van2017adaptive}). We need to come up with several new arguments along with sharper upper bounds on the posterior shrinkage coefficient and on the posterior variance of the mean vector in order to generalize the results of \cite{van2014horseshoe} and \cite{van2017adaptive} for the class of priors mentioned before.  
 We then study the asymptotic properties of a decision rule in the full Bayes approach based on global-local priors in a problem of simultaneous hypothesis testing. We employ a full Bayes approach using a broad class of prior densities on $\tau$. Within the asymptotic framework of \cite{bogdan2011asymptotic}, we have been able to prove that the Bayes risk induced by this broad class of priors attains the optimal Bayes risk, up to some multiplicative constant. Moreover, when the underlying model is extremely sparse, the decision rule is indeed ABOS for horseshoe-type priors. Our result also indicates that the upper bound of the Bayes risk can be substantially sharper than those of \cite{ghosh2016asymptotic} and \cite{datta2013asymptotic} for a wide range of sparsity levels and the whole class of priors. We hope that our results will reinforce the logic for using the global-local priors in sparse situations in place of their two-group counterparts.\\
\hspace*{0.5cm} Our chosen classes of priors are rich enough to include many interesting priors, but due to the assumption of boundedness on the slowly varying function,  Horseshoe+ prior is not included in this class. Thus our techniques can't be applied for the Horseshoe+ due to \cite{bhadra2017horseshoe+}. \cite{van2016conditions} proved that the optimal posterior concentration result also holds for Horseshoe+ using appropriate choice of $\tau$. No such result is still available in the literature for the same when either the empirical Bayes estimate of $\tau$ based on maximizing the marginal likelihood function is used or a full Bayes procedure is employed. We would like to consider this problem elsewhere.\\
\hspace*{0.5cm} {It is worth mentioning that though most of our results related to optimal contraction rate and asymptotic Bayes optimality are valid for a broad class of global-local priors, but the results related to the posterior concentration rate either based on MMLE or the full Bayes approach, have been proved only for the three parameter beta normal class of priors. The techniques used in this context does not be seem to be applicable to another class of priors, namely the Generalized double Pareto priors.
Thus a question which is still unanswered is whether the same results hold for the Generalized double Pareto priors. This might be considered as an interesting problem to be studied in the future.} \\
\hspace*{0.5cm} Throughout this work, we use the concept of slowly varying functions. \cite{van2016conditions} considered a broader class of global-local priors that contains our chosen class of priors, by using the concept of ``uniformly regularly varying functions". They established optimal posterior contraction rate for this class of priors when the level of sparsity is known, by choosing the global shrinkage parameter $\tau$ appropriately based on the level of sparsity. To the best of our knowledge, no such result is available in the literature when the level of sparsity is unknown. These priors have not been studied in the context of multiple testing either when the proportion of non-null means is unknown. These will be interesting problems to consider and we expect that some of the techniques employed in this work will be helpful in these studies also.

%% *** Frontmatter *** 

%\title{\support{}}
%\runtitle{Posterior Contraction Rate and Asymptotic Bayes Optimality-supp}

%\thankstext{<id>}{<text>}

%% ** Keywords **
%\begin{keyword}%[class=MSC]
%\kwd{asymptotic minimaxity, near-minimax rate, one-group prior, sparsity, ABOS}
%\kwd[]{}
%\end{keyword}

%% ** Mainmatter **

%\section{}\label{}

% \begin{figure} 
% \includegraphics{<eps-file>}% place <eps-file> in ./img  subfolder
% \caption{}
% \label{}
% \end{figure}

% \begin{table} 
% *****************
% \begin{tabular}{lll}
% \end{tabular}
% *****************
% \caption{}
% \label{}
% \end{figure}

%%%%%%%%%%%%%%%%%%%%%%%%%%%%%%%%%%%%%%%%%%%%%%
%% Supplementary Material, if any, should   %%
%% be provided in {supplement} environment  %%
%% with title and short description.        %%
%%%%%%%%%%%%%%%%%%%%%%%%%%%%%%%%%%%%%%%%%%%%%%
%\begin{supplement}
%\stitle{???}
%\sdescription{???.}
%\end{supplement}

%% ** The bibliograhy **
%\bibliographystyle{ba}
%\bibliography{<bib-data-file>}% place <bib-data-file> in ./bib folder 

% ** Acknowledgements **
%\begin{acks}[Acknowledgments]
%\end{acks}

%%%%%%%%%%%%%%%%%%%%%%%%%%%%%%%%%%%%%%%%%%%%%%
%% Supplementary Material, if any, should   %%
%% be provided in {supplement} environment  %%
%% with title and short description.        %%
%%%%%%%%%%%%%%%%%%%%%%%%%%%%%%%%%%%%%%%%%%%%%%
%\begin{supplement}
%\label{supp}
%\stitle{Proof of Theorem \ref{Thm-4}}
%\sdescription{Short description of Supplement A.}

\appendix
\section{Appendix}
\label{A}

\subsection{Results related to the contraction rate of Empirical Bayes Procedure--plug-in Estimator}

\begin{proof}[Proof of Theorem \ref{Thm-1}]
For proving Theorem 1, first we need to show that, the total posterior variance using the empirical Bayes estimator of $\tau$	corresponding to this class of priors satisfies, as $n \to \infty$,
	\begin{equation*} \label{eq:2.2.1}
	\sup_{\boldsymbol{\theta}_0 \in l_0[q_n]} \mathbb{E}_{\boldsymbol{\theta}_0} \sum_{i=1}^{n} Var(\theta_i|X_i,\widehat{\tau}) \lesssim q_n \log n, \tag{2}
	\end{equation*}
where $Var(\theta_i|X_i,\widehat{\tau})$ denotes $Var(\theta_i|X_i,{\tau})$ evaluated at $\tau=\widehat{\tau}$.
We make use of Markov's inequality along with \eqref{eq:2.2.1} and Theorem 2 of \cite{ghosh2017asymptotic}.\\
Now we move towards establishing \eqref{eq:2.2.1}.	Let us define $\Tilde{q}_n=\sum_{i=1}^{n}1_{\{\theta_{0i}\neq 0\}}$. Thus, $\Tilde{q}_n \leq q_n$. We then note
	%First we split $E_{\boldsymbol{\theta}_0} \sum_{i=1}^{n} Var(\theta_i|X_i,\widehat{\tau})$ as 
	\begin{equation}\label{eq:T-A1.1}
		\mathbb{E}_{\boldsymbol{\theta}_0} \sum_{i=1}^{n} Var(\theta_i|X_i,\widehat{\tau})= \sum_{i:\theta_{0i} \neq 0} \mathbb{E}_{{\theta}_{0i}}  Var(\theta_i|X_i,\widehat{\tau})+\sum_{i:\theta_{0i} = 0} \mathbb{E}_{{\theta}_{0i}}  Var(\theta_i|X_i,\widehat{\tau}) \hspace{0.1cm}. \tag{3}
	\end{equation}
	We now prove that $\sum_{i:\theta_{0i} \neq 0} \mathbb{E}_{{\theta}_{0i}}  Var(\theta_i|X_i,\widehat{\tau}) \lesssim \tilde{q}_n \log n$ and   
		$\sum_{i:\theta_{0i} = 0} \mathbb{E}_{{\theta}_{0i}}  Var(\theta_i|X_i,\widehat{\tau}) \lesssim  q_n \log n\hspace{0.1cm}$ in \textbf{Step-1} and \textbf{Step-2} below respectively. Combining these results we get the final result. Now let us prove \textbf{Step-1} and \textbf{Step-2}. \\\\
	\textbf{Proof of Step-1:} Fix any $c>1$ and choose $\rho>c$.
 Define $r_n= \sqrt{4a \rho^2 \log n}$. Fix any $i$ such that $\theta_{0i} \neq 0$. We split $\mathbb{E}_{{\theta}_{0i}}  Var(\theta_i|X_i,\widehat{\tau})$ as \begin{equation*}\label{eq:T-A1.2}
		\mathbb{E}_{{\theta}_{0i}}  Var(\theta_i|X_i,\widehat{\tau})=\mathbb{E}_{{\theta}_{0i}}[Var(\theta_i|X_i,\widehat{\tau}) 1_{\{|X_i|\leq r_n \}}] +\mathbb{E}_{{\theta}_{0i}}[Var(\theta_i|X_i,\widehat{\tau}) 1_{\{|X_i|> r_n\}}] \hspace{0.1cm}. \tag{4}
	\end{equation*}
	Since for any fixed $x_i \in \mathbb{R}$ and $\tau>0$,
	$Var(\theta_i|x_i,\tau) \leq 1+x^2_i$.
	\begin{equation*}\label{eq:T-A1.3}
		\mathbb{E}_{{\theta}_{0i}}[Var(\theta_i|X_i,\widehat{\tau}) 1_{\{|X_i|\leq \sqrt{4a \rho^2 \log n}\}}] \leq 1+ 4a \rho^2 \log n \hspace{0.1cm}. \tag{5}
	\end{equation*}
	Note that, for any fixed $x_i \in \mathbb{R}$, $x^2_iE(\kappa^2_i|x_i,\tau)=x^2_i E(\frac{1}{(1+\lambda^2_i\tau^2)^2}|x_i,\tau)$ is non-increasing in $\tau$. Also using Lemma A.1 in \cite{ghosh2017asymptotic}, $Var(\theta_i|x_i,\tau) \leq 1+x^2_iE(\kappa^2_i|x_i,\tau)$. Using these two results, 
	\begin{align*}
		Var(\theta_i|x_i,\widehat{\tau}) & \leq 1+x^2_iE(\kappa^2_i|x_i,\widehat{\tau}) \leq 1+x^2_iE(\kappa^2_i|x_i,\frac{1}{n}) \hspace{0.1cm} [\text{Since} \hspace{0.1cm} \widehat{\tau}\geq \frac{1}{n} ]  \hspace{0.1cm}.
	\end{align*}
	Using arguments similar to Lemma 3 of \cite{ghosh2017asymptotic}, one can show that for any fixed $\eta \in (0,1)$ and $\delta \in (0,1)$,
 the r.h.s. above can be bounded by a
	non-negative and measurable real-valued function $\tilde{h}(x_i,\tau, \eta, \delta)$ satisfying $\tilde{h}(x_i,\tau, \eta, \delta) = 1+ \tilde{h}_1(x_i,\tau)+\tilde{h}_2(x_i,\tau,\eta, \delta)$ with
	\begin{equation*}
		\tilde{h}_1(x_i,\tau) =C_{**} \bigg[ x^2_i  \int_{0}^{\frac{ x^2_i}{1+t_0}} \exp(-\frac{u}{2})u^{a+\frac{1}{2}-1} du\bigg]^{-1},
	\end{equation*}
	where $C_{**}$ is a global constant independent of both $x_i$ and $\tau$ and
	\begin{equation*}
		\tilde{h}_2(x_i,\tau,\eta, \delta)= x^2_i \frac{H(a,\eta,\delta)}{\Delta(\tau^2, \eta,\delta)}{\tau}^{-2a}e^{-\frac{\eta(1-\delta)x^2_i}{2}} \hspace*{0.05cm},
	\end{equation*}
	where $H(a,\eta,\delta)=\frac{(a+\frac{1}{2})(1-\eta \delta)^a}{K ({\eta \delta)}^{(a+\frac{1}{2})}}$ and 
		$\Delta(\tau^2, \eta,\delta)= \frac{\int_{\frac{1}{\tau^2}(\frac{1}{\eta \delta}-1)}^{\infty}t^{-(a+\frac{3}{2})}L(t) dt}{(a+\frac{1}{2})^{-1} {(\frac{1}{\tau^2}(\frac{1}{\eta \delta}-1))}^{-(a+\frac{1}{2})}}$.\\
	 Since, $\tilde{h}_1(x_i,\tau)$ is strictly decreasing in $|x_i|$ for any fixed $\tau$, one has,
	\begin{align*}
		\sup_{|x_i|>\sqrt{4a \rho^2 \log n}} \tilde{h}_1(x_i,\frac{1}{n}) & \leq C_{**} \bigg[ \rho^2 \log n \int_{0}^{\frac{4a \rho^2 \log n}{1+t_0}} \exp(-\frac{u}{2})u^{a+\frac{1}{2}-1} du\bigg]^{-1} \lesssim \frac{1}{\log n} \hspace*{0.05cm} \cdot
	\end{align*}
	Also noting that $\tilde{h}_2(x_i,\tau,\eta, \delta)$ is strictly decreasing in $|x_i|>C_1=\sqrt{\frac{2}{\eta(1-\delta)}}$, for any $\rho> C_1$,
	\begin{align*}
		\sup_{|x_i|>\sqrt{4a \rho^2 \log n}} \tilde{h}_2(x_i,\frac{1}{n},\eta, \delta) & \leq 4a \rho^2 \frac{H(a,\eta,\delta)}{\Delta(\frac{1}{n^2},\eta,\delta)} \log n \cdot n^{-2a(\frac{2\rho^2}{{C}^2_1}-1)} \hspace*{0.05cm} \cdot
	\end{align*}
	Using the definition of the slowly varying function $L(\cdot)$ and assumption (\hyperlink{A1}{A1}) in the right-hand side of above, we have, 
 \begin{equation*}
     \sup_{|x_i|>\sqrt{4a \rho^2 \log n}} \tilde{h}_2(x_i,\frac{1}{n},\eta, \delta) =o(1) \hspace*{0.05cm} \text{as} \hspace*{0.05cm} n \to \infty.
 \end{equation*}
 Therefore, we have,
	\begin{equation*}
		\sup_{|x_i|>\sqrt{4a \rho^2 \log n}} \tilde{h}(x_i,\frac{1}{n},\eta, \delta) \leq 1+\sup_{|x_i|>\sqrt{4a \rho^2 \log n}} \tilde{h}_1(x_i,\frac{1}{n})+\sup_{|x_i|>\sqrt{4a \rho^2 \log n}} \tilde{h}_2(x_i,\frac{1}{n},\eta, \delta) \lesssim 1 \hspace*{0.05cm} \cdot
	\end{equation*}
	Using above arguments, as $n \to \infty$,
	\begin{equation*}\label{eq:T-A1.4}
		\mathbb{E}_{{\theta}_{0i}}[Var(\theta_i|X_i,\widehat{\tau}) 1_{\{|X_i|> \sqrt{4a \rho^2 \log n}\}}] \lesssim 1\hspace{0.1cm} . \tag{6}
	\end{equation*}
	Combining \eqref{eq:T-A1.2},\eqref{eq:T-A1.3} and \eqref{eq:T-A1.4}, we obtain 
	\begin{equation*}\label{eq:T-A1.5}
		\mathbb{E}_{{\theta}_{0i}}  Var(\theta_i|X_i,\widehat{\tau}) \lesssim \log n \hspace{0.1cm} .\tag{7}
	\end{equation*}
 Noting that $\lesssim$ relation proved here actually holds uniformly in $i$'s such that  ${\theta}_{0i} \neq 0$, in that the corresponding constants in the upper bounds can be chosen to be the same for each $i$, we have
	%Since all of the above arguments hold true for any $i$ such that ${\theta}_{0i} \neq 0$, as $n \to \infty$,
	\begin{equation*}\label{eq:T-A1.6}
		\sum_{i:\theta_{0i} \neq 0} \mathbb{E}_{{\theta}_{0i}}  Var(\theta_i|X_i,\widehat{\tau}) \lesssim \tilde{q}_n \log n  \hspace{0.1cm} .\tag{8}
	\end{equation*}
\textbf{Proof of Step-2:} Fix any $i$ such that $\theta_{0i} = 0$. Define $v_n= \sqrt{c_1 \log n}$, where $c_1$ is defined in (5) of the main document. \\
\textbf{Case-1}
	First we consider the case when $a \in [\frac{1}{2},1)$. Again, we split $\mathbb{E}_{{\theta}_{0i}}  Var(\theta_i|X_i,\widehat{\tau})$ as
	\begin{equation*} \label{eq:T-A1.7}
		\mathbb{E}_{{\theta}_{0i}}  Var(\theta_i|X_i,\widehat{\tau})=\mathbb{E}_{{\theta}_{0i}} [ Var(\theta_i|X_i,\widehat{\tau})1_{\{|X_i|> v_n\}}]+\mathbb{E}_{{\theta}_{0i}} [ Var(\theta_i|X_i,\widehat{\tau})1_{\{|X_i|\leq v_n\}}] \hspace*{0.05cm} \cdot
  \tag{9}
	\end{equation*}
	Using $Var(\theta_i|x_i,\widehat{\tau}) \leq 1+x^2_i$ and the identity $x^2 \phi(x)=\phi(x)-\frac{d}{dx}[x \phi(x)]$, we obtain, 
	\begin{equation*} \label{eq:T-A1.8}
		\mathbb{E}_{{\theta}_{0i}} [ Var(\theta_i|X_i,\widehat{\tau})1_{\{|X_i|> v_n\}}] \leq 2 \int_{\sqrt{c_1 \log n}}^{\infty} (1+x^2)\phi(x) dx \lesssim \sqrt{\log n}\cdot n^{-\frac{c_1}{2}} \hspace*{0.1cm}. \tag{10}
	\end{equation*}
	Now let us choose some $\gamma>1$ such that $c_2\gamma-1>1$.
	Next, we decompose the second term as follows:
	\begin{align*} \label{eq:T-A1.9}
		\mathbb{E}_{{\theta}_{0i}}[Var(\theta_i|X_i,\widehat{\tau})1_{\{|X_i| \leq  v_n\}}] &= \mathbb{E}_{{\theta}_{0i}} [ Var(\theta_i|X_i,\widehat{\tau})
		1_{\{\widehat{\tau}>\gamma\frac{q_n}{n}\}}
		1_{\{|X_i| \leq  v_n\}}]+ \\
		\mathbb{E}_{{\theta}_{0i}} [ Var(\theta_i|X_i,\widehat{\tau})
		1_{\{\widehat{\tau}\leq \gamma\frac{q_n}{n}\}}
		1_{\{|X_i| \leq  v_n\}}] \hspace*{0.1cm}.\tag{11}
	\end{align*}
	Note that \begin{align*}
		\mathbb{E}_{{\theta}_{0i}} [ Var(\theta_i|X_i,\widehat{\tau})
		1_{\{\widehat{\tau}>\gamma\frac{q_n}{n}\}}
		1_{\{|X_i| \leq  v_n\}}] & \leq (1+c_1 \log n )\mathbb{P}_{\theta_0}[\widehat{\tau}>\gamma\frac{q_n}{n},|X_i| \leq  v_n]
	\end{align*}
 \begin{equation*} \label{eq:T-A1.10}
     \leq (1+c_1 \log n)\mathbb{P}_{\theta_0}[\frac{1}{c_2n}\sum_{j=1(\neq i)}^{n}1_{\{|X_j|> \sqrt{c_1 \log n}\}}>\gamma\frac{q_n}{n}]  \lesssim \frac{q_n}{n} \log n \hspace{0.1cm}. \tag{12}
\end{equation*}
 Inequality in the last line is due to employing similar arguments used for proving Lemma A.7 in \cite{van2014horseshoe}.\\
 %	Inequality in the last line is due to employing similar arguments used for proving Lemma A.7 in \cite{van2014horseshoe} we can show that,
	%\begin{equation*}
	%	\mathbb{P}_{\theta_0}[\frac{1}{c_2n}\sum_{j=1(\neq i)}^{n}1_{\{|X_j|> \sqrt{c_1 \log n}\}}>\gamma\frac{q_n}{n}] \lesssim \frac{q_n}{n}\hspace{0.1cm}.
	%\end{equation*}
	%Therefore we have
	%\begin{equation*} \label{eq:T-1.10}
	%	\mathbb{E}_{{\theta}_{0i}} [ Var(\theta_i|X_i,\widehat{\tau})	1_{\{\widehat{\tau}>\gamma\frac{q_n}{n}\}}	1_{\{|X_i| \leq  \sqrt{c_1 \log n}\}}] \lesssim \frac{q_n}{n} \log n \hspace{0.1cm}. \tag{2.18}
	%\end{equation*}
 We will now bound $\mathbb{E}_{\theta_{0i}}[ Var(\theta_i|X_i,\widehat{\tau})
	1_{\{\widehat{\tau}\leq \gamma\frac{q_n}{n}\}}
	1_{\{|X_i| \leq  \sqrt{c_1 \log n}\}}]$. Since for any fixed $x_i \in \mathbb{R}$  and $\tau>0$, $Var(\theta_i|x_i,\tau) \leq E(1-\kappa_i|x_i,\tau)+J(x_i,\tau)$ where $J(x_i,\tau)=x^2_iE[(1-\kappa_i)^2|x_i,\tau]$. Since $E(1-\kappa_i|x_i,\tau)$ is non-decreasing in $\tau$, so,
	$E(1-\kappa_i|x_i,\widehat{\tau}) \leq E(1-\kappa_i|x_i,\gamma\frac{q_n}{n} )$ whenever $\widehat{\tau}\leq \gamma\frac{q_n}{n}$. Using Lemma 2 and  A.2 of \cite{ghosh2017asymptotic}, 
	%\begin{equation*}
	%	E(1-\kappa_i|x_i,\gamma\frac{q_n}{n} ) \lesssim e^{\frac{x^2_i}{2}} (\frac{q_n}{n})^{2a}(1+o(1)),
	%\end{equation*}
	%where $o(1)$ depends only on $n$ such that $\lim_{n \to \infty}o(1)=0$ and $q_n \to \infty$ as $n \to \infty$ such that $q_n=o(n)$. Similarly, $J(x_i,\widehat{\tau}) \leq J(x_i,\gamma\frac{q_n}{n})$ whenever $\widehat{\tau}\leq \gamma\frac{q_n}{n}$. Using Lemma A.2 of \cite{ghosh2017asymptotic},
	%\begin{equation*}
	%	J(x_i,\gamma\frac{q_n}{n}) \lesssim (\frac{q_n}{n})^{2a}e^{\frac{x^2_i}{2}}(1+o(1)),
	%\end{equation*}
% where $o(1)$ depends only on $n$ such that $\lim_{n \to \infty}o(1)=0$ and $q_n \to \infty$ as $n \to \infty$ such that $q_n=o(n)$.
	%Using the above facts, we have,
	\begin{align*} \label{eq:T-1.11}
		\mathbb{E}_{\theta_{0i}}[ Var(\theta_i|X_i,\widehat{\tau})
		1_{\{\widehat{\tau}\leq \gamma\frac{q_n}{n}\}}
		1_{\{|X_i| \leq  \sqrt{c_1 \log n}\}}] & \lesssim (\frac{q_n}{n})^{2a} \int_{0}^{\sqrt{c_1 \log n}} e^{\frac{x^2}{2}} e^{-\frac{x^2}{2}} dx \\
		&= (\frac{q_n}{n})^{2a} \sqrt{c_1 \log n} \hspace{0.1cm}.\tag{13}
	\end{align*}
	Note that all these preceding arguments hold uniformly in $i$ such that $\theta_{0i}=0$. Combining all these results, for $a\in [\frac{1}{2},1)$ using \eqref{eq:T-A1.7}-\eqref{eq:T-1.11}, we have,
	\begin{align*}\label{eq:T-1.12}
		\sum_{i:\theta_{0i} = 0} \mathbb{E}_{{\theta}_{0i}}  Var(\theta_i|X_i,\widehat{\tau}) &\lesssim (n-\tilde{q}_n) [\sqrt{\log n}\cdot n^{-\frac{c_1}{2}}+\frac{q_n}{n}\cdot \log n+(\frac{q_n}{n})^{2a} \sqrt{\log n}] \\
		& \lesssim q_n \log n\hspace{0.1cm}.  \tag{14}
	\end{align*}
	The second inequality follows due to the fact that $\tilde{q}_n \leq q_n$ and $q_n=o(n)$ as $n \to \infty$.\\
	\textbf{Case-2} Now we assume $a \geq 1$ and split $\mathbb{E}_{{\theta}_{0i}}  Var(\theta_i|X_i,\widehat{\tau})$ as
	\begin{equation*} \label{eq:T-1.13}
		\mathbb{E}_{{\theta}_{0i}}  Var(\theta_i|X_i,\widehat{\tau})=\mathbb{E}_{{\theta}_{0i}} [ Var(\theta_i|X_i,\widehat{\tau})1_{\{|X_i|> {v_n}\}}]+\mathbb{E}_{{\theta}_{0i}} [ Var(\theta_i|X_i,\widehat{\tau})1_{\{|X_i|\leq {v_n }\}}] \tag{15},
	\end{equation*}
	where $v_n$ is the same as defined in \textbf{Case-1}.
 Using exactly the same arguments used for the case $a \in [\frac{1}{2},1)$, we have the following  for $a \geq 1$,
	\begin{equation*} \label{T-1.14}
		\mathbb{E}_{{\theta}_{0i}} [ Var(\theta_i|X_i,\widehat{\tau})1_{\{|X_i|> \sqrt{c_1  \log n}\}}] \lesssim \sqrt{\log n}\cdot n^{-\frac{c_1}{2}} \hspace*{0.1cm}. \tag{16}
	\end{equation*}
	and \begin{equation*} \label{T-1.15}
		\mathbb{E}_{{\theta}_{0i}} [ Var(\theta_i|X_i,\widehat{\tau})
		1_{\{\widehat{\tau}>\gamma\frac{q_n}{n}\}}
		1_{\{|X_i| \leq  \sqrt{c_1 \log n}\}}] \lesssim \frac{q_n}{n} \log n \hspace{0.1cm}. \tag{17}
	\end{equation*}
	For the part $	\mathbb{E}_{{\theta}_{0i}} [ Var(\theta_i|X_i,\widehat{\tau})
	1_{\{\widehat{\tau} \leq \gamma\frac{q_n}{n}\}}
	1_{\{|X_i| \leq  \sqrt{c_1 \log n}\}}]$, note that for fixed $x_i \in \mathbb{R}$ and any $\tau >0$,
	\begin{align*} \label{T-1.16}
		Var(\theta_i|x_i,\tau) & \leq E(1-\kappa_i|x_i,\tau)+x^2_i E[(1-\kappa_i)^2|x_i,\tau] \\
		& \leq E(1-\kappa_i|x_i,\tau)+x^2_i E(1-\kappa_i|x_i,\tau) \\
		& \leq E(1-\kappa_i|x_i,\tau) 1_{\{|x_i| \leq 1 \}} +2 x^2_i E(1-\kappa_i|x_i,\tau) \hspace{0.1cm}. \tag{18}
	\end{align*}

 Since for fixed $x_i \in \mathbb{R}$ and any $\tau >0$, $E(1-\kappa_i|x_i,\tau)$ is non-decreasing in $\tau$, we have, \begin{align*} \label{T-1.17}
		\mathbb{E}_{{\theta}_{0i}} [ Var(\theta_i|X_i,\widehat{\tau})
		1_{\{\widehat{\tau} \leq \gamma\frac{q_n}{n}\}}
		1_{\{|X_i| \leq  \sqrt{c_1 \log n}\}}] & \leq \mathbb{E}_{{\theta}_{0i}}[E(1-\kappa_i|X_i,\gamma \frac{q_n}{n})1_{\{|X_i| \leq 1 \}}] + \\
		2 \mathbb{E}_{{\theta}_{0i}}[X^2_i E(1-\kappa_i|X_i,\gamma \frac{q_n}{n})1_{\{|X_i| \leq \sqrt{c_1 \log n} \}}] \hspace{0.1cm}. \tag{19}
	\end{align*}
For bounding the first term in the r.h.s. of \eqref{T-1.17}, we use Lemma 1. We note that for any $\tau \in (0,1)$, 
	$\frac{t \tau^2}{1+t \tau^2} \cdot \frac{1}{\sqrt{1+t \tau^2}} t^{-a-1} \leq \tau t^{-(a+\frac{1}{2})}$ and  that $L(t)$ is bounded.  Using the fact that the second term in the upper bound in Lemma 1 can be bounded 
	\begin{align*} 
		A_2 & \leq \frac{2M}{(2a-1)} \tau e^{\frac{x^2_i}{2}} \hspace*{0.1cm},
	\end{align*}
where $A_2$ has been defined in the proof of Lemma 1 of the main document.
	So, we have
	\begin{align*}
		\mathbb{E}_{{\theta}_{0i}}[E(1-\kappa_i|X_i,\gamma \frac{q_n}{n})1_{\{|X_i| \leq 1 \}}] & \lesssim \frac{q_n}{n} \int_{0}^{1} e^{-\frac{x^2}{4}} dx+ \frac{q_n}{n} \hspace{0.1cm}. 
	\end{align*}
	Hence, \begin{equation*} \label{T-1.18}
		\mathbb{E}_{{\theta}_{0i}}[E(1-\kappa_i|X_i,\gamma \frac{q_n}{n})1_{\{|X_i| \leq 1 \}}] \lesssim \frac{q_n}{n} \hspace{0.1cm}.  \tag{20}
	\end{equation*}
	For the second term in the r.h.s. of \eqref{T-1.17} we shall use the upper bound of  $E(1-\kappa_i|x_i,\tau)$ of the form (7) of the main article and hence,
	\begin{equation*}
		\mathbb{E}_{{\theta}_{0i}}[X^2_i E(1-\kappa_i|X_i,\gamma \frac{q_n}{n})1_{\{|X_i| \leq \sqrt{c_1 \log n} \}}]  \lesssim  \frac{q_n}{n}  \int_{0}^{\sqrt{c_1 \log n}} x^2  e^{\frac{x^2}{4}}\phi(x)dx 
	\end{equation*}

	\begin{align*} \label{T-1.19}
		&+  \int_{0}^{\sqrt{c_1 \log n}} \int_{1}^{\infty}  \frac{t ({\gamma \frac{q_n}{n}})^2}{1+t ({\gamma \frac{q_n}{n}})^2} \cdot \frac{1}{\sqrt{1+t ({\gamma \frac{q_n}{n}})^2}} t^{-a-1} L(t) e^{\frac{x^2}{2}\cdot \frac{t ({\gamma \frac{q_n}{n}})^2}{1+t ({\gamma \frac{q_n}{n}})^2}}  x^2 \phi(x) dt dx \hspace{0.1cm}.  \tag{21}
	\end{align*}
	Note that the first integral is bounded by a constant. Using Fubini's theorem and the transformation $y=\frac{x}{\sqrt{1+t ({\gamma \frac{q_n}{n}})^2}}$, the second integral becomes
	\begin{equation*}
	\frac{1}{\sqrt{2 \pi}}	\int_{1}^{\infty} t^{-a-1} L(t) t ({\gamma \frac{q_n}{n}})^2  \bigg( \int_{0}^{\sqrt{\frac{c_1 \log n}{1+t ({\gamma \frac{q_n}{n}})^2}}} y^2 e^{-\frac{y^2}{2}} dy \bigg) dt \hspace{0.1cm}.  
	\end{equation*}
We handle the above integral separately for $a=1$ and $a>1$. For $a>1$, using the boundedness of $L(t)$ it is easy to show that, 
	\begin{equation*} \label{T-1.20}
	\frac{1}{\sqrt{2 \pi}}	\int_{1}^{\infty} t^{-a-1} L(t) t ({\gamma \frac{q_n}{n}})^2  \bigg( \int_{0}^{\sqrt{\frac{c_1 \log n}{1+t ({\gamma \frac{q_n}{n}})^2}}} y^2 e^{-\frac{y^2}{2}} dy \bigg) dt \lesssim \frac{q_n}{n} \hspace{0.1cm}.  \tag{22}
	\end{equation*}
	For $a=1$ note that,
 \begin{equation*}
     \int_{1}^{\infty} t^{-a-1} L(t) t ({\gamma \frac{q_n}{n}})^2  \bigg( \int_{0}^{\sqrt{\frac{c_1 \log n}{1+t ({\gamma \frac{q_n}{n}})^2}}} y^2 e^{-\frac{y^2}{2}} dy \bigg) dt 
 \end{equation*}
 %\begin{equation*}
 %    \leq 		\int_{1}^{\infty} t^{-a-1} L(t) t ({\gamma \frac{q_n}{n}})^2  \bigg( \int_{0}^{\sqrt{\frac{c_1 \log n}{t ({\gamma \frac{q_n}{n}})^2}}} y^2 e^{-\frac{y^2}{2}} dy \bigg) dt
% \end{equation*}
	\begin{align*} \label{T-1.21}
		\leq ({\gamma \frac{q_n}{n}})^2 \sqrt{2 \pi} \int_{1}^{\frac{c_1 \log n}{({\gamma \frac{q_n}{n}})^2 (2 \pi)^{\frac{1}{3}}}} \frac{1}{t}L(t) dt+  (\sqrt{c_1 \log n})^3 ({\gamma \frac{q_n}{n}})^2  \int_{\frac{c_1 \log n}{({\gamma \frac{q_n}{n}})^2 (2 \pi)^{\frac{1}{3}}}}^{\infty} \frac{L(t)}{t} \cdot \frac{1}{( t ({\gamma \frac{q_n}{n}})^2 )^{\frac{3}{2}}} dt \hspace{0.1cm}.  \tag{23}
	\end{align*}
	Here the division in the range of $t$ in \eqref{T-1.21} occurs due to the fact the integral \\ $( \int_{0}^{\sqrt{\frac{c_1 \log n}{t ({\gamma \frac{q_n}{n}})^2}}} y^2 e^{-\frac{y^2}{2}} dy)$ can be bounded by $(\frac{c_1 \log n}{t ({\gamma \frac{q_n}{n}})^2})^{\frac{3}{2}}$ when $t \geq \frac{c_1 \log n}{({\gamma \frac{q_n}{n}})^2 (2 \pi)^{\frac{1}{3}}}$ and by $\sqrt{2 \pi}$ when $t \leq \frac{c_1 \log n}{({\gamma \frac{q_n}{n}})^2 (2 \pi)^{\frac{1}{3}}}$. For the first term in \eqref{T-1.21} with the boundedness of $L(t)$,
	\begin{align*}
		({\gamma \frac{q_n}{n}})^2 \int_{1}^{\frac{c_1 \log n}{({\gamma \frac{q_n}{n}})^2 (2 \pi)^{\frac{1}{3}}}} \frac{1}{t}L(t) dt & \leq ({\gamma \frac{q_n}{n}})^2 M \log \bigg(\frac{c_1 \log n}{({\gamma \frac{q_n}{n}})^2 (2 \pi)^{\frac{1}{3}}} \bigg) \hspace{0.1cm}. 
	\end{align*}
	Hence for sufficiently large $n$ with $q_n \propto n^{\beta}, 0<\beta<1$,
	\begin{equation*} \label{T-1.22}
		({\gamma \frac{q_n}{n}})^2 \int_{1}^{ \frac{c_1 \log n}{({\gamma \frac{q_n}{n}})^2 (2 \pi)^{\frac{1}{3}}}} \frac{1}{t}L(t) dt \lesssim \frac{q_n}{n} \sqrt{\log n} \hspace{0.1cm}.  \tag{24}
	\end{equation*}
	Now for the second term in \eqref{T-1.21} again using the boundedness of $L(t)$,
	\begin{equation*} \label{T-1.23}
		(\sqrt{c_1 \log n})^3 ({\gamma \frac{q_n}{n}})^2  \int_{\frac{c_1 \log n}{({\gamma \frac{q_n}{n}})^2 (2 \pi)^{\frac{1}{3}}}}^{\infty} t \cdot t^{-2} L(t) \cdot \frac{1}{( t ({\gamma \frac{q_n}{n}})^2 )^{\frac{3}{2}}} dt \lesssim \frac{q_n}{n} \hspace{0.1cm}.  \tag{25}
	\end{equation*}
	So, using \eqref{T-1.21}-\eqref{T-1.23}, for $a=1$,
	\begin{equation*}\label{T-1.24}
		\int_{1}^{\infty} t^{-a-1} L(t) t ({\gamma \frac{q_n}{n}})^2  \bigg( \int_{0}^{\sqrt{\frac{c_1 \log n}{1+t ({\gamma \frac{q_n}{n}})^2}}} y^2 e^{-\frac{y^2}{2}} dy \bigg) dt  \lesssim \frac{q_n}{n} \sqrt{\log n} \hspace{0.1cm}.  \tag{26}
	\end{equation*}
 With the help of \eqref{T-1.20} and \eqref{T-1.24}, we have, for $a \geq 1$
	\begin{equation*} \label{T-1.25}
		\int_{1}^{\infty} t^{-a-1} L(t) t ({\gamma \frac{q_n}{n}})^2  \bigg( \int_{0}^{\sqrt{\frac{c_1 \log n}{1+t ({\gamma \frac{q_n}{n}})^2}}} y^2 e^{-\frac{y^2}{2}} dy \bigg) dt  \lesssim \frac{q_n}{n} \sqrt{\log n} \hspace{0.1cm}.  \tag{27}
	\end{equation*}
	Combining these facts, we finally have
	\begin{equation*} \label{T-1.26}
		\mathbb{E}_{{\theta}_{0i}}[X^2_i E(1-\kappa_i|X_i,\gamma \frac{q_n}{n})1_{\{|X_i| \leq \sqrt{c_1 \log n} \}}]  \lesssim  \frac{q_n}{n} \sqrt{\log n} \hspace{0.1cm}.  \tag{28}
	\end{equation*}
	Note that all these preceding arguments hold uniformly in $i$ such that $\theta_{0i}=0$. Combining all these results, for $a\geq 1$, using \eqref{eq:T-1.13}-\eqref{T-1.26}, we have,
	\begin{align*}\label{eq:T-1.27}
		\sum_{i:\theta_{0i} = 0} \mathbb{E}_{{\theta}_{0i}}  Var(\theta_i|X_i,\widehat{\tau}) &\lesssim (n-\tilde{q}_n) [\sqrt{\log n}\cdot n^{- \frac{c_1}{2}}+\frac{q_n}{n}\cdot \log n+\frac{q_n}{n} \sqrt{\log n}] \\
		& \lesssim q_n \log n\hspace{0.1cm}.  \tag{29}
	\end{align*}
	Using \eqref{eq:T-1.12} and \eqref{eq:T-1.27}, for $a \geq \frac{1}{2}$, we have,
	\begin{align*}\label{eq:T-1.28}
		\sum_{i:\theta_{0i} = 0} \mathbb{E}_{{\theta}_{0i}}  Var(\theta_i|X_i,\widehat{\tau}) &\lesssim  q_n \log n\hspace{0.1cm}.  \tag{30}
	\end{align*}
	Now using \eqref{eq:T-A1.1}, \eqref{eq:T-A1.6} and \eqref{eq:T-1.28}, for sufficiently large $n$,
	\begin{equation*}
		\mathbb{E}_{\boldsymbol{\theta}_0} \sum_{i=1}^{n} Var(\theta_i|X_i,\widehat{\tau}) \lesssim q_n \log n\hspace{0.1cm}.
	\end{equation*}
	Finally, taking supremum over all $\boldsymbol{\theta}_0 \in l_0[q_n]$, the result is obtained.	

\end{proof}

\begin{remark}
\label{rem-4}
 Note that, we have used different bounds on $Var(\theta_i|X_i, \tau)$ for two different ranges of $a$. When $a \in [\frac{1}{2},1)$, we have used the upper bound of $J(X_i, \tau)$ provided by \cite{ghosh2017asymptotic}. However using the same arguments when $a \geq 1$ yields an upper bound on $Var(\theta_i|X_i, \tau)$ such that 
 $	\mathbb{E}_{{\theta}_{0i}} [ Var(\theta_i|X_i,\widehat{\tau})
	1_{\{\widehat{\tau} \leq \gamma\frac{q_n}{n}\}}
	1_{\{|X_i| \leq  \sqrt{c_1 \log n}\}}]$ exceeds near minimax rate. As a result, we need to come up with a sharper upper bound, and hence \eqref{T-1.16} and Lemma 1 of the main document come in very handy.
\end{remark}

\subsection{Results related to the contraction rate of MMLE}
Recall that the logarithm of the marginal likelihood function $M_{\tau}(\mathbf{X})$ is of the form
\begin{align*}
    M_{\tau}(\mathbf{X}) &= \sum_{i=1}^{n} \log \bigg(\phi(x_i-\theta_i) g_{\tau}(\theta_i) d \theta_i \bigg), 
\end{align*}
where $g_{\tau}(\theta_i)$ is defined in the main paper.\\
\textbf{Main Idea of the Proof:-} In order to prove Lemma 2 of the main document, we are interested in finding an upper bound to the derivative of the logarithm of the marginal likelihood function $M_{\tau}(\mathbf{X})$. This is done with the help of a series of Lemmas stated below. In Lemma \ref{lem4}, we express the derivative as a sum of suitably defined i.i.d. random variables. In Lemmas \ref{lem5} and \ref{lem6}, we study the bounds of some quantities related to that derivative, which are defined in \eqref{eq:L-1.3} and \eqref{eq:L-1.4}. Lemma \ref{lem7} is about some bounds on those suitably defined i.i.d. random variables. Some bounds based on the moments of those random variables are provided in Lemmas \ref{lem8} and \ref{lem9}. Then, these Lemmas are used in order to obtain some corollaries, which are of similar form of those of \cite{van2017adaptive}. Finally, these Lemmas and corollaries are used for deriving Lemma 2 of the main document. We state and provide proof of the above-mentioned lemmas in the next subsection.

\begin{proof}[Proof of Lemma 2 of the main document]
    Using the definition of $M_{\tau}(\mathbf{X})$ along with \eqref{eq:L-1.1}, the derivative of the log-likelihood function is of the form
    \begin{align*} 
        \frac{d}{d \tau} M_{\tau}(\mathbf{X}) &= \frac{1}{\tau}  \sum_{i \in I_0} m_{\tau}(X_i)+
        \frac{1}{\tau}  \sum_{i \in I_1} m_{\tau}(X_i),
    \end{align*}
    where $I_0 := \{i: \theta_{0,i}=0 \}$ and $I_1 := \{i: \theta_{0,i} \neq 0 \}$. Next, for zero means, using Corollary \ref{cor3}, we see that the sum behaves like its expectation, uniformly in $\tau$. As $\tau \to 0$, the expression for $E_{\theta} m_{\tau}(X)$ is derived in  Lemma \ref{lem9}. On the other hand, for non-zero means, the quantity $m_{\tau}(x)$ is uniformly bounded in $\tau$, as proved in  (i) of Lemma \ref{lem7}. Using all of these arguments, we can show that, as $\tau \to 0$, the derivative of the log-likelihood function is bounded as
    \begin{align*} \label{eq:L-7.2}
        \frac{d}{d \tau} M_{\tau}(\mathbf{X}) & < -\frac{2K(n-q_n)}{\sqrt{2 \pi}\zeta_{\tau}} (1+o_P(1))+\frac{q_n}{\tau} \tag{31}.
    \end{align*}
    For $\frac{\tau}{\zeta_{\tau}} \gtrsim \frac{q_n}{(n-q_n)}$, the upper bound as obtained in \eqref{eq:L-7.2} is negative, and hence $\widehat{\tau}_n$ satisfies $\frac{\widehat{\tau}_n}{\zeta_{\widehat{\tau}_n}} \lesssim \frac{q_n}{(n-q_n)}$, i.e., $\widehat{\tau}_n \lesssim \tau_n(q_n) $, under the assumption $q_n=o(n)$.
\end{proof}

\begin{proof}[Proof of Theorem \ref{Thm1}]
    For proving Theorem 2, with the use of ({C1}) given in the main document, note that
    \begin{align*} \label{eq:T-1.1}
        &\mathbb{E}_{\boldsymbol{\theta}_0} \Pi_{\widehat{\tau}_n} (  \boldsymbol{\theta}:||\boldsymbol{\theta}-\boldsymbol{\theta}_0 ||^2>    M_n q_n \log n |\mathbf{X} \bigg) \\
        &= \mathbb{E}_{\boldsymbol{\theta}_0} \bigg[\Pi_{\widehat{\tau}_n} (  \boldsymbol{\theta}:||\boldsymbol{\theta}-\boldsymbol{\theta}_0 ||^2>    M_n q_n \log n |\mathbf{X} ) 1_{ \{\widehat{\tau}_n \in [\frac{1}{n}, C \tau_n(q_n)] \}} \bigg] \\ &+\mathbb{E}_{\boldsymbol{\theta}_0} \bigg[\Pi_{\widehat{\tau}_n} (  \boldsymbol{\theta}:||\boldsymbol{\theta}-\boldsymbol{\theta}_0 ||^2>    M_n q_n \log n |\mathbf{X} ) 1_{ \{\widehat{\tau}_n \notin [\frac{1}{n}, C \tau_n(q_n)] \}} \bigg] \\
        & \leq \mathbb{E}_{\boldsymbol{\theta}_0} \bigg[\Pi_{\widehat{\tau}_n} (  \boldsymbol{\theta}:||\boldsymbol{\theta}-\boldsymbol{\theta}_0 ||^2>    M_n q_n \log n |\mathbf{X} ) 1_{ \{\widehat{\tau}_n \in [\frac{1}{n}, C \tau_n(q_n)] \}} \bigg] +o(1) \tag{32}.
    \end{align*}
    Next, applying Markov's inequality to the first term in the r.h.s. of \eqref{eq:T-1.1}, it is sufficient to show that
    \begin{align*} \label{eq:T-1.2}
        \mathbb{E}_{\theta_0} \sum_{i=1}^{n} Var(\theta_i|X_i,\widehat{\tau}_n)1_{ \{\widehat{\tau}_n \in [\frac{1}{n}, C \tau_n(q_n)] \}} &\lesssim q_n \log n \tag{33}
    \end{align*}
    and
    \begin{align*} \label{eq:T-1.3}
        \mathbb{E}_{\theta_0} ||T_{\widehat{\tau}_n}(\mathbf{X}) -\boldsymbol{\theta}_0||^2 1_{ \{\widehat{\tau}_n \in [\frac{1}{n}, C \tau_n(q_n)] \}} &\lesssim q_n \log n .\tag{34}
    \end{align*}
    In order to prove \eqref{eq:T-1.2}, we need to show that
    \begin{align*} \label{eq:T-1.4}
        \sum_{i:\theta_{0i} \neq 0} \mathbb{E}_{{\theta}_{0i}}  Var(\theta_i|X_i,\widehat{\tau}_n)1_{ \{\widehat{\tau}_n \in [\frac{1}{n}, C \tau_n(q_n)] \}} &\lesssim \tilde{q}_n \log n \tag{35}
    \end{align*}
    and
    \begin{align*} \label{eq:T-1.5}
       \sum_{i:\theta_{0i} = 0} \mathbb{E}_{{\theta}_{0i}}  Var(\theta_i|X_i,\widehat{\tau}_n)1_{ \{\widehat{\tau}_n \in [\frac{1}{n}, C \tau_n(q_n)] \}} &\lesssim  q_n \log n . \tag{36}
    \end{align*}
\textbf{Step-1} Fix any $i$ such that $\theta_{0i} \neq 0$. Note that while proving the contraction rate of the empirical Bayes estimator of $\tau$ as proposed by \cite{van2014horseshoe}, we need to prove a statement same as that of \eqref{eq:T-1.4}. In that case, we used only that $\widehat{\tau} \geq \frac{1}{n}$ and some monotonicity of the shrinkage coefficient, which is true in this case, too. Hence, using the same arguments as used in Theorem 1 related to the non-null means, we can show \eqref{eq:T-1.4} holds.\\
\textbf{Step-2} Fix any $i$ such that $\theta_{0i} = 0$. Next, moving towards proving \eqref{eq:T-1.5}, using \eqref{eq:T-1.11} of Theorem 1 with $v_n=\sqrt{2 \log n}$, we have,
    \begin{align*} \label{eq:T-1.6}
        &\mathbb{E}_{{\theta}_{0i}}  Var(\theta_i|X_i,\widehat{\tau}_n)1_{ \{\widehat{\tau}_n \in [\frac{1}{n}, C \tau_n(q_n)] \}} 1_{ \{|X_i| \leq v_n\} } \\
        & \lesssim \tau_n(q_n) \int_{0}^{\sqrt{2 \log n}} e^{\frac{x^2}{2}} e^{-\frac{x^2}{2}} dx \\
        & \lesssim \frac{q_n}{n} \sqrt{\log (\frac{n}{q_n})} \sqrt{\log n} \\
        & \lesssim \frac{q_n}{n} \log n. \tag{37}
    \end{align*}
    On the other hand, using \eqref{eq:T-A1.8} of Theorem 1,
    \begin{align*} \label{eq:T-1.7}
         &\mathbb{E}_{{\theta}_{0i}}  Var(\theta_i|X_i,\widehat{\tau}_n)1_{ \{\widehat{\tau}_n \in [\frac{1}{n}, C \tau_n(q_n)] \}} 1_{ \{|X_i| > v_n\} } \\
         & \leq \mathbb{E}_{{\theta}_{0i}} [ Var(\theta_i|X_i,\widehat{\tau}_n)1_{\{|X_i|> v_n\}}] \\
         & \lesssim \frac{\sqrt{\log n}}{n}. \tag{38}
    \end{align*}
    Combining \eqref{eq:T-1.6} and \eqref{eq:T-1.7}, we can conclude that
    \begin{align*}
        \sum_{i:\theta_{0i} = 0} \mathbb{E}_{{\theta}_{0i}}  Var(\theta_i|X_i,\widehat{\tau}_n)1_{ \{\widehat{\tau}_n \in [\frac{1}{n}, C \tau_n(q_n)] \}} &\lesssim (n-\tilde{q}_n) [\frac{\sqrt{\log n}}{n}+\frac{q_n}{n}\cdot \log n] \\
		& \lesssim q_n \log n\hspace{0.1cm}.
    \end{align*}
    This completes the proof of \eqref{eq:T-1.5} and eventually \eqref{eq:T-1.2} is also established. \\
    Next, for \eqref{eq:T-1.3}, we aim to prove that
    \begin{align*} \label{eq:T-1.8}
        \sum_{i:\theta_{0i} \neq 0} \mathbb{E}_{{\theta}_{0i}} (T_{\widehat{\tau}_n}({X}_i) -\boldsymbol{\theta}_{0i})^2 1_{ \{\widehat{\tau}_n \in [\frac{1}{n}, C \tau_n(q_n)] \}} &\lesssim q_n \log n \tag{39}
    \end{align*}
    and
    \begin{align*} \label{eq:T-1.9}
        \sum_{i:\theta_{0i} = 0} \mathbb{E}_{{\theta}_{0i}} (T_{\widehat{\tau}_n}({X}_i) -\boldsymbol{\theta}_{0i})^2 1_{ \{\widehat{\tau}_n \in [\frac{1}{n}, C \tau_n(q_n)] \}} &\lesssim q_n \log n .\tag{40}
    \end{align*}
\textbf{Step-1} Fix any $i$ such that $\theta_{0i} \neq 0$. Again using the same set of arguments as used before, with Theorem 2 of \cite{ghosh2017asymptotic} for non-zero means, \eqref{eq:T-1.8} follows readily.\\
  \textbf{Step-2} Fix any $i$ such that $\theta_{0i} = 0$. Note that,
    \begin{align*} \label{eq:T-1.10}
       & \mathbb{E}_{{\theta}_{0i}} (T_{\widehat{\tau}_n}({X}_i) -\boldsymbol{\theta}_{0i})^2 1_{ \{\widehat{\tau}_n \in [\frac{1}{n}, C \tau_n(q_n)] \}}
       \leq \mathbb{E}_{{\theta}_{0i}}[T^2_{\tau_n(q_n)}(X_i)] \\
       &= \mathbb{E}_{{\theta}_{0i}}[T^2_{\tau_n(q_n)}(X_i)1_{\{|X_i| > \sqrt{2 \log (\frac{1}{\tau_n(q_n)})}\}}] + \mathbb{E}_{{\theta}_{0i}}[T^2_{\tau_n(q_n)}(X_i)1_{\{|X_i| \leq \sqrt{2 \log (\frac{1}{\tau_n(q_n)})}\}}] \\
       & \leq  \mathbb{E}_{{\theta}_{0i}}[X^2_i 1_{\{|X_i| > \sqrt{2 \log (\frac{1}{\tau_n(q_n)})}\}}] + \mathbb{E}_{{\theta}_{0i}}[T^2_{\tau_n(q_n)}(X_i)1_{\{|X_i| \leq \sqrt{2 \log (\frac{1}{\tau_n(q_n)})}\}}] \\
       & \lesssim e^{-\log (\frac{1}{\tau_n(q_n)})} \sqrt{\log (\frac{n}{q_n})}+ (\tau_n(q_n))^{2} \sqrt{\log (\frac{n}{q_n})} e^{ \log (\frac{1}{\tau_n(q_n)})} \\
       & \lesssim \frac{q_n}{n} \log (\frac{n}{q_n}) ,  \tag{41}
    \end{align*}
    where we use facts similar to (15) and (43) of \cite{ghosh2017asymptotic}. As a result,
    \begin{align*}
        \sum_{i:\theta_{0i} = 0} \mathbb{E}_{{\theta}_{0i}} (T_{\widehat{\tau}_n}({X}_i) -\boldsymbol{\theta}_{0i})^2 1_{ \{\widehat{\tau}_n \in [\frac{1}{n}, C \tau_n(q_n)] \}} &\lesssim (n-\Tilde{q}_n) \frac{q_n}{n} \log (\frac{n}{q_n}) \\
        & \lesssim q_n \log n.
    \end{align*}
    This completes the proof of \eqref{eq:T-1.9} and as a consequence of this, \eqref{eq:T-1.3} is also established. Combining \eqref{eq:T-1.1}-\eqref{eq:T-1.3} completes the proof of Theorem 2.
  % The proof of the first assertion is completed.  In the case of the second assertion, following the same steps mentioned before, we only need to show that, 
%    \begin{align*}
  %     E_{\theta_0} \sum_{i=1}^{n} Var(\theta_i|X_i,\widehat{\tau}_n)1_{ \{\widehat{\tau}_n \in [\frac{1}{n}, C \tau_n(q_n)] \}} &\lesssim q_n \log n, 
%    \end{align*}
  %  which is already proved in the first assertion.
\end{proof}

\subsection{Lemmas related to MMLE}

\begin{lemma}
\label{lem4}
The derivative of the log-likelihood function is of the form
\begin{equation*} \label{eq:L-1.1}
    \frac{d}{d \tau} M_{\tau}(\mathbf{X})= \frac{1}{\tau}\sum_{i=1}^{n}m_{\tau}(x_i) \tag{42},
\end{equation*}
where 
\begin{equation*}  \label{eq:L-1.2}
    m_{\tau}(x)= x^2 \frac{J_{\alpha+1,a}(x)-J_{\alpha+2,a}(x)}{I_{\alpha,a}(x)}- \frac{J_{\alpha+1,a}(x)}{I_{\alpha,a}(x)} \tag{43},
\end{equation*}
and $I_{\alpha,a}(x)$ and for $k=1,2, J_{\alpha+k,a}(x)$ are defined as follows
\begin{equation*} \label{eq:L-1.3}
    I_{\alpha,a}(x) = \int_{0}^{1} e^{\frac{x^2z}{2}} z^{\alpha-1} (1-z)^{a-\frac{1}{2}} \bigg(\frac{1}{\tau^2 +(1-\tau^2)z}\bigg)^{a+\alpha} d z \tag{44}
\end{equation*}
and
\begin{equation*}  \label{eq:L-1.4}
    J_{\alpha+k,a}(x) = \int_{0}^{1} e^{\frac{x^2z}{2}} z^{\alpha+k-1} (1-z)^{a-\frac{1}{2}} \bigg(\frac{1}{\tau^2 +(1-\tau^2)z}\bigg)^{a+\alpha} d z \tag{45}.
\end{equation*}
\end{lemma}
\begin{proof}
    Since, $X_i |\theta_i \simind \mathcal{N}(\theta_i,1)$ and $\theta_i|\lambda_i,\tau \simind \mathcal{N}(0,\lambda^2_i\tau^2)$, the marginal density of $X_i$ given $\tau$ is of this form
    \begin{align*}
        \psi_{\tau}(x) &= K \int_{0}^{\infty} \frac{e^{-\frac{1}{2} \frac{x^2}{(1+\lambda^2 \tau^2)}}}{\sqrt{1+\lambda^2 \tau^2} \sqrt{2 \pi}} (\lambda^2)^{-a-1}L(\lambda^2) d \lambda^2 \\
        &= K \frac{\tau^{2a}}{\sqrt{2 \pi}} \int_{0}^{1} e^{-\frac{x^2(1-z)}{2}} z^{-a-1} (1-z)^{a-\frac{1}{2}} L(\frac{1}{\tau^2} \frac{z}{1-z}) dz,
    \end{align*}
    where equality in the second line follows due to the substitution $1-z=(1+\lambda^2 \tau^2)^{-1}$. Next, noting that $L(t)=(1+\frac{1}{t})^{-(a+\alpha)}$, 
    \begin{align*} \label{eq:1}
        \psi_{\tau}(x) &= K \tau^{2a}  \int_{0}^{1} \frac{e^{-\frac{x^2(1-z)}{2}}} {\sqrt{2\pi}} z^{\alpha-1} (1-z)^{a-\frac{1}{2}} [\frac{1}{\tau^2+(1-\tau^2)z}]^{a+\alpha} dz \\
        &= K \tau^{2a} I_{\alpha,a}(x) \phi(x) \tag{46},
    \end{align*}
    where $\phi(x)$ denotes the standard normal density. Using \eqref{eq:1},
    \begin{align*} \label{eq:2}
        \frac{\dot{\psi}_{\tau}}{\psi_{\tau}} &= \frac{2a \tau^{2a-1}I_{\alpha,a}(x) +\tau^{2a} \dot{I}_{{\alpha,a}}(x)}{\tau^{2a}I_{\alpha,a}(x)} \\
        &=  \frac{2a I_{\alpha,a}(x) +\tau \dot{I}_{{\alpha,a}}(x)}{\tau I_{\alpha,a}(x)}  \tag{47} \\
        &= \frac{ \int_{0}^{1} e^{\frac{x^2z}{2}} z^{\alpha-1} (1-z)^{a-\frac{1}{2}} (\frac{1}{N(z)})^{a+\alpha+1}[2a N(z)-2\tau^2(\alpha+a)(1-z)] dz}{\tau I_{\alpha,a}(x)} ,
    \end{align*}
    where $N(z)=\tau^2+(1-\tau^2)z$.  On the other hand, using integration by parts,
    \begin{align*}
        & x^2 ({J_{\alpha+1,a}(x)-J_{\alpha+2,a}(x)}) = \int_{0}^{1} x^2 e^{\frac{x^2z}{2}} z^{\alpha} (1-z)^{a+\frac{1}{2}} (\frac{1}{N(z)})^{a+\alpha} dz \\
      %  &= -2 \int_{0}^{1} e^{\frac{x^2z}{2}}  d[\frac{z^{\alpha}(1-z)^{a+\frac{1}{2}}}{(N(z))^{\alpha+a}}]  \\
       & = \int_{0}^{1} e^{\frac{x^2z}{2}} \frac{ z^{\alpha-1} (1-z)^{a-\frac{1}{2}}}{(N(z))^{a+\alpha+1}} [N(z) \{-2 \alpha (1-z) +2(a+\frac{1}{2})z \} +2 z(1-z) (1-\tau^2) (\alpha+a)] dz
    \end{align*}
    As a consequence of this,
    \begin{align*} \label{eq:3}
     & x^2 ({J_{\alpha+1,a}(x)-J_{\alpha+2,a}(x)})- J_{\alpha+1,a}(x) \\
     &= \int_{0}^{1} e^{\frac{x^2z}{2}} \frac{z^{\alpha-1} (1-z)^{a-\frac{1}{2}}}{(N(z))^{a+\alpha+1}} [N(z) \{-2 \alpha (1-z) +2az\} +2(\alpha+a)(1-z)(N(z)-\tau^2)] dz \tag{48}.
    \end{align*}
    On Simplification the r.h.s. of \eqref{eq:3} matches with the numerator of \eqref{eq:2} and completes the proof.
\end{proof}
\begin{lemma}
\label{lem5} Let $\kappa_{\tau}$ be the solution to the equation $\frac{e^{x^2/2}}{x^2/2}=\frac{1}{\tau}$.  Choose any $B \geq 1$.
    There exist functions $R_{\tau}$ with $\sup_{x}|R_{\tau}(x)|={O}(\tau^{\frac{1}{2}})$ as $\tau \to 0$, such that, for $\alpha \geq \frac{1}{2}$,
    \begin{align*} \label{eq:L-2.1}
        I_{\alpha,\frac{1}{2}}(x) &= \begin{cases}
            \big(\frac{K^{-1}}{\tau} +\bigg(\frac{x^2}{2}\bigg)^{\frac{1}{2}} \int_{1}^{\frac{x^2}{2}}e^{v} v^{-\frac{3}{2}} dv \big) (1+R_{\tau}(x)), \text{uniformly} \hspace*{0.2cm} \text{in} \hspace*{0.2cm} |x| \leq  B\kappa_\tau ,\tag{49} \\ 
            \bigg(\frac{x^2}{2}\bigg)^{\frac{1}{2}} \int_{1}^{\frac{x^2}{2}}e^{v} v^{-\frac{3}{2}} dv  (1+R_{\tau}(x)) , \text{uniformly} \hspace*{0.2cm} \text{in} \hspace*{0.2cm} |x| > B\kappa_\tau .
        \end{cases}
    \end{align*}
 %   and
%     \begin{equation*} \label{eq:L-2.2}
%        I_{\alpha,\frac{1}{2}}(x)= \bigg(\frac{x^2}{2}\bigg)^{\frac{1}{2}} \int_{1}^{\frac{x^2}{2}}e^{v} v^{-\frac{3}{2}} dv  (1+R_{\tau}(x)) , \text{uniformly} \hspace*{0.2cm} \text{in} \hspace*{0.2cm} |x| > B\kappa_\tau . \tag{39}
%    \end{equation*}
   Further, given $\epsilon_{\tau} \to 0$, there exist functions $S_{\tau}$ with $\sup_{x\geq  1/\epsilon_{\tau}}|S_{\tau}(x)|={O}(\tau^{\frac{1}{2}}+\epsilon^2_{\tau})$, such that $\tau \to 0$,
   \begin{equation*} \label{eq:L-2.2}
     I_{\alpha,\frac{1}{2}}(x)=  \frac{e^{x^2/2}}{x^2/2} (1+S_{\tau}(x)) \tag{50},
   \end{equation*}
    where $K=\frac{\Gamma(\frac{1}{2}+\alpha)}{\sqrt{\pi}\Gamma(\alpha)} \hspace*{0.05cm} \cdot$
\end{lemma}
\begin{proof}
    We consider the cases separately when $|x| \leq B \kappa_{\tau}$ and $|x| > B\kappa_{\tau}$ where $\kappa_{\tau} \sim \zeta_{\tau}+\frac{2 \log \zeta_{\tau}}{\zeta_{\tau}}$ and
    $\zeta_{\tau}=\sqrt{2 \log (\frac{1}{\tau})}$.\\
    \textbf{Case 1}:- When $|x| \leq B \kappa_{\tau}$, the range of the integration in $I_{\alpha,\frac{1}{2}}(x)$ is divided into three parts, namely 
    \begin{equation*}
        I_1= \int_{0}^{\tau} e^{\frac{x^2z}{2}} z^{\alpha-1}  \bigg(\frac{1}{\tau^2 +(1-\tau^2)z}\bigg)^{\frac{1}{2}+\alpha} d z,
    \end{equation*}
    \begin{equation*}
        I_2= \int_{\tau}^{(\frac{2}{x^2})\wedge1} e^{\frac{x^2z}{2}} z^{\alpha-1}  \bigg(\frac{1}{\tau^2 +(1-\tau^2)z}\bigg)^{\frac{1}{2}+\alpha} d z
    \end{equation*}
    and
    \begin{equation*}
         I_3= \int_{(\frac{2}{x^2})\wedge1}^{1} e^{\frac{x^2z}{2}} z^{\alpha-1}  \bigg(\frac{1}{\tau^2 +(1-\tau^2)z}\bigg)^{\frac{1}{2}+\alpha} d z,
    \end{equation*}
    where $y_1\wedge y_2$ denotes minimum of $y_1$ and $y_2$. Next, making the substitution $z=u\tau^2$ in $I_1$,
    \begin{equation*} \label{eq:4}
        I_1= \tau^{-1} \int_{0}^{\frac{1}{\tau}} u^{\alpha-1} (1+u(1-\tau^2))^{-(\frac{1}{2}+\alpha)} e^{\frac{x^2 \tau^2 u}{2}} du \tag{51}.
    \end{equation*}
    Next, define 
    \begin{align*}
        I_1^{*} &= \tau^{-1} \int_{0}^{\frac{1}{\tau}} u^{\alpha-1} (1+u(1-\tau^2))^{-(\frac{1}{2}+\alpha)} du .
    \end{align*}
    Our target is to show that  $ \frac{I_1- I_1^{*}}{ I_1^{*}} \to 0$ as $\tau \to 0$.
    Now following the argument same as that used in Lemma C.9 of \cite{van2017adaptive}, for $|x| \leq B \kappa_{\tau}$, the exponent in the integral tends to $1$, uniformly in $u \leq \frac{1}{\tau}$. 
    Since, for any $y \geq 0$, $e^{y}-1 \leq y e^{y}$, replacing $e^{\frac{x^2 \tau^2 u}{2}}$ by 1 
    the difference between $I_1$ and $I_1^{*}$ is at most of the order of
    \begin{align*}
       &  \tau^{-1} \int_{0}^{\frac{1}{\tau}} u^{\alpha-1}  (1+u(1-\tau^2))^{-(\frac{1}{2}+\alpha)} [e^{\frac{x^2 \tau^2 u}{2}}-1]du \\
         & \lesssim x^2 \tau e^{\frac{x^2\tau}{2}} \int_{0}^{\frac{1}{\tau}} u^{-\frac{1}{2}} du  \lesssim {\tau}^{\frac{1}{2}} \log (\frac{1}{\tau}) (1+o(1)).
    \end{align*}
    
    Next, observe that, $\frac{1}{1+(1-\tau^2)u}=\frac{1}{(1+u)(1-\tau^2)}(1+O(\tau^2))$. Since,
    \begin{align*} \label{eq:5}
        \int_{0}^{\infty} u^{\alpha-1} (1+u)^{-(\alpha+\frac{1}{2})} du = \int_{0}^{\infty} u^{-\frac{3}{2}} (1+\frac{1}{u})^{-\frac{1}{2}-\alpha} du= K^{-1} \tag{52}
    \end{align*}
    and note that
   \begin{equation*} \label{eq:6}
       \int_{\frac{1}{\tau}}^{\infty} u^{\alpha-1} (1+u)^{-(\alpha+\frac{1}{2})} du \lesssim {\tau}^{\frac{1}{2}}\tag{53}.
   \end{equation*}
   As a result, we have as $\tau \to 0$,
   \begin{equation*}
        I_1^{*}= \frac{K^{-1}}{\tau}[1+O(\tau^{\frac{1}{2}})],
   \end{equation*}
   uniformly in $|x| \leq B \kappa_{\tau}$. This implies
   \begin{align*}
      0< \frac{I_1- I_1^{*}}{ I_1^{*}} & \lesssim {\tau}^{\frac{3}{2}} \log (\frac{1}{\tau}) (1+o(1)),
   \end{align*}
   and hence $ \frac{I_1- I_1^{*}}{ I_1^{*}} \to 0$ as $\tau \to 0$.
    Combining all these arguments along with \eqref{eq:5} and \eqref{eq:6}, we obtain
    \begin{equation*}  \label{eq:7}
        I_1= \frac{K^{-1}}{\tau}[1+O(\tau^{\frac{1}{2}})],\tag{54}
    \end{equation*}
    uniformly in $|x| \leq B \kappa_{\tau}$. Moving towards the second integral, first, we make a transformation $\frac{x^2z}{2}=v$ and hence
    \begin{equation*}
        I_2= \bigg(\frac{x^2}{2}\bigg)^{\frac{1}{2}} \int_{\frac{x^2\tau}{2}}^{\frac{x^2}{2} \wedge1} e^{v} v^{\alpha-1}  \bigg(\frac{\tau^2 x^2}{2}+(1-\tau^2)v \bigg)^{-(\frac{1}{2}+\alpha)} dv.
    \end{equation*}
    Now, we bound $\frac{\tau^2 x^2}{2}+(1-\tau^2)v$ below by $(1-\tau^2)v$ and observe that, the upper limit of the range of integration can be bounded by $1$ irrespective of whether $\frac{x^2}{2} \leq 1$ or not. Hence, we can show that
    \begin{equation*}
        I_2 \lesssim \tau^{-\frac{1}{2}}
    \end{equation*}
    and this contributes negligible compared to $I_1$. Finally, for $I_3$, again after the same transformation
    \begin{equation*}
        I_3= \bigg(\frac{x^2}{2}\bigg)^{\frac{1}{2}} \int_{1}^{\frac{x^2}{2} } e^{v} v^{\alpha-1} \bigg(\frac{\tau^2 x^2}{2}+(1-\tau^2)v\bigg)^{-(\frac{1}{2}+\alpha)} dv.
    \end{equation*}
    Here observe that, the integral contributes nothing when $\frac{x^2}{2} \leq 1$, hence we are only interested when $\frac{x^2}{2} >1$. Next, we define 
    \begin{equation*}
        I_3^{*}= \bigg(\frac{x^2}{2}\bigg)^{\frac{1}{2}} \int_{1}^{\frac{x^2}{2}}e^{v} v^{-\frac{3}{2}} dv.
    \end{equation*}
   Now our target is to show that the difference between $I_3$ and $I_3^{*}$ is negligible compared to $I_3^{*}$ as $\tau \to 0$. In order to prove that, first note, 
    \begin{equation*} \label{eq:7a}
      \bigg ( \frac{1}{\frac{\tau^2 x^2}{2}+(1-\tau^2)v} \bigg)^{\frac{1}{2}+\alpha} \leq \bigg(\frac{1}{v} \bigg)^{\frac{1}{2}+\alpha} \bigg[1+\frac{\tau^2(v+{x^2})}{v(1-\tau^2)}\bigg]^{\frac{1}{2}+\alpha}. \tag{55}
    \end{equation*}
    Now, first, consider the case when $\alpha+\frac{1}{2}$ is a positive integer. 
     Hence using the Binomial theorem, we have
     \begin{align*}
          \bigg ( \frac{1}{\frac{\tau^2 x^2}{2}+(1-\tau^2)v} \bigg)^{\frac{1}{2}+\alpha} -v^{-(\frac{1}{2}+\alpha)} & \leq v^{-(\frac{1}{2}+\alpha)} \sum_{j=1}^{\frac{1}{2}+\alpha} {\frac{1}{2}+\alpha \choose j} \bigg[\frac{\tau^2}{(1-\tau^2)} \bigg(1+\frac{x^2}{v}\bigg)\bigg]^{j} 
     \end{align*}
    Next, observing that for $1 \leq v \leq \frac{x^2}{2}, 2\leq \frac{x^2}{v} \leq x^2$, for $|x| \leq B \kappa_{\tau}$, the difference between $I_3$ and $I_3^{*}$ can be bounded as
    \begin{align*}
         I_3-I_3^{*} 
         & \lesssim \bigg(\frac{\tau^2}{1-\tau^2}\bigg) \log \bigg(\frac{1}{\tau}\bigg) (1+o(1)) \bigg(\frac{x^2}{2}\bigg)^{\frac{1}{2}}\int_{1}^{\frac{x^2}{2}} e^v v^{-\frac{3}{2}} dv,
    \end{align*}
    which implies
    \begin{equation*}
        I_3-I_3^{*}  \lesssim \bigg(\frac{\tau^2}{1-\tau^2}\bigg) \log \bigg(\frac{1}{\tau}\bigg) (1+o(1)) I_3^{*} .
    \end{equation*}
     On the other hand, observe that,
    \begin{align*}
        I_3-I_3^{*} 
         & = \bigg(\frac{x^2}{2}\bigg)^{\frac{1}{2}} \int_{1}^{\frac{x^2}{2}}e^{v} v^{-\frac{3}{2}} \bigg[(\frac{\tau^2x^2}{2}+1-\tau^2)^{-(\alpha+\frac{1}{2})} -1 \bigg] dv \\
         & \gtrsim \bigg[(\frac{\tau^2\zeta^2_{\tau}}{2}(1+o(1))+1-\tau^2)^{-(\alpha+\frac{1}{2})} -1 \bigg] I_3^{*}.
    \end{align*}
    These two bounds ensure $\frac{I_3-I_3^{*}}{I_3^{*}} \to 0$ as $\tau \to 0$. Next, consider the case when $\alpha+\frac{1}{2}$ is a fraction. When $\alpha+\frac{1}{2}$ is a fraction, then there exists another fraction $b>0$ such that $\alpha+\frac{1}{2}+b$ is a positive integer. Hence, in this case, $\bigg[1+\frac{\tau^2(v+{x^2})}{v(1-\tau^2)}\bigg]^{\frac{1}{2}+\alpha} \leq \bigg[1+\frac{\tau^2(v+{x^2})}{v(1-\tau^2)}\bigg]^{\alpha+\frac{1}{2}+b}$. Now, applying exactly the same set of arguments on $\alpha+\frac{1}{2}+b$ in place of $\alpha+\frac{1}{2}$, we again can show that, $\frac{I_3-I_3^{*}}{I_3^{*}} \to 0$ as $\tau \to 0$.
    This completes the proof for $|x| \leq B \kappa_{\tau}$.\\
    \textbf{Case 2}:- When $|x| > B \kappa_{\tau}$, choose any $A \in (0,1)$. In this case, 
    the range of the integration in $I_{\alpha,\frac{1}{2}}(x)$ is divided into two parts, namely 
    \begin{equation*}
        I_4= \int_{0}^{A} e^{\frac{x^2z}{2}} z^{\alpha-1}  \bigg(\frac{1}{\tau^2 +(1-\tau^2)z}\bigg)^{\frac{1}{2}+\alpha} d z
    \end{equation*}
    and
    \begin{equation*}
        I_5= \int_{A}^{1} e^{\frac{x^2z}{2}} z^{\alpha-1}  \bigg(\frac{1}{\tau^2 +(1-\tau^2)z}\bigg)^{\frac{1}{2}+\alpha} d z.
    \end{equation*}
    Note that
    \begin{align*}
       I_4  & \leq e^{\frac{x^2A}{2}} \bigg(\frac{1}{\tau^2}\bigg)^{\frac{1}{2}+\alpha}  \int_{0}^{A}  z^{\alpha-1}  dz \\
        & \lesssim \tau^{-1-2\alpha} e^{\frac{x^2A}{2}}.
    \end{align*}
    One can choose $B \geq 1$ and $A \in (0,1)$ such that $B^2(1-A)>2\alpha+\frac{3}{2}$, which implies $I_4 \ll \tau^{\frac{1}{2}} \frac{e^{x^2/2}}{x^2/2}$, and hence the contribution is negligible compared to the second term in the expression of $I_{\alpha,\frac{1}{2}}(x)$ as given in \eqref{eq:L-2.1}. Finally, for $I_5$, we use that, for $z \geq A$, $\frac{1}{\tau^2+(1-\tau^2)z}=\frac{1}{z}[1+O(\tau^2)]$.
    This implies $I_5$ is of the form
    \begin{align*} \label{eq:8}
        I_5 & =  \int_{A}^{1} z^{-\frac{3}{2}} e^{\frac{x^2z}{2}} dz [1+O(\tau^2)]^{\frac{1}{2}+\alpha} \tag{56}.
    \end{align*}
    Next, note that, with the transformation $\frac{x^2z}{2}=v$,
    \begin{align*} \label{eq:9}
        \int_{A}^{1} z^{-\frac{3}{2}} e^{\frac{x^2z}{2}} dz &= \bigg(\frac{x^2}{2}\bigg)^{\frac{1}{2}} \int_{\frac{x^2A}{2}}^{\frac{x^2}{2}} e^v v^{-\frac{3}{2}} dv \\
        &= \bigg(\frac{x^2}{2}\bigg)^{\frac{1}{2}} \bigg[\int_{1}^{\frac{x^2}{2}}-\int_{1}^{\frac{x^2A}{2}} \bigg] e^v v^{-\frac{3}{2}} dv  \tag{57}.
    \end{align*}
    Now, using the first assertion of Lemma C.8 of \cite{van2017adaptive}, the second integral is bounded above by a multiple of $(x^2/2)^{-1} e^{x^2A/2}$, which is negligible compared to the first (this is of the order of $(x^2/2)^{-1} e^{x^2/2}$). Hence, combining \eqref{eq:8} and \eqref{eq:9}, we immediately have
    \begin{equation*} \label{eq:10}
        I_5 = \bigg(\frac{x^2}{2}\bigg)^{\frac{1}{2}} \int_{1}^{\frac{x^2}{2}}e^{v} v^{-\frac{3}{2}} dv \big[1+O(\tau^2)\big] \tag{58}.
    \end{equation*}
    Combining \eqref{eq:7} and \eqref{eq:10}, the r.h.s. of \eqref{eq:L-2.1} is established. \\
   On the other hand, expanding the integral given in \eqref{eq:L-2.1} with the help of Lemma C.8 of \cite{van2017adaptive} provides \eqref{eq:L-2.2}.
\end{proof}

\begin{lemma}
\label{lem6}
    There exist functions $R_{\tau,1}$ with $\sup_{x}|R_{\tau,1}(x)|={O}(\tau^{\frac{1}{2}})$ as $\tau \to 0$, such that for $\alpha \geq \frac{1}{2}$,
    \begin{equation*} \label{eq:L-3.1}
        J_{\alpha+1,\frac{1}{2}}(x) = \bigg(\frac{x^2}{2}\bigg)^{-\frac{1}{2}} \int_{0}^{\frac{x^2}{2}} e^v v^{-\frac{1}{2}} dv (1+R_{\tau,1}(x)) \lesssim (1 	\wedge x^{-2}) e^{\frac{x^2}{2}} \tag{59}
    \end{equation*}
    and
    \begin{equation*} \label{eq:L-3.2}
        J_{\alpha+1,\frac{1}{2}}(x)- J_{\alpha+2,\frac{1}{2}}(x) = \bigg(\frac{x^2}{2}\bigg)^{-\frac{1}{2}} \int_{0}^{\frac{x^2}{2}} e^v v^{-\frac{1}{2}} \bigg(1-\frac{2v}{x^2}\bigg) dv (1+R_{\tau,1}(x)) \lesssim (1 	\wedge x^{-4}) e^{\frac{x^2}{2}}. \tag{60}
    \end{equation*}
\end{lemma}
\begin{proof}
    Recall that $J_{\alpha+1,\frac{1}{2}}(x)$ is defined as
    \begin{equation*} \label{eq:11}
        J_{\alpha+1,\frac{1}{2}}(x) = \int_{0}^{1} e^{\frac{x^2z}{2}} z^{\alpha}  \bigg(\frac{1}{\tau^2 +(1-\tau^2)z}\bigg)^{\frac{1}{2}+\alpha} d z \tag{61}.
    \end{equation*}
    Next, we split the range of the integration as $[0, \tau]$ and $[\tau,1]$. Note that, the contribution of the first integral obtained from \eqref{eq:11} is bounded by
    \begin{align*} \label{eq:12}
         \int_{0}^{ \tau} e^{\frac{x^2z}{2}} z^{\alpha}  \bigg(\frac{1}{\tau^2 +(1-\tau^2)z}\bigg)^{\frac{1}{2}+\alpha} d z  & \lesssim e^{\frac{x^2 \tau}{2}}  \tau^{\frac{1}{2}} \tag{62}.
    \end{align*}
    On the other hand, note that, when $z \geq \tau$, $\frac{1}{\tau^2+(1-\tau^2)z}=\frac{1}{z}[1+O(\tau)]$, hence, we have
    \begin{align*} \label{eq:13}
         \int_{ \tau}^{1} e^{\frac{x^2z}{2}} z^{\alpha}  \bigg(\frac{1}{\tau^2 +(1-\tau^2)z}\bigg)^{\frac{1}{2}+\alpha} d z &=  \int_{ \tau}^{1} e^{\frac{x^2z}{2}} z^{-\frac{1}{2}} dz [1+O(\tau)]^{\frac{1}{2}+\alpha} \\
         & \gtrsim \exp(\frac{x^2}{2}\tau) {(1-\tau^{\frac{1}{2}})} [1+O(\tau)] \tag{63}.
    \end{align*}
    Combining \eqref{eq:12} and \eqref{eq:13} states that the contribution of the first integral is negligible as compared to the second one. Hence, one has
    \begin{align*}
        J_{\alpha+1,\frac{1}{2}}(x) = \int_{ \tau}^{1} e^{\frac{x^2z}{2}}  z^{-\frac{1}{2}} dz (1+O(\tau)+O(\tau^{\frac{1}{2}})).
    \end{align*}
      Also note that, $\tau=O(\tau^{\frac{1}{2}})$.  Observe that, \eqref{eq:12} and \eqref{eq:13} also imply that, as $\tau \to 0$,
      \begin{align*}
          \int_{\tau}^{1}  e^{\frac{x^2z}{2}} z^{-\frac{1}{2}} dz &= \int_{0}^{1} e^{\frac{x^2z}{2}}  z^{-\frac{1}{2}} dz[1+O(\tau^{\frac{1}{2}})].
      \end{align*}
      As a result,
      \begin{align*}
          J_{\alpha+1,\frac{1}{2}}(x) = \int_{0}^{1} e^{\frac{x^2z}{2}}  z^{-\frac{1}{2}} dz (1+O(\tau^{\frac{1}{2}})).
      \end{align*}
      Then the equality follows due to the change of variable $\frac{x^2z}{2}=v$. The second one is due to the exactly same set of arguments used in Lemma C.10 of \cite{van2017adaptive}. \\
      Proof of the second assertion follows using similar set of arguments.
\end{proof}

\begin{lemma}
    \label{lem7}
    The function $x \mapsto m_{\tau}(x)$ is symmetric about $0$ and non-decreasing in $[0,\infty)$ with \\
    $(i) -2\alpha \leq m_{\tau}(x) \leq 2a$, for all $x \in \mathbb{R}$ and $\tau \in (0,1)$.\\
    $(ii) |m_{\tau}(x)| \lesssim \tau e^{x^2/2}(x^{-2} \wedge 1)$, as $\tau \to 0$, for every $x$.   
\end{lemma}
\begin{proof}
   The symmetric behavior follows from the definition of $m_{\tau}(x)$ as given in \eqref{eq:L-1.1}.\\
    For monotonicity, using \eqref{eq:2} of Lemma \ref{lem4}, it readily follows that
    \begin{align*} \label{eq:L-4.1}
        m_{\tau}(x) &= 2a+ \tau \frac{\dot{I}_{{\alpha,a}(x)}}{I_{\alpha,a}(x)} \\
        &= 2a -2 \tau^2(a+\alpha) \frac{\int_{0}^{1} e^{\frac{x^2z}{2}} z^{\alpha-1} (1-z)^{a+\frac{1}{2}}  (\frac{1}{N(z})^{a+\alpha+1} dz  }{\int_{0}^{1} e^{\frac{x^2z}{2}} z^{\alpha-1} (1-z)^{a-\frac{1}{2}}  (\frac{1}{N(z})^{a+\alpha} dz  } \tag{64} \\
        &= 2a +2 \tau^2 (a+\alpha) \frac{\int_{0}^{1} e^{\frac{x^2z}{2}} (\frac{z-1}{\tau^2+(1-\tau^2)z})  z^{\alpha-1} (1-z)^{a-\frac{1}{2}}  (\frac{1}{N(z})^{a+\alpha} dz}{\int_{0}^{1} z^{\alpha-1} e^{\frac{x^2z}{2}} (1-z)^{a-\frac{1}{2}}  (\frac{1}{N(z})^{a+\alpha} dz}\\
        &= 2a+2 \tau^2(a+\alpha) \int_{0}^{1} \bigg(\frac{z-1}{\tau^2+(1-\tau^2)z}\bigg) g_x(z) dz,
    \end{align*}
    where $z \mapsto g_x(z)$ is the probability density function on $[0,1]$ with $g_x(z) \propto z^{\alpha-1} e^{\frac{x^2z}{2}} (1-z)^{a-\frac{1}{2}}  \big(\frac{1}{N(z)}\big)^{a+\alpha}$. Next, following the same set of arguments as used in Lemma C.7  of  \cite{van2017adaptive}. \\
    (i) The upper bound is obvious by using \eqref{eq:L-4.1}. For the lower bound, note that
\begin{equation*}
    \tau \frac{\dot{I}_{{\alpha,a}(x)}}{I_{\alpha,a}(x)}  = -2 (\alpha+a)  \int_{0}^{1} \bigg(\frac{(1-z)\tau^2}{\tau^2+(1-\tau^2)z}\bigg) g_x(z) dz.
\end{equation*}
Since, for any $0 \leq z \leq 1$, $\tau^2+(1-\tau^2)z \geq \tau^2$ implies $\frac{(1-z)\tau^2}{\tau^2+(1-\tau^2)z} \leq 1$, the lower bound follows from it. \\
(ii) First, using the definition of $m_{\tau}(x)$ as given in \eqref{eq:L-1.1} and then followed by the triangle inequality, an upper bound on $|m_{\tau}(x)|$ is obtained as
\begin{align*}
    |m_{\tau}(x)| & \leq x^2 \frac{|J_{\alpha+1,a}(x)-J_{\alpha+2,a}(x)|}{I_{\alpha,a}(x)} +\frac{|J_{\alpha+1,a}(x)|}{I_{\alpha,a}(x)}.
\end{align*}
The assertion is proved by using \eqref{eq:L-3.1} and \eqref{eq:L-3.2} of Lemma \ref{lem6}, \eqref{eq:L-2.1} of Lemma \ref{lem5} and finally noting that $\frac{e^{x^2/2}}{x^2/2} \geq \frac{1}{\tau}$ for $|x| \geq B\kappa_{\tau}$ for $B \geq 1$.
\end{proof}
As an immediate consequence of Lemma \ref{lem7}, we have the following corollary.
\begin{corollary}
\label{cor1}
    Let $X \sim \mathcal{N}(\theta,1)$. Then, as $\tau \to 0$,
    \begin{align*}
        E_{\theta}m^2_{\tau}(X) =\begin{cases}
            O(\tau \zeta^{-2}_{\tau}), |\theta| \lesssim \zeta^{-1}_{\tau}, \\
            O(\tau^{\frac{1}{16}} \zeta^{-2}_{\tau}), |\theta| \leq \frac{\zeta_{\tau}}{4}.
        \end{cases}
    \end{align*}
\end{corollary}
\begin{proof}
    Noting that the upper bound of the absolute value of $m_{\tau}(x)$ as obtained in (ii) of Lemma \ref{lem7} matches with that of (vii) of Lemma C.7 of \cite{van2017adaptive}, the proof is immediate using the same set of arguments used in Lemma C.5 of \cite{van2017adaptive}.
\end{proof}

\begin{lemma}
    \label{lem8}
    Let $X \sim \mathcal{N}(\theta,1)$. For $|\theta| \lesssim \zeta^{-1}_{\tau}$, and $\tau \leq \tau_1 <\tau_2 $ and $\tau_2 \to 0$,
    \begin{align*} \label{eq:L-5.1}
        E_{\theta} \bigg( \frac{\zeta_{\tau_1}}{\tau_1} m_{\tau_1}(X) -\frac{\zeta_{\tau_2}}{\tau_2} m_{\tau_2}(X) \bigg)^2  \lesssim (\tau_2-\tau_1)^2 {\tau}^{-3}_1 \tag{54}.
    \end{align*}
    Further, for $|\theta| \leq \frac{\zeta_{\tau}}{4}$, and $\epsilon=\frac{1}{16}$, and $\tau \leq \tau_1 <\tau_2 $ and $\tau_2 \to 0$,
     \begin{align*}\label{eq:L-5.2}
        E_{\theta} \bigg( \frac{\zeta_{\tau_1}}{\tau^{\epsilon}_1} m_{\tau_1}(X) -\frac{\zeta_{\tau_2}}{\tau^{\epsilon}_2} m_{\tau_2}(X) \bigg)^2  \lesssim (\tau_2-\tau_1)^2 {\tau}^{-2-\epsilon}_1 \tag{66}.
    \end{align*}
\end{lemma}
\begin{proof}
    Using Lemma C.11 of \cite{van2017adaptive} with $V_{\tau}=\frac{\zeta_{\tau}}{\tau} m_{\tau}(X)$, the l.h.s. of \eqref{eq:L-5.1} can be upper bounded as
    \begin{align*}\label{eq:L-5.3}
      &  E_{\theta} \bigg( \frac{\zeta_{\tau_1}}{\tau_1} m_{\tau_1}(X) -\frac{\zeta_{\tau_2}}{\tau_2} m_{\tau_2}(X) \bigg)^2  \leq (\tau_2-\tau_1)^2 \sup_{\tau \in [\tau_1,\tau_2]} E_{\theta} \bigg( \frac{\zeta_{\tau}}{\tau} \dot{m}_{\tau}(X)- \frac{\zeta_{\tau}+\zeta^{-1}_{\tau}}{\tau^2} m_{\tau}(X) \bigg)^2 \\
        & \leq 2 (\tau_2-\tau_1)^2 \bigg[\sup_{\tau \in [\tau_1,\tau_2]}  E_{\theta}\bigg(\frac{\zeta_{\tau}}{\tau} \dot{m}_{\tau}(X) \bigg)^2  + \sup_{\tau \in [\tau_1,\tau_2]}  E_{\theta} \bigg(\frac{\zeta_{\tau}+\zeta^{-1}_{\tau}}{\tau^2} m_{\tau}(X) \bigg)^2 \bigg] \tag{67}.
    \end{align*}
    Note that, with the help of the first part of  Corollary \ref{cor1}, the second term in the r.h.s. of \eqref{eq:L-5.3} is bounded above by a constant times $ \sup_{\tau \in [\tau_1,\tau_2]}  \tau^{-3} \lesssim \tau^{-3}_1$. Hence, in order to show that \eqref{eq:L-5.1} holds, it is enough to show that the first term in the r.h.s. of \eqref{eq:L-5.1} is bounded above by a constant times $\tau^{-3}_1$. \\
    For the first term, observe that using \eqref{eq:L-1.2} with $a= \frac{1}{2}$,
    \begin{align*} \label{eq:L-5.4}
        \dot{m}_{\tau}(x) &= \frac{I_{\alpha,\frac{1}{2}}(x)[(x^2-1) \dot{J}_{\alpha+1,\frac{1}{2}}(x)-x^2 \dot{J}_{\alpha+2,\frac{1}{2}}(x)]-\dot{I}_{\alpha,\frac{1}{2}}(x)m_{\tau}(x)I_{\alpha,\frac{1}{2}}(x)}{(I_{\alpha,\frac{1}{2}}(x))^2} \\
        &= (x^2-1) \frac{\dot{J}_{\alpha+1,\frac{1}{2}}(x)}{I_{\alpha,\frac{1}{2}}(x)} -x^2 \frac{\dot{J}_{\alpha+2,\frac{1}{2}}(x)}{I_{\alpha,\frac{1}{2}}(x)}-\frac{\dot{I}_{\alpha,\frac{1}{2}}(x)}{I_{\alpha,\frac{1}{2}}(x)} m_{\tau}(x) \tag{68}.
    \end{align*}
    Now, note that, using the definition of $J_{\alpha+1,\frac{1}{2}}(x)$, $\dot{J}_{\alpha+1,\frac{1}{2}}(x)=2\tau(\alpha+\frac{1}{2}) (H_{\alpha+2,\frac{1}{2}}(x)-H_{\alpha+1,\frac{1}{2}} (x))$ and $\dot{J}_{\alpha+2,\frac{1}{2}}(x)=2\tau(\alpha+\frac{1}{2}) (H_{\alpha+3,\frac{1}{2}}(x)-H_{\alpha+2,\frac{1}{2}} (x))$ where $H_{\alpha+k,\frac{1}{2}}(x)=\int_{0}^{1} e^{\frac{x^2z}{2}} z^{\alpha+k-1} (1-z)^{a-\frac{1}{2}} \bigg(\frac{1}{\tau^2 +(1-\tau^2)z}\bigg)^{\alpha+\frac{3}{2}} d z$. Next, note that, $H_{\alpha+k,\frac{1}{2}}(x)$ is a decreasing function of $k$. Also, observe that, $H_{\alpha+1,\frac{1}{2}} (x) \leq I_{\alpha,\frac{1}{2}}(x)$. Finally, for the third term in the r.h.s. of \eqref{eq:L-5.4}, the definition of $ I_{\alpha,\frac{1}{2}}(x)$ implies
    \begin{equation*}
        -\frac{\dot{I}_{\alpha,\frac{1}{2}}(x)}{I_{\alpha,\frac{1}{2}}(x)} \leq 2\tau (\alpha+\frac{1}{2}). \frac{1}{\tau^2}.
    \end{equation*}
    Combining all these arguments provides an upper bound for the r.h.s. of \eqref{eq:L-5.4} as
    \begin{equation*}
        \dot{m}_{\tau}(x) \leq 2\tau (\alpha+\frac{1}{2}) [1+x^2+ \frac{1}{\tau^2} m_{\tau}(x)] .
    \end{equation*}
    As a consequence of this, 
    \begin{align*} \label{eq:L-5.5}
        E_{\theta}\dot{m}^2_{\tau}(X) & \lesssim \tau^2[1+E_{\theta}X^4+\frac{1}{\tau^4}E_{\theta}m^{2}_{\tau}(X)] \tag{69}.
    \end{align*}
    Note that, in this case, $E_{\theta}X^4$ is bounded by a constant and from the first part of Corollary \ref{cor1}, $E_{\theta}m^{2}_{\tau}(X)$ is bounded by $\tau \zeta^{-2}_{\tau}$. This results in $(\frac{\zeta_{\tau}}{\tau})^2 E_{\theta}\dot{m}^2_{\tau}(X)$ being bounded above by a multiple of $\tau^{-3}$ and \eqref{eq:L-5.1} is established. \\
    For proving \eqref{eq:L-5.2}, we again use Lemma C.11 of \cite{van2017adaptive}, but with $V_{\tau}=\frac{\zeta_{\tau}}{\tau^{\epsilon}} m_{\tau}(X)$. As a result, we have
     \begin{align*}
      & E_{\theta} \bigg( \frac{\zeta_{\tau_1}}{\tau^{\epsilon}_1} m_{\tau_1}(X) -\frac{\zeta_{\tau_2}}{\tau^{\epsilon}_2} m_{\tau_2}(X) \bigg)^2   \\
        & \leq 2 (\tau_2-\tau_1)^2 \bigg[\sup_{\tau \in [\tau_1,\tau_2]}  E_{\theta}\bigg(\frac{\zeta_{\tau}}{\tau^\epsilon} \dot{m}_{\tau}(X) \bigg)^2  + \sup_{\tau \in [\tau_1,\tau_2]}  E_{\theta} \bigg(\frac{\epsilon \zeta_{\tau}+\zeta^{-1}_{\tau}}{\tau^{1+\epsilon}} m_{\tau}(X) \bigg)^2 \bigg] .
    \end{align*}
    Using the second part of  Corollary \ref{cor1}, the second term in the r.h.s. of the above inequality is bounded above by a constant times $\tau^{-2-\epsilon}$.
     Next, we follow the same steps as before and obtain an expression same as \eqref{eq:L-5.5}. In this case, $E_{\theta}X^4$ is bounded by $\zeta^{4}_{\tau}$ and from second part of Corollary \ref{cor1}, $E_{\theta}m^{2}_{\tau}(X)$ is bounded by $\tau^{\epsilon} \zeta^{-2}_{\tau}$. This implies $(\frac{\zeta_{\tau}}{\tau^{1+\epsilon}})^2 E_{\theta}\dot{m}^2_{\tau}(X)$ being bounded above by a multiple of $\tau^{-2-\epsilon}$ and completes the proof of \eqref{eq:L-5.2}.
\end{proof}
The next corollary follows immediately as a byproduct of Lemmas mentioned previously.
%\begin{corollary}
 %   \label{cor2}
 %   For any $\epsilon_{\tau} \to 0$ and uniformly in $I_0 \subseteq \{i: |\theta_{0,i}| \leq \zeta^{-1}_{\tau} \}$ with $|I_0| \gtrsim n$, 
 %   \begin{align*}
%        \sup_{\frac{1}{n} \leq \tau \leq \epsilon_{\tau}} \frac{1}{|I_0|} |\sum_{i \in I_0} m_{\tau}(X_i)\frac{\zeta_{\tau}}{\tau}-\sum_{i \in I_0} E_{\theta_0} m_{\tau}(X_i)\frac{\zeta_{\tau}}{\tau} | \xrightarrow{P_{\theta_0}} 0.
 %   \end{align*}
%    Similarly, uniformly in $I_1 \subseteq \{i: |\theta_{0,i}| \leq \frac{\zeta_{\tau}}{4} \}$,
  %   \begin{align*}
 %       \sup_{\frac{1}{n} \leq \tau \leq \epsilon_{\tau}} \frac{1}{|I_1|} |\sum_{i \in I_1} m_{\tau}(X_i)\frac{\zeta_{\tau}}{\tau^{\frac{1}{32}}}-\sum_{i \in I_1} E_{\theta_0} m_{\tau}(X_i)\frac{\zeta_{\tau}}{\tau^{\frac{1}{32}}} | \xrightarrow{P_{\theta_0}} 0.
 %   \end{align*}
%\end{corollary}
%\begin{proof}
%    This proof follows from using the same set of arguments used in Lemma C.3 of \cite{van2017adaptive} along with Corollary \ref{cor1} and Lemma \ref{lem8} of this work.
%\end{proof}

\begin{corollary}
    \label{cor3}
    If the cardinality of $I_0 := \{i: \theta_{0,i}=0 \}$ tends to infinity, then
    \begin{align*}
        \sup_{\frac{1}{n} \leq \tau \leq \beta_n} \frac{1}{|I_0|} |\sum_{i \in I_0} m_{\tau}(X_i)-\sum_{i \in I_0} E_{\theta_0} m_{\tau}(X_i) | \xrightarrow{P_{\theta_0}} 0.
    \end{align*}
\end{corollary}
\begin{proof}
    The proof of this result follows using a similar set of arguments to those of Lemma C.6 of \cite{van2017adaptive} along with Lemma \ref{lem8}.
\end{proof}

\begin{lemma}
    \label{lem9}
    Let $X \sim \mathcal{N}(\theta,1)$. Then as $ \tau \to 0$
    \begin{align*}
        E_{\theta} m_{\tau}(X) =\begin{cases}
            -\frac{2}{K^{-1}\sqrt{2\pi}} \frac{\tau}{\zeta_{\tau}} (1+o(1)), & |\theta| =o(\zeta^{-2}_{\tau}), \\
             o(\tau^{\frac{1}{16}} \zeta^{-1}_{\tau}), & |\theta| \leq \frac{\zeta_{\tau}}{4}.
        \end{cases}
    \end{align*}
\end{lemma}
\begin{proof}
For proving the first assertion, following the steps of Proposition C.2 of \cite{van2017adaptive}, we can show that both
$\int_{|x| \geq \kappa_{\tau}} m_{\tau}(x) \phi(x-\theta) dx=o(\frac{\tau}{\zeta_{\tau}})$ and \\ $\int_{ \zeta_{\tau} \leq |x| \leq \kappa_{\tau}} m_{\tau}(x) \phi(x-\theta) dx=o(\frac{\tau}{\zeta_{\tau}})$,where $\kappa_{\tau} \sim \zeta_{\tau} + \frac{2 \log \zeta_{\tau}}{\zeta_{\tau}}$. The remaining argument also follows using similar sets of arguments used in Proposition C.2 of \cite{van2017adaptive} along with some algebraic manipulations. However, for the sake of completeness, we present all the steps.
\\
Now, note that, from \eqref{eq:L-2.1} of Lemma \ref{lem5}, when $\frac{x^2}{2} \leq 1$, $I_{\alpha,\frac{1}{2}}(x)=\frac{K^{-1}}{\tau}[1+O(\sqrt{\tau})]$. On the other hand, since, $\frac{e^{\frac{x^2}{2}}}{x^2}$ is increasing for large values of $x$ and attains the value $\frac{\tau}{\zeta^2_{\tau}}$ at $x=\zeta_{\tau}$, 
by \eqref{eq:L-2.1} of Lemma \ref{lem5}, $I_{\alpha,\frac{1}{2}}(x)=\frac{K^{-1}}{\tau}[1+O(\frac{1}{\zeta^2_{\tau}})]$, when $1 \leq \frac{x^2}{2} \leq \log (\frac{1}{\tau})$. Combining these two facts, we get that, $I_{\alpha,\frac{1}{2}}(x)=\frac{K^{-1}}{\tau}[1+O(\frac{1}{\zeta^2_{\tau}})]$,
uniformly in $x \in (0,\zeta_{\tau})$. Hence,
\begin{align*} \label{eq:L-6.1}
    \int_{|x| \leq \zeta_{\tau}} m_{\tau}(x) \phi(x-\theta) dx &= \int_{0}^{\zeta_{\tau}}  \frac{x^2(J_{\alpha+1,\frac{1}{2}}(x)-J_{\alpha+2,\frac{1}{2}}(x))-J_{\alpha+1,\frac{1}{2}}(x)}{\frac{K^{-1}}{\tau}} \phi(x) dx +R_{\tau}  \tag{70},
\end{align*}
where the absolute value of $R_{\tau}$ is bounded by $\int_{0}^{\zeta_{\tau}}|x^2(J_{\alpha+1,\frac{1}{2}}(x)-J_{\alpha+2,\frac{1}{2}}(x))-J_{\alpha+1,\frac{1}{2}}(x)|\phi(x) dx$
times $\sup_{|x| \leq \zeta_{\tau}} |\frac{\phi(x-\theta)}{I_{\alpha,\frac{1}{2}}(x) \phi(x)}- \frac{K}{\tau^{-1}}|$. By using Lemma \ref{lem6}, the integrand is bounded above by a constant for $x$ near $0$ and by a multiple of $x^{-2}$ otherwise, which makes the integral to be bounded. Next, observe that,
\begin{align*}
    \sup_{|x| \leq \zeta_{\tau}} |\frac{\phi(x-\theta)}{I_{\alpha,\frac{1}{2}}(x) \phi(x)}- \frac{K}{\tau^{-1}}| & = \sup_{|x| \leq \zeta_{\tau}} \frac{\tau K}{I_{\alpha,\frac{1}{2}}(x)} |\tau^{-1} K^{-1} (e^{x \theta-\frac{\theta^2}{2}})-I_{\alpha,\frac{1}{2}}(x)| \\
    &= \sup_{|x| \leq \zeta_{\tau}} \tau  K |\frac{e^{x \theta-\frac{\theta^2}{2}}}{1+O(\frac{1}{\zeta^2_{\tau}})}-1| 
     \lesssim \tau [\frac{1}{\zeta^2_{\tau}}+e^{\zeta_{\tau} \theta-\frac{\theta^2}{2}}-1].
\end{align*}
Observe that, for $|\theta|=o(\zeta^{-2}_{\tau}), \zeta_{\tau}|\theta|-\frac{\theta^2}{2}=o(\zeta^{-1}_{\tau})$ and using the fact $e^{y}-1 \sim y$ as $y \to 0$,  we have
\begin{align*}
    \sup_{|x| \leq \zeta_{\tau}} |\frac{\phi(x-\theta)}{I_{\alpha,\frac{1}{2}}(x) \phi(x)}- \frac{K}{\tau^{-1}}| & \lesssim \tau[\frac{1}{\zeta^2_{\tau}}+o(\zeta^{-1}_{\tau})] =o(\frac{\tau}{\zeta_{\tau}}).
\end{align*}
These two arguments show that $R_{\tau}$ is negligible compared to $\frac{\tau}{\zeta_{\tau}}$.\\
 Next, using Fubini's Theorem, the above integral in \eqref{eq:L-6.1} can be rewritten as
    \begin{align*} \label{eq:L-6.2}
       & \int_{0}^{\zeta_{\tau}}  \frac{x^2(J_{\alpha+1,\frac{1}{2}}(x)-J_{\alpha+2,\frac{1}{2}}(x))-J_{\alpha+1,\frac{1}{2}}(x)}{\frac{K^{-1}}{\tau}} \phi(x) dx \\
        &= K \tau \int_{0}^{1} \int_{0}^{\zeta_{\tau}} z^{\alpha} \bigg(\frac{1}{N(z)}\bigg)^{\frac{1}{2}+\alpha}  e^{\frac{x^2z}{2}} [x^2(1-z)-1] \frac{e^{-\frac{x^2}{2}}}{\sqrt{2 \pi}} dx \hspace*{0.1cm} dz  \tag{71}\\
        &= K \tau \int_{0}^{1} z^{\alpha} \bigg(\frac{1}{N(z)}\bigg)^{\frac{1}{2}+\alpha}  \int_{0}^{\zeta_{\tau}} [x^2(1-z)-1] \frac{e^{-\frac{x^2(1-z)}{2}}}{\sqrt{2 \pi}} dx \hspace*{0.1cm} dz .
    \end{align*}
    Note that the inner integral becomes zero if the range of integration is over $(0, \infty)$ instead of $(0, \zeta_{\tau})$. Hence, we have
    \begin{align*}
        & \int_{0}^{\zeta_{\tau}} \frac{x^2(J_{\alpha+1,\frac{1}{2}}(x)-J_{\alpha+2,\frac{1}{2}}(x))-J_{\alpha+1,\frac{1}{2}}(x)}{\frac{K^{-1}}{\tau}} \phi(x) dx \\
        &=- K \tau \int_{0}^{1} z^{\alpha} \bigg(\frac{1}{N(z)}\bigg)^{\frac{1}{2}+\alpha} \int_{\zeta_{\tau}}^{\infty}[x^2(1-z)-1] \frac{e^{-\frac{x^2(1-z)}{2}}}{\sqrt{2 \pi}} dx \hspace*{0.1cm} dz \\
        &= - K \tau \int_{0}^{1} z^{\alpha} \bigg(\frac{1}{N(z)}\bigg)^{\frac{1}{2}+\alpha}  \frac{\zeta_{\tau}}{\sqrt{2\pi}} e^{-\frac{\zeta^2_{\tau}(1-z)}{2}} dz   .
    \end{align*}
    Here, in the last line, we use the fact $\int_{y}^{\infty} [(vb)^2-1]\phi(vb) dv= y \phi(yb)$. Next, similar to Proposition C.2 of \cite{van2017adaptive}, we split the range of integration as $(0,\frac{1}{2}]$ and $(\frac{1}{2},1)$. When $0 \leq z\leq \frac{1}{2}$, the absolute value of the integral is bounded as
    \begin{align*}
        K \tau \int_{0}^{\frac{1}{2}} z^{\alpha} \bigg(\frac{1}{N(z)}\bigg)^{\frac{1}{2}+\alpha}  \frac{\zeta_{\tau}}{\sqrt{2\pi}} e^{-\frac{\zeta^2_{\tau}(1-z)}{2}} dz & \lesssim K \frac{e^{-\frac{\zeta^2_{\tau}}{4}}\tau \zeta_{\tau}}{\sqrt{2 \pi} (1-\tau^2)^{\frac{1}{2}+\alpha}} \int_{0}^{\frac{1}{2}} z^{-\frac{1}{2}} dz\\
        &= O(e^{-\frac{\zeta^2_{\tau}}{4}}\tau \zeta_{\tau}) = o(\frac{\tau}{\zeta_{\tau}}).
    \end{align*}
    On the other hand, when $\frac{1}{2} \leq z \leq 1$, we again use, $\frac{1}{\tau^2+(1-\tau^2)z}=\frac{1}{z}[1+O(\tau^2)]$. This implies
    \begin{align*} 
       &  - K \tau \int_{\frac{1}{2}}^{1} z^{\alpha} \bigg(\frac{1}{N(z)}\bigg)^{\frac{1}{2}+\alpha}  \frac{\zeta_{\tau}}{\sqrt{2\pi}} e^{-\frac{\zeta^2_{\tau}(1-z)}{2}} dz  \\
        &=  -\frac{K \tau \zeta_{\tau}}{\sqrt{2\pi}} \int_{\frac{1}{2}}^{1} z^{-\frac{1}{2}} e^{-\frac{\zeta^2_{\tau}(1-z)}{2}} dz [1+O(\tau^2)] \\
         &=  -\frac{K}{\sqrt{2\pi}} \frac{\tau}{\zeta_{\tau}} \int_{0}^{\frac{\zeta^2_{\tau}}{2}} e^{-\frac{u}{2}} \frac{1}{(1-\frac{u}{\zeta^2_{\tau}})^{\frac{1}{2}}} du[1+O(\tau^2)],
    \end{align*}
    where the equality is due to the substitution $\zeta^2_{\tau}(1-z)=u$. Finally, following the same argument as used in Proposition C.2 of \cite{van2017adaptive}, the integral tends to
    $\int_{0}^{\infty} e^{-\frac{u}{2}} du =2$, and completes the proof of the first assertion.\\
    For proving the second statement, we follow the steps mentioned in Proposition C.2 of \cite{van2017adaptive}.
\end{proof}

\subsection{Results related to the contraction rate of Hierarchical Bayes}

% \begin{proof}[Proof of Lemma 2]
 %    Here we use the same set of arguments used in Lemma 3.6 of \cite{van2017adaptive}, with $C_e=\frac{2^{\frac{3}{2}}K}{\sqrt{\pi}} \tau_n(q_n)$ and obtain that,  with $P_{\theta_0}$-probability tending to one, 
%     \begin{align*}
%       \Pi(\tau \geq 5t_n |\mathbf{X}) & \lesssim \frac{e^{-aq_n}}{\int_{\frac{t_n}{2}}^{t_n} \pi(\tau) d \tau},
%     \end{align*}
%     which tends to zero using (\hyperlink{C3}{C3}).     As a result, for $b_n=(\frac{q_n}{n})^{M_1}$ and $h(\mathbf{X})= \Pi(\tau \geq 5t_n |\mathbf{X})$,
 %    \begin{align*}
%         E_{\theta_0} h(\mathbf{X}) &= E_{\theta_0}[h(\mathbf{X}) 1_{\{h(\mathbf{X}) \lesssim b_n\} }] +E_{\theta_0}[h(\mathbf{X}) 1_{\{h(\mathbf{X}) \gtrsim b_n\} }] \\
 %        & \lesssim b_n +P_{\theta_0}[h(\mathbf{X}) \gtrsim b_n]
%     \end{align*}
%     The proof is completed since both of the terms go to zero.
 %\end{proof}

 \begin{proof}[Proof of Theorem \ref{Thm-3}]
     In order to prove Theorem 3 with the use of Lemma 3 of the main document,
     \begin{align*} \label{eq:T-2.1}
        & \mathbb{E}_{\boldsymbol{\theta}_0} \Pi \bigg(  \boldsymbol{\theta}:||\boldsymbol{\theta}-\boldsymbol{\theta}_0 ||^2>    M_n q_n \log n |\mathbf{X} \bigg) \\
        &= \mathbb{E}_{\boldsymbol{\theta}_0}\bigg[\int_{\frac{1}{n}}^{1} \Pi \big(  \boldsymbol{\theta}:||\boldsymbol{\theta}-\boldsymbol{\theta}_0 ||^2>    M_n q_n \log n |\mathbf{X},\tau \big) \pi(\tau|\mathbf{X}) d \tau
        \bigg] \\
        &= \mathbb{E}_{\boldsymbol{\theta}_0} \bigg[\big(\int_{\frac{1}{n}}^{5t_n} +\int_{5t_n}^{1} \big) \Pi \big(  \boldsymbol{\theta}:||\boldsymbol{\theta}-\boldsymbol{\theta}_0 ||^2>    M_n q_n \log n |\mathbf{X},\tau \big) \pi(\tau|\mathbf{X}) d \tau \bigg] \\
        & \leq  \mathbb{E}_{\boldsymbol{\theta}_0} \bigg[\int_{\frac{1}{n}}^{5t_n} \Pi \big(  \boldsymbol{\theta}:||\boldsymbol{\theta}-\boldsymbol{\theta}_0 ||^2>    M_n q_n \log n |\mathbf{X},\tau \big) \pi(\tau|\mathbf{X}) d \tau \bigg] +o(1) \\
        & \leq  \mathbb{E}_{\boldsymbol{\theta}_0}[\sup_{\frac{1}{n} \leq \tau \leq 5t_n} \Pi \big(  \boldsymbol{\theta}:||\boldsymbol{\theta}-\boldsymbol{\theta}_0 ||^2>    M_n q_n \log n |\mathbf{X},\tau \big) ]+o(1)
        \tag{72}.
     \end{align*}
     Next, applying Markov's inequality to the first term in the r.h.s. of \eqref{eq:T-2.1}, it is sufficient to show that
     \begin{align*} \label{eq:T-2.2}
\mathbb{E}_{\boldsymbol{\theta}_0}\bigg[\sup_{\frac{1}{n} \leq \tau \leq 5t_n} \sum_{i=1}^{n} Var(\theta_i|\mathbf{X},\tau)   \bigg] & \lesssim q_n \log n \tag{73}
     \end{align*}
     and
     \begin{align*} \label{eq:T-2.3}
\mathbb{E}_{\boldsymbol{\theta}_0}\bigg[\sup_{\frac{1}{n} \leq \tau \leq 5t_n} ||\boldsymbol{\theta}_0- T_{\tau}(\mathbf{X})||^2   \bigg] & \lesssim q_n \log n \tag{74}.
\end{align*}
Noting that the posterior distribution of $\theta_i$ given $(\boldsymbol{X},\tau)$ depends on $(X_i,\tau)$ only, we are left to prove that
 \begin{align*} \label{eq:T-2.4}
\mathbb{E}_{\boldsymbol{\theta}_0}\bigg[\sup_{\frac{1}{n} \leq \tau \leq 5t_n} \sum_{i=1}^{n} Var(\theta_i|X_i,\tau)   \bigg] & \lesssim q_n \log n \tag{75}.
 \end{align*}
In order to prove \eqref{eq:T-2.4}, we need to show that
\begin{align*} \label{eq:T-2.6}
    \sum_{i:\theta_{0i} \neq 0} \mathbb{E}_{{\theta}_{0i}} \bigg[\sup_{\frac{1}{n} \leq \tau \leq 5t_n} Var(\theta_i|X_i,\tau)   \bigg] & \lesssim q_n \log n \tag{76}
\end{align*}
and
\begin{align*} \label{eq:T-2.7}
    \sum_{i:\theta_{0i} = 0} \mathbb{E}_{{\theta}_{0i}} \bigg[\sup_{\frac{1}{n} \leq \tau \leq 5t_n} Var(\theta_i|X_i,\tau)   \bigg] & \lesssim q_n \log n \tag{77}
\end{align*}
\textbf{Step-1} Fix any $i$ such that $\theta_{0i} \neq 0$. Define $r_n=\sqrt{2 \rho^2 \log n}$ where $\rho$ is defined in the proof of Theorem 1. Since, for any fixed $x_i \in \mathbb{R}$ and $\tau>0$, $Var(\theta_i|x_i,\tau) \leq 1+x^2_i$,
\begin{align*}  \label{eq:T-2.8}
 \mathbb{E}_{{\theta}_{0i}}\bigg( \sup_{\frac{1}{n} \leq \tau \leq 5t_n} Var(\theta_i|X_i,\tau) 1_{\{|X_i|\leq \sqrt{2 \rho^2 \log n}\}}   \bigg) 
		& \leq (1+ 2 \rho^2 \log n) \hspace{0.1cm}. \tag{78}
\end{align*}
	Now, we want to bound $\mathbb{E}_{{\theta}_{0i}}\bigg( \sup_{\frac{1}{n} \leq \tau \leq 5t_n} Var(\theta_i|X_i,\tau) 1_{\{|X_i|> \sqrt{2 \rho^2 \log n}\}}  \bigg)$.
	%Note that, for any fixed $x \in \mathbb{R}$, $x^2E(\kappa^2|x,\tau)=x^2 E(\frac{1}{(1+\lambda^2\tau^2)^2}|x,\tau)$ is non-increasing in $\tau$. 
 Again we use the fact that for any fixed $x_i \in \mathbb{R}$ and $\tau>0$,
	$Var(\theta_i|x_i,\tau) \leq 1+x^2_i$. Using this fact and the monotonicity of  $E(\kappa^2_i|x_i,{\tau})$ for any fixed $x_i$, we have as in 
\textbf{Step-1} of Theorem 1, for that any $\tau \in [\frac{1}{n},5t_n]$
	\begin{align*}
		Var(\theta_i|x_i,{\tau}) &\leq 1+x^2_iE(\kappa^2_i|x_i,{\tau}) \\
		& \leq 1+x^2_iE(\kappa^2_i|x_i,\frac{1}{n}) \hspace{0.1cm}   \\
		& = \tilde{h}(x_i,\frac{1}{n}) \hspace{0.1cm}.
	\end{align*}
	Hence, using the same arguments used in \textbf{Step-1} of Theorem 1,
	
	\begin{equation*} 
		\sup_{|x_i|>\sqrt{2 \rho^2 \log n}} Var(\theta_i|x_i,{\tau}) \leq \sup_{|x_i|>\sqrt{2 \rho^2 \log n}} \tilde{h}(x_i,\frac{1}{n}) 
	\end{equation*}
 \begin{equation*}
     \leq 1+\sup_{|x_i|>\sqrt{2 \rho^2 \log n}} \tilde{h}_1(x_i,\frac{1}{n})+\sup_{|x_i|>\sqrt{2\rho^2 \log n}} \tilde{h}_2(x_i,\frac{1}{n}) \lesssim 1.
 \end{equation*}
	Using the above arguments, 
	\begin{align*}\label{eq:T-2.9}
		 & \mathbb{E}_{{\theta}_{0i}}\bigg( \sup_{\frac{1}{n} \leq \tau \leq 5t_n} Var(\theta_i|X_i,\tau) 1_{\{|X_i| > r_n\}}   \bigg) \lesssim 1\hspace{0.1cm} . \tag{79}
	\end{align*}
	Combining \eqref{eq:T-2.8} and \eqref{eq:T-2.9}, we obtain,
	\begin{equation*}
		\mathbb{E}_{{\theta}_{0i}}\bigg( \sup_{\frac{1}{n} \leq \tau \leq 5t_n} Var(\theta_i|X_i,\tau)  \bigg) \lesssim \log n \hspace{0.1cm} .
	\end{equation*}
	Since all of the above arguments hold uniformly for any $i$ such that ${\theta}_{0i} \neq 0$, as $n \to \infty$,
	\begin{equation*}\label{eq:T-2.10}
		\sum_{i:\theta_{0i} \neq 0} \mathbb{E}_{{\theta}_{0i}}\bigg( \sup_{\frac{1}{n} \leq \tau \leq 5t_n} Var(\theta_i|X_i,\tau)   \bigg) \lesssim \tilde{q}_n \log n  \hspace{0.1cm} .\tag{80}
	\end{equation*}
 \textbf{Step-2} Fix any $i$ such that $\theta_{0i} = 0$. Define $s_n= \sqrt{2 \log(\frac{1}{t_n})}$. We again decompose $\mathbb{E}_{{\theta}_{0i}} \bigg[\sup_{\frac{1}{n} \leq \tau \leq 5t_n} Var(\theta_i|X_i,\tau)  \bigg]$ as
	\begin{align*} 
		\mathbb{E}_{{\theta}_{0i}} \bigg[\sup_{\frac{1}{n} \leq \tau \leq 5t_n} Var(\theta_i|X_i,\tau)  \bigg] &= \mathbb{E}_{{\theta}_{0i}} \bigg[\sup_{\frac{1}{n} \leq \tau \leq 5t_n} Var(\theta_i|X_i,\tau) 1_{
			\{ |X_i| \leq s_n \}} \bigg] 
	\end{align*}
 \begin{equation*} \label{eq:T-4.14}
     + \mathbb{E}_{{\theta}_{0i}} \bigg[\sup_{\frac{1}{n} \leq \tau \leq 5t_n} Var(\theta_i|X_i,\tau) 1_{
			\{ |X_i| > s_n\}} \bigg] \hspace{0.1cm}. \tag{81}
 \end{equation*}
	Since for any fixed $x_i \in \mathbb{R}$ and $\tau>0$, $Var(\theta_i|x_i,\tau) \leq E(1-\kappa_i|x_i,\tau)+J(x_i,\tau)$ and $E(1-\kappa_i|x_i,\tau)$ is non-decreasing in $\tau$ and using Lemma 2 and Lemma A.2 of \cite{ghosh2017asymptotic}, for any fixed $x_i \in \mathbb{R}$ and $\tau>0$,
	\begin{align*} \label{eq:T-4.15}
		\mathbb{E}_{{\theta}_{0i}} \bigg[\sup_{\frac{1}{n} \leq \tau \leq 5t_n} Var(\theta_i|X_i,\tau) 1_{
			\{ |X_i| \leq \sqrt{2 \log(\frac{1}{t_n})} \}} \bigg] & \leq {t}_n \int_{0}^{\sqrt{2 \log (\frac{1}{t_n})}} e^{\frac{x^2}{2}} \phi(x) dx\\
		& \lesssim \frac{q_n}{n} \log(\frac{n}{q_n}). \tag{82}
	\end{align*}
	 Using the fact that for any $\tau>0$, $Var(\theta_i|x_i,\tau) \leq 1+x^2_i$ and the identity $x^2 \phi(x)=\phi(x)-\frac{d}{dx}[x \phi(x)]$, we obtain, 
	\begin{equation*} \label{eq:T-4.16}
		\mathbb{E}_{{\theta}_{0i}} \bigg[ \sup_{\frac{1}{n} \leq \tau \leq 5t_n} Var(\theta_i|X_i,\tau) 1_{
			\{ |X_i| > \sqrt{2 \log(\frac{1}{t_n})} \}} \bigg]   \lesssim \frac{q_n}{n} \log(\frac{n}{q_n}) \hspace{0.1cm}. \tag{83}
	\end{equation*}
	On combining \eqref{eq:T-4.14}-\eqref{eq:T-4.16}, for sufficiently large $n$,
	\begin{equation*} \label{eq:T-4.17}
		\mathbb{E}_{{\theta}_{0i}} \bigg[\sup_{\frac{1}{n} \leq \tau \leq 5t_n} Var(\theta_i|X_i,\tau)   \bigg]  \lesssim \frac{q_n}{n} \log(\frac{n}{q_n}) \hspace{0.1cm}. \tag{84}
	\end{equation*}
 Note that all these preceding arguments hold uniformly in $i$ such that $\theta_{0i} = 0$. Hence, as $n \to \infty$,
	\begin{align*} \label{eq:T-4.22}
		\sum_{i:\theta_{0i} = 0} \mathbb{E}_{{\theta}_{0i}}\bigg( \sup_{\frac{1}{n} \leq \tau \leq 5t_n} Var(\theta_i|X_i,\tau)   \bigg)  & \lesssim (n-\tilde{q}_n)\frac{q_n}{n} \log(\frac{n}{q_n}) 
		 \lesssim q_n \log(\frac{n}{q_n})\hspace{0.1cm}. \tag{85}
	\end{align*}
	The second inequality follows due to the fact that $\tilde{q}_n \leq q_n$ and $q_n=o(n)$ as $n \to \infty$.
 This completes the proof for \eqref{eq:T-2.4}. \\
 Moving towards proving \eqref{eq:T-2.3}, here, we need to prove that
 \begin{align*} \label{eq:T-2.11}
     \sum_{i:\theta_{0i} \neq 0} \mathbb{E}_{{\theta}_{0i}} \bigg[\sup_{\frac{1}{n} \leq \tau \leq 5t_n} (T_{\tau}(X_i)-\theta_{0i} )^2   \bigg] & \lesssim q_n \log n \tag{86}
 \end{align*}
 and
 \begin{align*} \label{eq:T-2.12}
     \sum_{i:\theta_{0i} = 0} \mathbb{E}_{{\theta}_{0i}} \bigg[\sup_{\frac{1}{n} \leq \tau \leq 5t_n} (T_{\tau}(X_i)-\theta_{0i} )^2   \bigg] & \lesssim q_n \log n \tag{87}.
 \end{align*}
 \textbf{Step-1} Fix any $i$ such that $\theta_{0i} \neq 0$. First, note that,
 \begin{align*}\label{eq:T-2.13}
      (T_{\tau}(X_i)-\theta_{0i} )^2 & \leq 2[ (T_{\tau}(X_i)-X_i )^2 + (X_i-\theta_{0i} )^2] \tag{88}.
 \end{align*}
 Observe that, the second term in the r.h.s. of \eqref{eq:T-2.13} is independent of choice of $\tau$. Also, using the fact, for any $x_i \in \mathbb{R}, |T_{\tau}(x_i)-x_i| \leq x_i$, we have
 \begin{align*}
     (T_{\tau}(X_i)-X_i )^2 1_{ \{ |X_i| \leq \rho \zeta_{\tau}\}} & \lesssim \zeta^2_{\tau}.
 \end{align*}
 Also, using Lemma 3 of \cite{ghosh2017asymptotic}, there exists a non-negative real-valued function $h$ such that $|T_{\tau}(x_i)-x_i| \leq h(x_i,\tau)$ for all $x_i \in \mathbb{R}$ satisfying for any $\rho >c>2$,
 \begin{align*} \label{eq:T-2.14}
     \lim_{\tau \to 0} \sup_{|x_i| >\rho \zeta_{\tau}} h(x_i,\tau) &=0.
 \end{align*}
 Combining these two arguments, we conclude, as $\tau \to 0$,
 \begin{align*}
     \sup_{\frac{1}{n} \leq \tau \leq 5t_n} (T_{\tau}(X_i)-X_i )^2  & \lesssim  \sup_{\frac{1}{n} \leq \tau \leq 5t_n} \zeta^2_{\tau}(1+o(1)) \lesssim \log n \tag{89}.
 \end{align*}
 Using \eqref{eq:T-2.13}, \eqref{eq:T-2.14} and noting that $\mathbb{E}_{{\theta}_{0i}}(X_i-\theta_{0i} )^2=1$, we obtain
 \begin{align*}
     \mathbb{E}_{{\theta}_{0i}} \bigg[\sup_{\frac{1}{n} \leq \tau \leq 5t_n} (T_{\tau}(X_i)-\theta_{0i} )^2   \bigg] & \lesssim  \log n.
 \end{align*}
 Since all of the above arguments hold uniformly for any $i$ such that ${\theta}_{0i} \neq 0$, as $n \to \infty$,
  \begin{align*} 
     \sum_{i:\theta_{0i} \neq 0} \mathbb{E}_{{\theta}_{0i}} \bigg[\sup_{\frac{1}{n} \leq \tau \leq 5t_n} (T_{\tau}(X_i)-\theta_{0i} )^2   \bigg] & \lesssim q_n \log n. 
 \end{align*}
 \textbf{Step-2} Fix any $i$ such that $\theta_{0i} = 0$. In this case,
 \begin{align*}
     \mathbb{E}_{{\theta}_{0i}}\bigg[\sup_{\frac{1}{n} \leq \tau \leq 5t_n} (T_{\tau}(X_i)-\theta_{0i} )^2 \bigg] & \leq \mathbb{E}_{{\theta}_{0i}}[T^2_{t_n}(X_i)] \lesssim \frac{q_n}{n}\log (\frac{n}{q_n}),
 \end{align*}
 where the inequality in the last line follows due to the use of an argument similar to that of \eqref{eq:T-1.10}.  Note that all these preceding arguments hold uniformly in $i$ such that $\theta_{0i} = 0$. Hence, as $n \to \infty$,
	\begin{align*} \label{eq:T-2.15}
		\sum_{i:\theta_{0i} = 0} \mathbb{E}_{{\theta}_{0i}}\bigg( \sup_{\frac{1}{n} \leq \tau \leq 5t_n}    (T_{\tau}(X_i)-\theta_{0i} )^2 \bigg)  & \lesssim (n-\tilde{q}_n)\frac{q_n}{n} \log(\frac{n}{q_n}) 
		 \lesssim q_n \log(\frac{n}{q_n})\hspace{0.1cm}. \tag{90}
	\end{align*}
	The second inequality follows due to the fact that $\tilde{q}_n \leq q_n$ and $q_n=o(n)$ as $n \to \infty$.
 This completes the proof for \eqref{eq:T-2.3}. As a result, Theorem 3 is also established using \eqref{eq:T-2.1}-\eqref{eq:T-2.3}.
% In the case of the second assertion, following the same steps mentioned before, we only need to show that, 
 %   \begin{align*}
% \mathbb{E}_{\boldsymbol{\theta}_0}\bigg[\sup_{\frac{1}{n} \leq \tau \leq 5t_n} \sum_{i=1}^{n} Var(\theta_i|\mathbf{X},\tau)   \bigg] & \lesssim q_n \log n, 
%    \end{align*}
 %   which is already proved in the first assertion.
 \end{proof}

 { %Though we have not provided any upper bound on the second term of $E(1-\kappa_i|x_i,\tau)$ in Lemma \ref{lem1}, later either we will use \eqref{L-1.1} directly or an \eqref{L-1.3} and an upper bound on $A_2$ based on the situations as required. This Lemma will be used in \textbf{Case-2} in the proof of Theorems 1 and \ref{Thm-3} when $a \geq 1$. Another important point that we need to mention here is that the upper bound on the posterior shrinkage coefficient is independent of any assumption on $L(\cdot)$, but assumptions (\hyperlink{A1}{A1})and (\hyperlink{A2}{A2}) will definitely be used in proving the subsequent results.}
\subsection{Results related to the ABOS of the full Bayes procedure}

\begin{proof}[Proof of Theorem \ref{Thm-6}]
	First, note that
	\begin{align*} \label{eq:T-6.1}
		E(1-\kappa_i|\mathbf{X}) &= \int_{\frac{1}{n}}^{\alpha_n} E(1-\kappa_i|\mathbf{X},\tau) \pi_{2n}(\tau|\mathbf{X}) d \tau \\
		&= \int_{\frac{1}{n}}^{\alpha_n} E(1-\kappa_i|{X}_i,\tau) \pi_{2n}(\tau|\mathbf{X}) d \tau  \leq E(1-\kappa_i|{X}_i,\alpha_n) \hspace*{0.05cm} \cdot \tag{91}
	\end{align*}
	For proving the inequality above we first use the fact that given any $x_i \in \mathbb{R}$, $E(1-\kappa_i|x_i,\tau)$ is non-decreasing in $\tau$ and next we use (\hyperlink{C4}{C4}). Now, using Theorem 4 of \cite{ghosh2016asymptotic}, we have for $a \in (0,1)$, as $n \to \infty$
	\begin{equation*} \label{eq:T-6.2}
		E(1-\kappa_i|{X}_i,\alpha_n) \leq \frac{ KM}{a(1-a)} e^{\frac{X^2_i}{2}} {\alpha^{2a}_n} (1+o(1)). \tag{92}
	\end{equation*}
	Here $o(1)$ depends only on $n$ such that $\lim_{n \to \infty}o(1)=0$ and is independent of $i$. Two constants $K$ and $M$ are defined in \eqref{eq:1.4} and \hyperlink{assumption1}{Assumption 1}, respectively.
%	\hspace*{0.5cm} From \eqref{eq:2.2} it is evident  that  for any $\tau >0$, the posterior distribution of $\kappa_i$ given $X_i$ and $\tau$ is the same for each $i, i=1,2,\cdots,n$. Also note that under $H_{0i}$, $X_i \simiid \mathcal{N}(0,1)$. So, the distribution of $X_i$ is independent of the choice of $i$ and hence the same for all $i$. As a result of this, the probability of  type-I error induced by the decision rule \eqref{eq:4.9}, $t^{\text{FB}}_{1i}= P_{H_{0i}}(E(1-\kappa_i|\mathbf{X})> \frac{1}{2})$ is same for each $i$ and denoted as $t^{\text{FB}}_{1}$. Using these facts along with \eqref{eq:T-6.1} and \eqref{eq:T-6.2}, 
Hence, we have the following :
	\begin{align*} 
		t^{\text{FB}}_{1i} &= P_{H_{0i}}(E(1-\kappa_i|\mathbf{X})> \frac{1}{2}) \\
		& \leq P_{H_{0i}}\bigg(\frac{X^2_i}{2} > 2a \log (\frac{1}{\alpha_n})+ \log \big(\frac{a(1-a)}{2 A_0 KM} \big)- \log (1+o(1)) \bigg) \\
		&=  P_{H_{0i}} \bigg(|X_i| > \sqrt{4a \log (\frac{1}{\alpha_n})}\bigg)(1+o(1)) \hspace*{0.05cm}.
	\end{align*}
% Hence, we have the following :
	%\begin{align*} 
	%	t^{\text{FB}}_{1i} &= P_{H_{0i}}(E(1-\kappa_i|\mathbf{X})> \frac{1}{2}) \\
	%	& \leq P_{H_{0i}}(\frac{A_0 KM}{a(1-a)} e^{\frac{X^2_i}{2}} {\alpha^{2a}_n} (1+o(1)) > \frac{1}{2}) \\
	%	&= P_{H_{0i}}\bigg(\frac{X^2_i}{2} > 2a \log (\frac{1}{\alpha_n})+ \log \big(\frac{a(1-a)}{2 A_0 KM} \big)- \log (1+o(1)) \bigg) \\
	%	&=  P_{H_{0i}} \bigg(|X_i| > \sqrt{4a \log (\frac{1}{\alpha_n})}\bigg)(1+o(1)) \hspace*{0.05cm}.
	%\end{align*}
 Note that, as $\alpha_n <1$ for all $n \geq 1$, $4a \log (\frac{1}{\alpha_n}) >0$ for all $n \geq 1$. Next, using the fact that under $H_{0i}, X_i \simiid \mathcal{N}(0,1)$ with $1-\Phi(t) < \frac{\phi(t)}{t}$ for $t >0$, we get
 \begin{align*} \label{eq:T-6.3}
     t^{\text{FB}}_{1i} &  \leq 2 \frac{\phi(\sqrt{ 4 a \log (\frac{1}{\alpha_n})})}{\sqrt{4a \log (\frac{1}{\alpha_n})}}(1+o(1)) = \frac{1}{\sqrt{\pi a}} \frac{{\alpha}^{2a}_n}{\sqrt{\log (\frac{1}{\alpha^2_n}) }} (1+o(1)) \hspace*{0.05cm}. \tag{93}
 \end{align*}
Here $o(1)$ depends only on $n$ such that $\lim_{n \to \infty} o(1)=0$. This completes the proof of Theorem 4.
\end{proof}

\begin{proof}[Proof of Theorem \ref{Thm-7}]

%	Here we use exactly similar argument as used in Theorem \ref{Thm-6} about the equality of $t^{\text{FB}}_{1i}$ for $i=1,2,\cdots,n$. In this case, we note that, under $H_{1i}$, $X_i \sim \mathcal{N}(0,1+\psi^2)$, i.e. the distribution of $X_i$ is independent of the choice of $i$ and hence same for all $i$. As a result of this, 
 
 In order to provide an upper bound on the probability of  type-II error induced by the decision rule (4.9) of the main paper, $t^{\text{FB}}_{2i}= P_{H_{1i}}(E(\kappa_i|\mathbf{X})> \frac{1}{2})$, 
 we first note that,
	\begin{align*} \label{eq:T-7.4}
		E(\kappa_i|\mathbf{X}) &= \int_{\frac{1}{n}}^{\alpha_n} E(\kappa_i|\mathbf{X},\tau) \pi_{2n}(\tau|\mathbf{X}) d \tau \\
		&= \int_{\frac{1}{n}}^{\alpha_n} E(\kappa_i|{X}_i,\tau) \pi_{2n}(\tau|\mathbf{X}) d \tau 
		 \leq E(\kappa_i|{X}_i,\frac{1}{n}), \tag{94}
	\end{align*} 
	where inequality in the last line follows due to the fact that given any $x_i \in \mathbb{R}$, 
 $E(\kappa_i|x_i,\tau)$ is non-increasing in $\tau$. Using arguments similar to Lemma 3 of \cite{ghosh2017asymptotic}, for each fixed $\eta \in (0,1)$ and $\delta \in (0,1)$, $E(\kappa_i|x_i,\tau)$ can be bounded above by a real-valued function $g(x_i,\tau, \eta, \delta)$, depending on $\eta$ and $\delta$, is given by,
	$g(x_i,\tau,\eta, \delta)= g_1(x_i, \tau)+g_2(x_i, \tau,\eta, \delta)$ where 
 $g_1(x_i,\tau)=C(a,\eta,L)[x_i^2\int_{0}^{ \frac{x_i^2}{1+t_0} } e^{-\frac{u}{2}} u^{a-\frac{1}{2}} du]^{-1}$ and \\ $g_2(x_i,\tau,\eta, \delta )= \tau^{-2a} c_0 H(a,\eta,\delta) e^{-\frac{\eta(1-\delta) x_i^2}{2}}$ where $C(a,\eta,L)$ and $ H(a,\eta,\delta)$ are two global constants independent of any $i$ and $n$ and $c_0$ is defined in \hyperlink{assumption1}{Assumption 1}. 
 It may be noted that for $\rho > \frac{2}{\eta (1-\delta)}$,
	\begin{align*} \label{eq:T-7.5}
		\lim_{n \to \infty} \sup_{|x_i| > \sqrt{\rho \log(n^{2a })}} g(x_i, \frac{1}{n},\eta,\delta) &=0 \hspace*{0.05cm}. \tag{95}
	\end{align*}
 It is evident that there exists $N$ (depending on $\eta, \delta$) such that for all $n \geq N$, 
\begin{equation*}
     \sup_{|x_i| > \sqrt{ \rho \log (n^{2a})}} g(x_i,\frac{1}{n},\eta,\delta) \leq \frac{1}{2}.
\end{equation*}
	Now define $B_n$ and $C_n$ as
	$B_n =\{g(X_i, \frac{1}{n},\eta, \delta) >\frac{1}{2} \}$ and $C_n= \{ |X_i| > \sqrt{\rho \log (n^{2a})} \}$ . Hence
	\begin{align*} \label{eq:T-7.6}
		t^{\text{FB}}_{2i} &= P_{H_{1i}}(E(\kappa_i|\mathbf{X}) > \frac{1}{2}) \\
		& \leq P_{H_{1i}}(B_n)   
		 \leq P_{H_{1i}}(B_n\cap C_n)+P_{H_{1i}}(C^{c}_n), \tag{96}
	\end{align*}
 \text{for any fixed} $(\eta,\delta) \in (0,1)\times (0,1)$ \text{and} \hspace*{0.05cm} $\rho >\frac{2}{\eta(1-\delta)} $.
	Recall that, from the definition of $g_1(\cdot)$ and $g_2(\cdot)$, it is clear that $g(\cdot)$ is monotonically decreasing in $|x_i|$ for $x_i \neq 0$. 
	Next using the definition of $B_n$ and $C_n$ and with the use of \eqref{eq:T-7.5}, we  immediately have,
	\begin{align*} \label{eq:T-7.7}
		  P_{H_{1i}}(B_n\cap C_n)=0\tag{97}
	\end{align*}
 for all sufficiently large $n$, depending on $\rho, \eta, \delta$.
	Further, note that under $H_{1i}, X_i \simiid \mathcal{N}(0,1+\psi^2)$ and hence
	\begin{align*} \label{eq:T-7.8}
		P_{H_{1i}}(C^{c}_n) &= P_{H_{1i}}(|X_i| \leq \sqrt{\rho \log (n^{2a})}) \\
		&= P_{H_{1i}} \bigg( \frac{|X_i|}{\sqrt{1+\psi^2}} \leq \sqrt{a \rho} \sqrt{ \frac{2 \log n}{1+\psi^2}} \bigg) = \bigg[2 \Phi \bigg(\sqrt{\frac{a \rho C}{\epsilon}} \bigg)-1\bigg] (1+o(1)), \tag{98}
	\end{align*}
 where the last equality follows from \hyperlink{assumption2}{Assumption 2}. 
 Here the term $o(1)$ is independent of $i$ and $\lim_{n \to \infty} o(1)=0$. 
   From \eqref{eq:T-7.6}-\eqref{eq:T-7.8} the desired result follows.	
\end{proof}

\begin{proof}[Proof of Theorem \ref{Thm-8}]
	The Bayes risk corresponding to the decision rule 
 (4.9) of the main document,  denoted $R^{\text{FB}}_{\text{OG}}$ is of the form
 \begin{align*} \label{eq:T-8.1}
		R^{\text{FB}}_{\text{OG}} 
  &=\sum_{i=1}^{n} [(1-p)t^{\text{FB}}_{1i}+pt^{\text{FB}}_{2i}] 
  = p\sum_{i=1}^{n} \bigg[\bigg(\frac{1-p}{p}\bigg)t^{\text{FB}}_{1i}+t^{\text{FB}}_{2i} \bigg] \hspace{0.05cm}.  \tag{99}
	\end{align*}
	%since as argued before, $t^{\text{FB}}_{1i}=t^{\text{FB}}_{1}$ and $t^{\text{FB}}_{2i}=t^{\text{FB}}_{2}$ for $i=1,2,\cdots,n$. 
Now, for $a \in [0.5,1)$, using \eqref{eq:T-6.3} we have
	\begin{equation*}
		t^{\text{FB}}_{1i} \lesssim \frac{\alpha_n}{\sqrt{\log (\frac{1}{\alpha_n})}} (1+o(1)) \hspace*{0.05cm}.
	\end{equation*}
	Hence, for $p_n \propto n^{-\epsilon}, 0< \epsilon \leq 1$, 
 with the choice of $\alpha_n$ as $\log (\frac{1}{\alpha_n})=\log n-\frac{1}{2} \log \log n+c_n$, where $c_n =o(\log \log n)$ as $n \to \infty$,
 we have, 	
 \begin{align*}\label{eq:T-8.2}
     \frac{1}{n} \sum_{i=1}^{n} \bigg(\frac{1-p}{p}\bigg)t^{\text{FB}}_{1i}
     & \lesssim \frac{n^{\epsilon}\alpha_n}{\sqrt{\log (\frac{1}{\alpha_n})}} (1+o(1)) =o(1) \hspace*{0.05cm}.\tag{100}
 \end{align*}
Next, note that, from Theorem 5, for each $i$, the upper bound of $t^{\text{FB}}_{2i}$ is independent of $i$. Therefore, we obtain
\begin{align*}\label{eq:T-8.3}
    \sum_{i=1}^{n}t^{\text{FB}}_{2i} \leq n \bigg[2 \Phi \bigg(\sqrt{ \frac{a \rho C}{\epsilon}} \bigg)-1 \bigg](1+o(1))  \hspace*{0.05cm}.\tag{101}
\end{align*}
 As a result, using \eqref{eq:T-8.1}-\eqref{eq:T-8.3}, for any fixed $\eta \in (0,1)$, $\delta \in (0,1)$ and $\rho >\frac{2}{\eta(1-\delta)}$ and any $\epsilon \in (0,1]$
 \begin{equation*} \label{eq:T-8.4}
     R^{\text{FB}}_{\text{OG}} \leq np \bigg[2 \Phi \bigg(\sqrt{ \frac{a \rho C}{\epsilon}} \bigg)-1 \bigg] (1+o(1))  \hspace*{0.05cm}. \tag{102}
 \end{equation*}
  Since, by definition, $\frac{ R^{\text{FB}}_{\text{OG}}}{ R^{\text{BO}}_{\text{Opt}}} \geq 1$, using (4.6) of the main paper, for any fixed  $\eta \in (0,1)$, $\delta \in (0,1)$, $\rho > \frac{2}{\eta (1-\delta)}$, and any $\epsilon \in (0,1]$, 
\begin{equation*} \label{eq:T-8.5}
		1\leq \frac{ R^{\text{FB}}_{\text{OG}}}{ R^{\text{BO}}_{\text{Opt}}} \leq \frac{\bigg[2 \Phi \bigg(\sqrt{ \frac{a \rho C}{\epsilon}} \bigg)-1 \bigg]}{\bigg[2 \Phi(\sqrt{{C}})-1 \bigg]} (1+o(1)) ,\tag{103}
	\end{equation*}
	%\begin{equation*} \label{eq:T-8.3}
		%\frac{\bigg[2 \Phi(\sqrt{2a{C}})-1 \bigg]}{\bigg[2 \Phi(\sqrt{{C}})-1 \bigg]}(1+o(1)) \leq \frac{ R^{\text{FB}}_{\text{OG}}}{ R^{\text{BO}}_{\text{Opt}}} \leq \frac{\bigg[2 \Phi \bigg(\sqrt{ \frac{a \rho C}{\epsilon}} \bigg)-1 \bigg]}{\bigg[2 \Phi(\sqrt{{C}})-1 \bigg]} (1+o(1)) \tag{T-8.3}
	%\end{equation*}
where $o(1)$ term is independent of $i$ and depends on $n, \eta, \delta$ such that $\lim_{n \to \infty}o(1)=0$.
 For extremely sparse situation, i.e.  for $\epsilon=1$, taking limit inferior and limit superior in \eqref{eq:T-8.5}, we have,

\begin{equation*} \label{eq:T-8.6}
		1 \leq \liminf_{n \to \infty} \frac{ R^{\text{FB}}_{\text{OG}}}{ R^{\text{BO}}_{\text{Opt}}} \leq \limsup_{n \to \infty} \frac{ R^{\text{FB}}_{\text{OG}}}{ R^{\text{BO}}_{\text{Opt}}} \leq \frac{\bigg[2 \Phi \bigg(\sqrt{ {a \rho C}} \bigg)-1 \bigg]}{\bigg[2 \Phi(\sqrt{{C}})-1 \bigg]}. \tag{104}
	\end{equation*}
	%\begin{equation*} \label{eq:T-8.4}
	%	\frac{\bigg[2 \Phi(\sqrt{2a{C}})-1 \bigg]}{\bigg[2 \Phi(\sqrt{{C}})-1 \bigg]} \leq \liminf_{n \to \infty} \frac{ R^{\text{FB}}_{\text{OG}}}{ R^{\text{BO}}_{\text{Opt}}} \leq \limsup_{n \to \infty} \frac{ R^{\text{FB}}_{\text{OG}}}{ R^{\text{BO}}_{\text{Opt}}} \leq \frac{\bigg[2 \Phi \bigg(\sqrt{ {a \rho C}} \bigg)-1 \bigg]}{\bigg[2 \Phi(\sqrt{{C}})-1 \bigg]}. \tag{T-8.4}
	%\end{equation*}
Note that the decision rule does not depend on how $\eta \in (0,1),\delta \in (0,1)$ and $\rho >\frac{2}{\eta (1-\delta)}$ are chosen. Hence, the ratio $\frac{ R^{\text{FB}}_{\text{OG}}}{ R^{\text{BO}}_{\text{Opt}}}$ is also independent to the choices of $\eta, \delta$ and $\rho$ for any $n \geq 1$. As a result, limit inferior and limit superior in \eqref{eq:T-8.6} is also free of any particular value of 
	$\eta, \delta$ and $\rho$. Taking infimum of the r.h.s. of \eqref{eq:T-8.6} over all possible choices of $(\eta,\delta)\in (0,1) \times (0,1)$ and $\rho >\frac{2}{\eta(1-\delta)}$ and finally using continuity of $\Phi(\cdot)$, we can conclude that, for extremely sparse case with $a=0.5$,
	\begin{equation*}
		\lim_{n \to \infty}	\frac{ R^{\text{FB}}_{\text{OG}}}{ R^{\text{BO}}_{\text{Opt}}}=1.
	\end{equation*}\end{proof}
}

%% ** The bibliograhy **
 %place <bibliography.bib> in ./bib folder 
\end{document}